\documentclass[12pt,a4paper,english,reqno]{amsart}

\usepackage{amsmath,amsfonts,amssymb,mathrsfs}
\usepackage[english]{babel}
\usepackage[utf8x]{inputenc}
\usepackage[T1]{fontenc}
\usepackage{graphicx}
\usepackage{tabularx}
\usepackage{hyperref}
\usepackage{mathtools}
\usepackage[margin=1in]{geometry}
\usepackage{cleveref}
\usepackage{comment}
\usepackage[all]{xy}
\usepackage[usenames,dvipsnames]{color}
\usepackage{tikz}
\usepackage{hyperref}

\newcommand{\nocontentsline}[3]{}
\newcommand{\tocless}[2]{\bgroup\let\addcontentsline=\nocontentsline#1{#2}\egroup}
\DeclareUnicodeCharacter{0177}{\^u}

\DeclareMathOperator{\rank}{rank}
\DeclareMathOperator{\End}{End}
\DeclareMathOperator{\Hom}{Hom}

\DeclareMathOperator{\SL}{SL}
\DeclareMathOperator{\GL}{GL}
\DeclareMathOperator{\SO}{SO}

\DeclareMathOperator{\Aut}{Aut}

\DeclareMathOperator{\Vol}{Vol}
\DeclareMathOperator{\tr}{tr}

\DeclareMathOperator{\disc}{disc}

\DeclareMathOperator{\Ad}{Ad}

\DeclareMathOperator{\Hodge}{Hdg}

\numberwithin{equation}{subsection}

\theoremstyle{plain}

\newtheorem{theorem}{Theorem}[section]

\newtheorem{proposition}[theorem]{Proposition}
\newtheorem{lemme}[theorem]{Lemma}

\newtheorem{lemma}[theorem]{Lemma}

\newtheorem{corollaire}[theorem]{Corollary}
\newtheorem{conjecture}[theorem]{Conjecture}
\newtheorem{question-nr}[theorem]{Question}
\newtheorem*{question}{Question}
\theoremstyle{definition}
\newtheorem{definition}[theorem]{Definition}

\theoremstyle{remark}
\newtheorem{remarque}[theorem]{Remark}
\newtheorem{example}[theorem]{Example}

\newcommand{\R}{\mathbb{R}}

\newcommand{\Z}{\mathbb{Z}}
\newcommand{\Q}{\mathbb{Q}}

\newcommand{\N}{\mathbb{N}}

\newcommand{\V}{\mathbb{V}}

\newcommand{\C}{\mathbb{C}}

\newcommand{\D}{\mathcal{D}}

\newcommand{\cA}{\mathcal{A}}

\newcommand{\cS}{\mathcal{S}}
\newcommand{\cF}{\mathcal{F}}

\newcommand{\cV}{\mathcal{V}}
\newcommand{\cU}{\mathcal U}

\renewcommand{\H}{\mathrm H}

\newcommand{\fg}{\mathfrak{g}}
\newcommand{\fh}{\mathfrak{h}}
\newcommand{\fk}{\mathfrak{k}}
\newcommand{\fl}{\mathfrak{l}}

\newcommand{\fp}{\mathfrak{p}}

\newcommand{\so}{\mathfrak{so}}

\newcommand{\p}{\overline{p}}

\renewcommand{\d}{\mathrm{d}}

\newcommand{\tuple}[1]{\underline{#1}}

\renewcommand{\hat}[1]{\widehat{#1}}
\renewcommand{\tilde}[1]{\widetilde{#1}}

\newcommand{\equaldef}{\overset{\textrm{def}}{=}}

\title[]{Equidistribution of Hodge loci II}
\author{Salim Tayou}
\address{Institute for Advanced Study, 1 Einstein drive, Princeton, New  Jersey 08540, USA}
\email{tayou@math.harvard.edu}

\author{Nicolas Tholozan}
\address{DMA – UMR8553, École Normale Supérieure, CNRS – PSL Research University, 45 rue d’Ulm, 75230 Paris Cedex 5, France}
\email{nicolas.tholozan@ens.fr}

\date\today

\begin{document}

\begin{abstract}
Let $\mathbb V$ be a polarized variation of Hodge structure over a smooth complex quasi-projective variety $S$. In this paper, we give a complete description of the typical Hodge locus for such variations. We prove that it is either empty or equidistributed with respect to a natural differential form, \emph{the pull-push form}. In particular, it is always analytically dense when the pull-push form does not vanish. When the weight is $2$, the Hodge numbers are $(q,p,q)$ and the dimension of $S$ is least $rq$, we prove that the typical locus where the Picard rank is at least $r$ is equidistributed in $S$ with respect to the volume form $c_q^r$, where $c_q$ is the $q$\textsuperscript{th} Chern form of the Hodge bundle. We obtain also several equidistribution results of the typical locus in Shimura varieties: a criterion  for the density of the typical Hodge loci of a variety in $\mathcal{A}_g$, equidistribution of certain families of CM points and equidistribution of Hecke translates of curves and surfaces in $\mathcal A_g$. 

These results are proved in the much broader context of dynamics on homogeneous spaces of Lie groups which are of independent interest. The pull-push form appear in this greater generality and we provide several tools to determine it and we compute it in many examples.

\end{abstract}
\thanks{}
\maketitle
\setcounter{tocdepth}{1}

\tableofcontents

\section{Introduction}

Let $G$ be a semi-simple Lie group and let $\Gamma\subset G$ be a lattice. Homogeneous dynamics is traditionally interested in the equidistribution properties of the orbits of a Lie subgroup $H$ acting on $\Gamma \backslash G$ by right multiplication and, dually, on the dynamics of the left action of $\Gamma$ on $G/H$. Classifying the closure of orbits of such actions is the subject of an extensive literature with far reaching applications to number theory and ergodic theory.

In this paper, our first purpose is to provide a fairly general answer to the following question:
\begin{question}
Assume that a sequence of closed $H$-orbits is equidistributed in $\Gamma \backslash G$. Can we deduce an equidistribution result for the intersection of these $H$-orbits with a fixed analytic subvariety $V\subset \Gamma \backslash G$?
\end{question}

We consider more generally the following setting: let $G$ be a real semi-simple Lie group, $\Gamma$ a lattice in $G$, $H$ a semi-simple subgroup of $G$, $K$ a compact subgroup of $G$ (which is not assumed to be maximal), and $L= H\cap K$. Denote by $p$ the projection map $G/L\rightarrow G/H$ and by $\pi$ the projection map $G/L\rightarrow G/K$. We fix compatible choices of invariant volume forms $\omega_{G}$, $\omega_{G/H}$ and $\omega_{H}$ on $G$, $G/H$ and $H$ respectively (see \Cref{notations}). 

Let $(\mathcal O_n)_{n\in \N}$ be a sequence of finite unions of closed $H$-orbits in $\Gamma \backslash G$. Since $H$ is semi-simple, $\mathcal O_n$ has finite volume with respect to the volume form $\omega_H$ along $H$-orbits for every $n\in \N$ . We say that the sequence $(\mathcal O_n)_{n\in \N}$ is equidistributed in $\Gamma \backslash G$ if the normalized integration measure\footnote{See \Cref{normalization-convention} for our convention on the normalization of this measure.}  on $(\mathcal O_n,\omega_H)$ converges weakly to the Haar measure $\omega_G$ on $\Gamma \backslash G$.

Let now $S$ be a real analytic subvariety of $\Gamma \backslash G/K$ of codimension $\dim(H/L)$ whose smooth locus is oriented. Denote by $\mu_{S\cap \pi(\mathcal O_n)}$ the \emph{transverse intersection measure} of $V$ and $\mathcal O_n$, which counts (with an orientation sign and multiplicity) the transverse intersection points between $S$ and $\pi(\mathcal O_n)$ (see \Cref{countingpoints} for precisions). We prove the following:

\begin{theorem} \label{main}
Assume that the sequence $(\mathcal O_n)_{n\in \N}$ is equidistributed in $\Gamma \backslash G$. Then the sequence of signed measures $\frac{1}{\Vol(\mathcal O_n)}\mu_{S\cap \pi(\mathcal O_n)}$ on $S$ converges weakly to the restriction of the $G$-invariant form
\[\frac{1}{\Vol(\Gamma \backslash G/L)}\pi_*p^*\omega_{G/H}~.\]
\end{theorem}

This general result has countless potential applications, some of which will be detailed in the paper. What is interesting about this theorem is that the \emph{pull-push form} $\pi_*p^*\omega_{G/H}$ is not easy, in general, to determine precisely and depends greatly on the subgroup $H$. These forms were studied extensively by the second author in \cite{tholozan2016volume} with a very different motivation. Building on this previous work, we will give tools to analyze this pull-push form and characterize it in various examples.\\

It sometimes happens that the form $\pi_*p^*\omega_{G/H}$ vanishes (this vanishing played an important role in \cite{tholozan2016volume}). In that case, our theorem only asserts that positive and negative intersection points ``cancel each other'' asymptotically.

The theorem is stronger when $G/K$, $H/L$ and $S$ are complex analytic. Then, all intersections are counted positively and the form $\pi_*p^*\omega_{G/H}$ does not vanish (see \Cref{complexnonvanishing}). It vanishes in restriction to $S$ if and only if all the intersections $S\cap \mathcal O_n$ have ``exceptional dimension'' (i.e. complex dimension $\geq 1$). 

In the applications we develop in the next sections, $H/L\subset G/K$ will be Mumford--Tate domains of Hodge structures and $S\rightarrow \Gamma \backslash G/K$ is the period map of a polarized variation of Hodge structure. In \Cref{generictransversality}, we give an algebraic characterization of when the form $\pi_*p^*\omega_{G/H}$ is non-zero in restriction to some analytic variation of Hodge structure.

\subsection{Equidistribution of typical Hodge loci}

One of the main motivations of \Cref{main} is its application to the equidistribution of Hodge loci of variations of Hodge structure. In fact, the present paper is a continuation of the first author's previous work \cite{tayouequi}, which studied the particular case of the Noether--Lefschetz locus of a one-parameter family of K3 surfaces.

Let $\mathbb V=\{\mathbb{V}_\Z,\mathcal{F}^\bullet\mathcal{V},B\}$ be a polarized variation of Hodge structure ($\Z$-PVHS) of weight $2k$ over a complex quasi-projective algebraic variety $S$ of dimension $d\geq 1$ (see \Cref{definitionvhs}). The \emph{group of Hodge classes} $\Hodge(s)$ at a point $s$ is the free abelian group $\mathbb{V}_{\Z,s} \cap \mathcal F^k\cV_s$. An important class of examples of $\Z$-PVHS is provided by those of \emph{geometric origin}: starting from a smooth projective morphism $f:X\rightarrow S$ where $S$ is a smooth complex quasi-projective variety, the $2k^{th}$ cohomology groups of the fibers, modulo torsion, endowed with their Hodge structure give rise to a $\Z$-PVHS of weight $2k$ on $S$, see \cite[Partie III]{voisin} for more details. 

More generally, let $\mathbb{V}^{\otimes}$ be the countable direct sum of of $\Z$-PVHS $\bigoplus_{a,b\geq 0} \mathbb{V}^{a}\otimes\left(\mathbb{V}^{\vee}\right)^{b}$, where $\mathbb{V}^{\vee}$ is the dual $\Z$-PVHS. Then the Hodge locus $HL(S,\mathbb{V}^{\otimes})$ is defined as the subset $s\in S$ where $\mathbb{V}_s$ has more Hodge tensors than the very general fiber $\mathbb{V}_{s'}$. It is a countable union of algebraic subvarieties by  \cite{cattanidelignekaplan,bkt}.

In this paper, we give a precise description of the typical part of this Hodge locus. More precisely, let $G$ be the generic Mumford--Tate group of the variation. The Hodge locus can also be seen as the locus of points of $S$ where the Mumford--Tate group is strictly contained in $G$.

For $H\subset G$ a sub-Mumford--Tate group, the \emph{typical Hodge locus for $H$} is the set of points $s\in HL(S,\mathbb{V}^{\otimes})$ whose Mumford--Tate group is contained in $H$ and such that $\pi^*p_*\omega_{G/H,s}\neq 0$. The typical Hodge locus is then the union of typical Hodge loci over all Mumford--Tate subgroups $H$. 

We prove then the following theorem, see \Cref{generictransversality}. 
\begin{theorem}\label{T:typical}
The following statements are equivalent: 
\begin{enumerate}
    \item There exists $H\subseteq G$ such that the typical Hodge locus for $H$ is equidistributed with respect to $\pi_*p^* \omega_{G/H}$ and in particular analytically dense. 
    \item There exists $H\subseteq G$ and one point $s\in S$ such that $\pi_*p^* \omega_{G/H}$ is non-zero at $x$.
    \item The typical Hodge locus is non-empty. 
\end{enumerate}
\end{theorem}
\begin{remarque}
The equidistribution assertion in the above theorem as well as in all the subsequent ones in this article, when not explicitly specified, should be understood in the sense of \Cref{main}.  
\end{remarque}
\begin{remarque}
In the case of the Noether--Lefschetz loci, the above theorem is a strengthening of the classical criterion of Green, see \cite[Prop. 17.20]{voisin}. 
\end{remarque}

The following proposition gives a criterion for the emptiness of the typical Hodge locus, see \Cref{p:level2}.
\begin{proposition}\label{p:level2_intro}
If for every sub-Mumford--Tate group $H\subseteq G$ the Hodge structure $\mathfrak{g}/\mathfrak{h}$ satisfies: \[(\mathfrak{g}/\mathfrak{h})^{-p,p}\neq 0\quad \textrm{for some}\quad |p|\geq 2,\] then the typical Hodge locus is empty. 
\end{proposition}

\Cref{T:typical} will be applied to situations where we know how to compute the form $\pi_*p^* \omega_{G/H}$. These applications will be explained in the next section. 
\medskip

\Cref{T:typical} and \Cref{p:level2_intro} have also been independently studied by Baldi--Klingler--Ullmo in \cite{baldiklinglerullmo}, see also the prior work of Klingler--Otwinowska \cite{klinglerotwinowska}. Moreover, the authors prove in \cite[Theorem 2.3]{baldiklinglerullmo} that the condition in \Cref{p:level2_intro} is always satisfied whenever $\mathbb{V}$ has level more than $3$. 

\subsection{Applications}
We explain now further applications of our main theorem. They correspond to situations where we know how to compute the pull-push form and they are hence far from exhaustive. 
\subsubsection{Refined Noether--Lefschetz loci}\label{section:hodgeloci}

Let $\mathbb{V}$ be a polarized variation of Hodge structure of weight $2$ over a complex quasi-projective algebraic variety $S$ of dimension $d\geq 1$. Let $(q,p,q)$ be the Hodge numbers.

Without loss of generality, one can assume that the $\Z$-PVHS is \emph{simple}, so that the group of Hodge classes at a generic point is $0$. The \emph{Noether--Lefschetz locus} is then defined as the subset of elements of $S$ which admit non-trivial Hodge classes. More generally, we define the \emph{refined Noether--Lefschetz locus} of rank $r$ as the subset where the group of Hodge classes has rank at least $r$, and denote it by $NL^{\geq r}(S)$.

Let $(V_\Z,B)$ be the fiber of $\mathbb V$ at a point $s\in S$. The period domain associated to $\mathbb V$ is the homogeneous space $\mathcal D= G/K$, where $G$ is the real group $\SO(B)$ and $K$ the stabilizer of the Hodge structure at $s$, and $S$ has a \emph{period map} to the quotient $\Gamma \backslash \mathcal D$, where $\Gamma$ is the subgroup of $G$ preserving $V_\Z$. For our purposes, we will assume that $\mathbb V$ has generically immersive period map. Otherwise, we can replace $S$ by its image by the period map, which still has an algebraic structure by \cite{bakker-brunebarbe-tsimerman}. 

The refined Noether--Lefschetz locus of rank $r$ is a countable union of algebraic subvarieties which are the intersection of $S$ with certain Mumford--Tate subdomains obtained as projections of right $H$-orbits for the subgroup $H$ of $G$ stabilizing a set of $r$ integral elements in $V_\Z$ with positive intersection matrix. Applying \Cref{main} in this setting, we obtain equidistribution results for refined Noether--Lefschetz loci.

For all $n\in \N_{>0}$, define $NL^{\geq r}(n)$ as the set of points $s$ such that $(\Hodge(s),B)$ contains a primitive sublattice of rank $r$ in restriction to which $B$ has discriminant at most $n$. The set $ NL^{\geq r}(n)$ is an algebraic subvariety of $S$. It has expected dimension $d-rq$ but can contain higher dimensional components. 

\begin{theorem}\label{maink3-alt}
Let $\{\mathbb{V}_\Z,\mathcal{F}^\bullet\mathcal{V},B\}$ be a simple $\Z$-PVHS of weight $2$ and Hodge numbers $(q,p,q)$ over a complex analytic variety $S$ of dimension $d=rq$ and which has generically immersive period map. If $r\leq p$, then there is a constant $\lambda>0$ such that, for every relatively compact open subset $\Omega\subset S$ with boundary of measure $0$, we have
\[n^{-\frac{p+2q}{2}}\left \vert \{(s,P), s\in \Omega, P\subseteq \Hodge(s), \rank(P)=r, \mathrm{disc}(P)\leq n\} \right \vert \underset{n\to+\infty}{\longrightarrow} \lambda \int_\Omega c_q(\mathcal F^2\cV)^r~,\]
where $c_q$ denotes the $q$\textsuperscript{th} \emph{Chern form} of the bundle $\mathcal F^2\cV$ endowed with the Hodge metric.
\end{theorem}

This theorem relies on an ``elementary'' equidistribution result for positive definite sublattices of a quadratic lattice that we prove in \Cref{refinedloci}. Using a more refined equidistribution result of Eskin--Oh \cite{eskinoh} based on Ratner theory, one can get a more precise equidistribution theorem for the locus where the Néron--Severi group is a fixed quadratic lattice. For a positive definite matrix $M$, we denote by $\mu_1(M)$ the square-root of the smallest non-zero value integrally represented by $M$. Say also that $M$ is primitively represented by $(V_\Z,B)$ if there exists a primitive sublattice of $V_\Z$ of rank $r$ having a basis with intersection matrix equal to $M$.

\begin{theorem}\label{maink3}
Let $\{\mathbb{V}_\Z,\mathcal{F}^\bullet\mathcal{V},B\}$ be a simple $\Z$-PVHS of weight $2$ and Hodge numbers $(q,p,q)$ over a complex analytic variety $S$ of dimension $d=rq$ and which has generically immersive period map. Assume that $p,2q\geq 2$ and $rq<p$. 

Let $(M_n)_{n\in\N}$ be a sequence of positive definite integral matrices of rank $r$ which are primitively represented by $(V_\Z,B)$ and such that $\mu_{1}(M_n)\rightarrow \infty$, as $n\rightarrow \infty$. Then there exists a sequence $\left(a(M_n)\right)_{n\in \N}$ of positive real numbers such that, for every relatively compact open subset $\Omega \subset S$ with boundary of measure $0$, we have:
    \[\frac{1}{a(M_n)}|\{(s,\lambda_1,\cdots,\lambda_r))\in \Omega\times \V^r_{\Z,s},(B(\lambda_i.\lambda_j))=M_n, \lambda_i \in \mathrm{Hdg}(s)\}|\underset{n\rightarrow \infty}{\longrightarrow} \int_{\Omega}c_q(\mathcal{F}^2\cV)^r.\]
\end{theorem}

\begin{remarque}
In the above theorems, one needs to exclude the points belonging to exceptional components of $NL^{\geq r}(S)$ of dimension $\geq 1$ from the counting, i.e., the non-typical ones. The precise versions of these theorems are given in \Cref{sss:ProofEquidistributionLN}.
\end{remarque}

\begin{remarque}
The Siegel--Weil formula gives an arithmetic expression for the asymptotic behaviour of the sequence $a(M_n)$ in \Cref{maink3}. The precise expression and how it grows with $\det(M_n)$ are discussed in \Cref{totalvolume}.
\end{remarque}

\begin{remarque}
A base of dimension $rq$ is the minimal dimension for which a refined Noether--Lefschetz locus is expected to exist for dimension reasons. If the dimension of the base $S$ is greater than $rq$, then Theorem \ref{maink3} gives the equidistribution of $NL^{\geq r}(S)$ towards $c_q(\mathcal F^2 \cV)^r$ in terms of currents (see \Cref{currentsremark}).
\end{remarque}

\begin{remarque}
As soon as $q\geq 2$, Griffiths' transversality combined with the integrability condition of the tangent space to $S$ imply that the dimension of $S$ is at most $\frac{pq}{2}$ by \cite[Theorem 1.1]{carlsonbounds}.
\end{remarque}

\begin{remarque}
The above theorems are already interesting when $r=1$, for which they give the equidistribution of the Noether--Lefschetz locus. The case $r=q=1$ was treated by the first author in \cite{tayouequi}.
\end{remarque}

\begin{remarque}
If the base $S$ of the variation $\mathbb V$ is a complex projective variety, one can apply Theorem \ref{main} to $\Omega=S$ and get an asymptotic estimate of the ``growth'' of the Noether--Lefschetz locus of $S$. We conjecture that the same estimate holds when $S$ is quasi-projective of arbitrary dimension (which implies that $\int_S c_q(\mathcal F^2\cV)^r < +\infty$). The case $q=r=1$ has been settled in \cite{tayouequi}. One could hopefully obtain this global estimate by working with the cohomology of an appropriate compactification of  $\Gamma \backslash G/K$. This raises more general questions that are beyond the scope of this paper.
\end{remarque}

\Cref{main} applies to families of algebraic varieties whenever one has a generic local Torelli theorem. This holds for families of abelian varieties, $K3$ surfaces, hyperkähler manifolds and for projective hypersurfaces by the general result of \cite{donagi}, which yields the following examples:

\begin{itemize}
    \item Smooth quintic surfaces in $\mathbf P^3$ have Hodge numbers
    \[h^{2,0} = 4~, \quad h^{1,1} = 45~.\]
    The moduli space of quintic surfaces has dimension $40$ and satisfies a generic Torelli theorem by \cite{donagi}. Thus Theorems \ref{maink3-alt} and  \ref{maink3} give equidistribution results for the refined Noether--Lefschetz loci on the moduli space of quintic surfaces up to $r=10$.
    
    \item Cubic hypersurfaces in $\mathbf P^7$ have Hodge numbers
    \[h^{6,0}=h^{5,1}=0~,\quad h^{4,2}=8~,\quad h^{3,3} = 178\]
    So their cohomology in degree 6 is a Hodge structure of weight $2$. The moduli space of cubic hypersurfaces in $\mathbf P^7$ has dimension $56$ and satisfies the generic Torelli theorem \cite{donagi}. Thus again Theorems \ref{maink3-alt} and \ref{maink3} give equidistribution results for refined Hodge loci on the moduli space of cubic hypersurfaces of $\mathbf P^7$ up to $r=7$.
    
\end{itemize}
\subsubsection{Hodge loci in Shimura varieties}
For $g\geq 2$, let $\mathcal{A}_g$ be the moduli space of principally polarized complex abelian varieties of dimension $g$. It is well-known that the smallest codimension of a special subvariety is $g-1$ and it is realized for example by  $\mathcal{A}_{g-1}\times\mathcal{A}_{1}$. It is then expected, see \cite[Remark 2.16]{baldiklinglerullmo}, that the Hodge locus is analytically dense in any Hodge generic subvariety of dimension at least $g-1$. 

As a partial answer, we have the following results. Let $\cF^{1}\rightarrow \mathcal{A}_g$ be the Hodge bundle and let $\left(c_{n}(\cF^{1})\right)_{0\leq n\leq g}$ be its Chern forms with respect to the Hodge metric. For $1\leq k\leq g$, define \footnote{This is a particular example of a Schur polynomial.} \[s_k=\det\left(((-1)^{j-i}c_{g-k+j-i}(\cF^{1}))_{1\leq i,j\leq k}\right).\]  Then $s_k$ is a semi-positive form. We prove then the following result.
\begin{theorem}\label{T:densityabelianvarieties}
Let $X\subseteq \mathcal{A}_g$ be a smooth subvariety. If the restriction of $s_k$ to $X$ is non-zero, then the locus of elements in $X$ parameterizing abelian varieties containing a sub-abelian variety of dimension $k$ is analytically dense and equidistributed with respect to $s_k$. In particular, if $X$ is compact and has dimension at least $\frac{(g-1)(g-2)}{2}$, then the Hodge locus is analytically dense in $X$.
\end{theorem}

If $k=1$, then the above theorem yields that the Hodge locus is dense if $c_{g-1}$ is non-zero by restriction to $X$. This prompts the following question.
\begin{question-nr}\label{c:abelian-var}
Let $X$ be a Hodge generic subvariety of $\mathcal{A}_g$ of dimension at least $g-1$. Is the restriction of $c_{g-1}$ to $X$ always non-zero?
\end{question-nr}
\begin{remarque}\label{rq:a3-a4}
The assumption of Hodge genericity is necessary. Indeed, if $X=\mathcal{A}_2\times \{pt\}\subseteq \mathcal{A}_4$. Then $X$ has dimension $3\geq 4-1=3$ but the the restriction of $c_3$ to $X$ vanishes. Indeed, the restriction of the vector bundle $\cF^1$ to $X$ splits as a direct sum of two vector bundles of rank $2$, one of them being trivial. Hence $c_{3}$ vanishes on $X$.
\end{remarque}

\Cref{T:densityabelianvarieties} admits the following generalization. 
\begin{theorem}\label{t:shimura_simple}
Let $S$ be a connected Shimura variety associated to a connected Shimura datum $(G,\D)$. Assume that there exists a Shimura sub-datum $(H,\D_H)$ such that $\pi_*p^* \omega_{G/H}$ is positive of type $(k,k)$. Then the Hodge locus is dense in any subvariety of dimension at least $k$ and equidistributed with respect to $\pi_*p^* \omega_{G/H}$. In particular, if $G$ is absolutely simple and has a Shimura curve associated to a Shimura subgroup $H$, then the Hodge locus is dense in any hypersurface and equidistributed w.r.t $\pi_*p^* \omega_{G/H}$.
\end{theorem}

In particular, for unitary Shimura varieties, we obtain the following.
\begin{corollaire}\label{C:unitary-shimura}
Let $S$ be a Shimura variety of unitary type $(n,1)$. Then the typical Hodge locus is dense and equidistributed in any subvariety of $S$ of positive dimension.
\end{corollaire}
Let $S$ be a Shimura variety associated to a connected Shimura datum $(G,\D_H)$ and let $k$ be the minimal integer for which there exists a sub-Shimura datum $(H,\D_H)$ such that $\pi_*p^*\omega_{G/H}$ is of type $(k,k)$. Then \Cref{c:abelian-var} has the following generalization.
\begin{question-nr}
Let $X\subset S$ be a Hodge generic subvariety of dimension at least $k$. Is the restriction of $\pi_*p^*\omega_{G/H}$ to $X$ always non-zero?  
\end{question-nr}

\subsubsection{Equidistribution of families of CM points in Shimura varieties}
Another application of \Cref{main} is an equidistribution result for families of CM points in some Shimura varieties. Several results about the equidistribution of CM points are known, see for example \cite{dukehyperbolic,zhang,khayutin} and in general, the following conjecture is widely open, see \cite[Conjecture 2.6]{yafaev} for more details on this conjecture.
\begin{conjecture} \label{conj:EquidistributionCM}
Let $S$ be a Shimura variety over $\overline{\Q}$ and let  $(x_n)_{n\in\N}$ be a generic sequence of CM points in $S$. Then the sequence of Galois orbits $\mathrm{Aut}(\overline{\Q}/\Q)\cdot x_n$ is equidistributed in $S(\C)$. 
\end{conjecture}
In what follows, we will state the result we prove in the simplest case of Siegel Shimura varieties, referring to \Cref{equidistributionCM} for the most general statement. 

Let $g\geq 1$. A polarized isogeny $f:(A_1,\omega_1)\rightarrow (A_2,\omega_2)$ of similitude factor $N\geq 1$ between two principally polarized abelian varieties of dimension $g$ is an isogeny which satisfies $f^*\omega_2=N\,\omega_1$. If $A_1=A_2$,  we  say moreover that $f$ is regular if the centralizer of the homological realization of $f$ in $\mathrm{GSp}_{2g}(H^{1}(A,\R))$ is a torus. This implies that $A_1$ is a \emph{CM abelian variety}\footnote{Short for ``abelian variety with complex multiplication''}, meaning that $\mathrm{End}(A)_\Q$ is a commutative algebra of degree $2g$ over $\Q$, i.e., the maximal possible dimension. Conversely, every CM polarized Abelian variety admits a regular self-isogeny (\Cref{transverselocusCM}). An equivalent characterization of CM abelian varieties is that their Mumford--Tate group $MT(A)$ (which is contained in the centralizer of the isogeny $f$) is a torus, see \Cref{CM}.

Let $\omega$ be the first Chern form of the Hodge bundle on $\mathcal{A}_g$. It is well know that $\omega$ is a K\"ahler form. If $A$ is a principally polarized abelian variety, we denote by $\dagger:\mathrm{End}(A)_\Q\rightarrow \mathrm{End}(A)_\Q$ the Rosati involution. As an application of \Cref{main}, and of the main result of \cite{clozelohullmo} (see also \cite{eskinoh1}), we get an equidistribution result for CM abelian varieties admitting self-isogenies of fixed degree.
\begin{theorem}\label{CMabelianvarieties}
There exists a sequence $b(N)$ such that, for every relatively compact open subset $\Omega \subset \mathcal A_g$ with boundary of measure $0$, we have:
\[\left|\{(A,f), A\in \Omega, f\in\mathrm{End}(A),\, f^\dagger\circ f=N\mathrm{Id},\, \textrm{and $f$ regular}\}\right|\underset{N\to +\infty}{\sim}b(N) \int_\Omega \omega^{\frac{g(g+1)}{2}}~.\]

\end{theorem}

This theorem does not answer \Cref{conj:EquidistributionCM} because, as $N$ grows, our equidistributing sets are the union of an increasing number of Galois orbits. It is however sharper than other more elementary equidistribution results. For comparison, in the case of $g=1$, Conjecture \ref{conj:EquidistributionCM} is answered positively by Duke's equidistribution theorem for CM elliptic curves with fundamental discriminant $N$ \cite[Thereom 1]{dukehyperbolic} and by Clozel--Ullmo in general \cite[Théorème 2.4]{clozelullmohecke}, while an elementary counting argument easily gives the equidistribution of CM curves with discriminant $\leq N$. Our theorem lies in-between: it asserts that the set of CM elliptic curves with discriminant of the form $N-4a^2$ for all integers $0\leq a\leq \sqrt{N}$ is equidistributed when $N$ goes to $+\infty$.

\subsubsection{Equidistribution of Hecke translates}
We mention two further applications of \Cref{main} that we obtain. The first one is related to the dynamics of Hecke translates in $\mathcal{A}_g$.

Let $S$ and $D$ be two subvarieties of $\mathcal{A}_g$ of complimentary dimensions such that $S$ has dimension $d\leq 2$. Let $\omega$ be as before the first Chern form of the Hodge bundle on $\mathcal{A}_g$. Also, if $(s,d)\in S\times D$, with corresponding abelian varieties $A_s$ and $A_d$, we denote by $\mathrm{Isog}^N(A_s,A_d)$ the set of isogenies from $A_s$ to $A_d$ of similitude factor $N$. An isogeny $f:A_{s}\rightarrow A_d$ is said to be \emph{transverse}, if it does not admit first order deformation in $S\times D$. Then we prove the following, see
\Cref{hecketranslatestorelli} for more details.
\begin{theorem}\label{Torelliequidistribution}

There exists a sequence $(c(N))_{N\geq 1}$ of positive real numbers such that for every relatively compact open subsets $\Omega \subset S$, $\Omega'\subset D$ with boundary of measure $0$, we have:
\[|\{(s,d,f)|\, (s,d)\in \Omega\times \Omega', f\in \mathrm{Isog}^N(A_s,A_d)\,\textrm{regular}\}|\underset{N\rightarrow \infty}{\sim} c(N)\int_\Omega \omega^d\int_{\Omega'} \omega^{\frac{g(g+1)}{2}-d}~.\]
In particular, the locus of points in $S$ isogenous to a point in $D$ is analytically dense and equidistributed in $S$.
\end{theorem}
We have the following corollary.
\begin{corollaire}\label{Corollary:Torelli}
Let $S\subset \mathcal{A}_4$ be a curve. Then the locus of points in $S$ isogenous to the Jacobian of a curve is analytically dense and equidistributed in $S$.
\end{corollaire}

\subsubsection{Equidistribution in cohomology}\Cref{main} yields that cohomology classes of $\Gamma \backslash G/K$ represented by an equidistributing sequence of locally homogeneous submanifolds converge after normalization to the cohomology class of a locally invariant form. (See \Cref{cohomologyequidistribution} for a precise statement.)

 To illustrate this in a specific example, we use the same notation as in \Cref{section:hodgeloci}. Then we have a family of \emph{special cycles} in $\Gamma\backslash G/K\simeq \Gamma\backslash \mathcal D$ where $\mathcal D$ is the period domain, defined as follows : for $r\geq 1$, $\underline{\lambda}_0\in V_\Z^r$, $H=\mathrm{Stab}(\underline{\lambda}_0)$, $M\in M_r(\Z)$ semi-positive definite matrix with rank $r(M)$, and $V_M\equaldef \{\underline{\lambda}\in V_\Z^r,\,\left(B(\lambda_i.\lambda_j)\right)_{1\leq i,j\leq r}=M\}$, let
\[\mathcal{Z}(M)\equaldef \Gamma\backslash\left(\cup_{\underline{\lambda}\in V_M}\{x\in D, x\bot \lambda_i,\, \forall i=1,\dots,r\}\right)\hookrightarrow \Gamma \backslash \mathcal D.\]
Let $c_q(\mathcal{F}^2\mathcal{V})$ be as before the top Chern form of the vector bundle $\mathcal{F}^2\mathcal{V}$. 
\begin{proposition}\label{kudlamillson}
Let $(M_n)_{n\in\N}$ be a sequence of positive definite matrices primitively represented by $(V_\Z,B)$ such that $\mu_1(M_n)\rightarrow \infty$, as $n\rightarrow \infty$. Then \[\mathcal{Z}(M_n)\underset{n\rightarrow \infty}{\sim}  a(M_n)\, c_q(\mathcal{F}^2\mathcal{V})^r.\]
\end{proposition}

This result is reminiscent of the work of Kudla--Millson on modularity of special cycles, see \cite{kudlamillson}, see also \cite{luisgarcia} for a recent approach using superconnections. Indeed, in both these papers, it is proved that the formal generating series:
\[\sum_{M\geq 0}\mathcal{Z}(M)\cup c_q(\mathcal{F}^2\mathcal{V})^{r-r(M)}e^{\tr(2i\pi M\tau)},\, \tau \in\mathfrak{H}_r\]
is a Siegel modular form valued in $H^{2qr}(\Gamma\backslash \mathcal D,\R)$. Here $\mathfrak{H}_r$ is the Siegel upper half space. Hence the knowledge of the structure of the space of Siegel modular forms allows in principle to give an asymptotic formula of $\mathcal{Z}(M)$ in terms of the constant term $c_q(\mathcal{F}^2\mathcal{V})^r$. As $r$ grows, the structure of the space of Siegel modular forms becomes complicated to analyze and much work is needed to derive formulas similar to ours through this approach. Our method yields a straightforward estimate on the asymptotic growth of $\mathcal{Z}(M)$ without using these results. One might even hope that this asymptotic estimate could help understand the cusp structure of Kudla--Millson's modular forms.

\medskip
\subsection{Related work}
The distribution of the Hodge locus has been investigated independently and concomitantly by Baldi--Klingler--Ullmo in \cite{baldiklinglerullmo}. Additionally to striking results on the atypical Hodge locus (see also \cite{klinglerotwinowska} for prior work), they prove several properties about the typical Hodge locus that echo the present work, namely that the typical Hodge locus is either empty or dense, and is always empty when the level is at least $3$. 

Several results analogous to \Cref{maink3} for algebraic families parameterized by Shimura varieties have been settled in arithmetic situations over rings of integers of number fields and over curves defined over finite fields: indeed Charles proved \cite{charles1} that there are infinitely many places where the reduction of two elliptic curves are isogenous ; Shankar and Tang \cite{shankartang} proved that an abelian surface over a number field with real multiplication has infinitely many specializations which are isogenous to the self-product of an elliptic curve and in collaboration with Maulik in \cite{maulikshankartang1} they derive similar results for ordinary abelian surface over the function field of a curve over a finite field. Finally the analogous statement of \Cref{maink3} for K3 type variations of Hodge structures over curves has been proved in the number field setting in \cite{sstt,tayoubad} and over curves over finite fields in \cite{maulikshankartang2}. It is thus interesting to further explore other analogous statements of \Cref{main}, \Cref{maink3} and \Cref{Torelliequidistribution} over number fields and function fields situations.

\subsection{Organization of the paper}
In \Cref{ratner}, we introduce the setting of homogeneous dynamics and explain how to reformulate an equidistribution theorem in terms of currents. In \Cref{equiintersectionpoints}, we deduce our general theorem for transverse intersections of locally homogeneous subspaces with a fixed analytic subvariety, proving \Cref{main}. \Cref{s:pullpush} is devoted to the study of the pull-push form. In particular, we explain how to interpret its cohomology class via compact duality of homogeneous spaces. In \Cref{definitionvhs}, we discuss the pull-push forms in the case of period domains of variations of Hodge structures and relate them to Chern classes of Hodge bundles, allowing us to compute these forms explicitly in the setting of a variation of Hodge structure of weight $2$ and Shimura varieties. Finally, in \Cref{vhs}, we discuss our applications : the study of refined Noether--Lefschetz loci, the equidistribution of some families of CM points in Shimura varieties and the equidistribution of intersection points of Hecke translates of the Torelli locus with a curve and a surface in $\mathcal{A}_g.$ 

\subsection{Acknowledgments}
We thank Bruno Klingler for suggesting the use of o-minimality and Gabriele Mondello for pointing out the Giambelli--Porteous--Thom formula. The application to Hecke translates of the Torelli locus was raised during a fruitful discussion with Ananth Shankar and Yunqing Tang, to whom we are grateful. We also thank Olivier Benoist, Nicolas Bergeron, François Charles, Phillip Griffiths and Shou-Wu Zhang for many useful discussions. S.T. also thanks CRM in Montréal and IAS in Princeton for excellent working conditions. We also thank the referees for their careful reading of the paper and useful comments.

\section{Equidistribution in terms of currents}\label{ratner}
In this section, we recall some background on convergence of measures, currents and homogeneous spaces. Then we reformulate equidistribution results in homogeneous dynamics in terms of weak convergence of currents. Finally, we recall some equidistribution results from Ratner's work. 

\subsection{Convergence of measures}

Let us start by recalling a few classical facts in measure theory which will be used mainly in \Cref{equiintersectionpoints}.

Let $S$ be an analytic subset of dimension $d$ of a manifold $M$ whose smooth locus is oriented, and let $\omega$ be a smooth form of degree $d$ on $M$. Then the restriction of $\omega$ to the smooth locus of $S$ defines a signed Radon measure on $S$. This measure is \emph{regular} in the sense that:
\begin{enumerate}
    \item For every open subset $U$ of $S$ and every sequence of compact sets $(C_n)_{n\in\N}$ with $C_n \subset \mathring{C}_{n+1}$ and $\bigcup_{n\in \N} C_n = U$, we have
    \[\int_U \omega = \lim_{n\to +\infty} \int_{C_n} \omega~.\]
    \item For every compact subset $C\subset S$ and every sequence of open sets $(U_n)_{n\in\N}$ with $\bar U_{n+1} \subset U_n$ and $\bigcap_{n\in \N} U_n=C$, we have
    \[\int_C \omega = \lim_{n\to +\infty} \int_{U_n} \omega~.\]
\end{enumerate}

In the absence of any precision, we say that a set $A\subset S$ has measure $0$ if the intersection of $A$ with the smooth locus of $S$ has Lebesgue measure $0$ in any coordinate chart. This implies that its measure with respect to $\omega$ is $0$.

Given a sequence of signed Radon measures $\mu_n$, we have the following equivalent characterizations of the convergence of $\mu_n$ to $\omega$:
\begin{itemize}
    \item[(i)] for every continuous function $f:S\to \R$ with compact support,
    \[\mu_n(f) \underset{n\to +\infty}{\longrightarrow} \int_S f \omega~,\]
    \item[(ii)] for every relatively compact open subset $\Omega$ of $S$ with boundary of measure $0$,
    \[\mu_n(\Omega) \underset{n\to +\infty}{\longrightarrow} \int_\Omega\omega~.\]
\end{itemize}

We then say that $\mu_n$ converges weakly to $\omega$ and we write \[\mu_n\underset{n\to +\infty}{\rightharpoonup}\omega.\] 

We will say that $\mu_n$ converges weakly to $\omega$ \emph{on an open subset $U$} if the restriction of $\mu_n$ to $U$ converges weakly to the restriction of $\omega$. Equivalently, $\mu_n$ converges weakly to $\omega$ on $U$ if Property (i) above holds for any function with compact support inside $U$.

The following standard facts will be useful in proving weak convergence of measures:

\begin{proposition} \label{p:ConvergenceMeasureCovering}
Let $(U_i)_{i\in I}$ be an open covering of $S$. Then $\mu_n$ converges weakly to $\omega$ on $S$ if and only if converges weakly to $\omega$ on $U_i$ for all $i\in I$.
\end{proposition}

\begin{proposition} \label{p:ConvergenceMeasureComplement}
Let $Z$ be a closed subset of $S$ of measure $0$. Assume that ${\mu_n}$ converges weakly to $\omega$ on $Z^c$. Then the following are equivalent:
\begin{itemize}
    \item[(i)] $\mu_n$ converges weakly to $\omega$ on $S$,
    \item[(ii)] for very compact subset $C \subset Z$ and every $\epsilon >0$, there exists an open neighborhood $U_\epsilon$ of $C$ such that
    \[\limsup_{n\to +\infty} \vert \mu_n\vert (U_\epsilon) \leq \epsilon~.\]
\end{itemize}
\end{proposition}
Here, $\vert \mu_n\vert = \mu_n^+ + \mu_n^-$ denotes the total variation measure of $\mu$.

\subsection{Currents on manifolds}\label{currents}
For more details on this section, we refer to \cite[Chapter 3, section 1]{griffithsharris}.

Let $M$ be a real manifold of dimension $n$. For $k\geq 0$, let $\Omega_{c}^{k}(M)$ be the vector space of $C^{\infty}$ differential forms on $M$ of degree $k$ with compact support. It is endowed with its natural topological space structure making it a Fréchet space. 
\begin{definition}
A current of degree $k$ on $M$ is a continuous linear form on $\Omega_{c}^{n-k}(M)$. The space of currents of degree $k$ on $M$ is denoted by $\mathcal{D}^k(M)$.
\end{definition}
\begin{example}

    If $N\hookrightarrow M$ is an oriented properly immersed submanifold of codimension $k$ of $M$, the \emph{integration current} $T_N\in \mathcal{D}^k(M) $ is defined by
\[T_N(\beta) = \int_N \beta, \quad \beta\in \Omega_{c}^{n-k}(M).\]
\end{example}

\begin{example} \label{ex:Form->Current}
A differential $k$-form $\alpha$ induces a $k$-dimensional current $T_\alpha$ defined by
\[T_\alpha(\beta) = \int_M \alpha \wedge \beta, \quad \beta\in \Omega_{c}^{n-k}(M).\]
\end{example}
The exterior derivative $d$ on differential forms induces a map \[d:\mathcal{D}^k(M)\rightarrow \mathcal{D}^{k+1}(M)\] 
defined by \[dT(\phi)=(-1)^{k+1}T(d\phi),\quad \phi\in\Omega_{c}^{n-k-1}(M).\]
A current $T$ is \emph{closed} if $dT=0$.

The exterior derivative defines a cochain complex structure on $(\mathcal{D}^\bullet(M))$, and Example~\ref{ex:Form->Current} gives a morphism of cochain complexes \[(\Omega^\bullet(M),d)\rightarrow (\mathcal{D}^\bullet(M),d).\]

The previous morphism is in fact a quasi-isomorphism (i.e. it induces isomorphisms at the level of cohomology groups, see \cite[p.382]{griffithsharris}).

The space $\mathcal{D}^k(M)$ is naturally a topological vector space when equipped with the weak topology: a sequence $T_n$ of degree $k$ currents \emph{converges weakly} to a current $T$ (which we write $T_n \rightharpoonup T$) if
\[T_n(\beta)\underset{n\to +\infty}{\longrightarrow} T(\beta)\]
for all $\beta \in \Omega_c^{n-k}(M)$.

Assume now that $M$ is a complex manifold of complex dimension $n$. Then the complex $\D^\bullet(M)$ admits a bigrading
\[\D^k(M) = \bigoplus_{p+q = k} \D^{p,q}(M)~,\]
where $\D^{p,q}(M)$ is the topological dual of the complex vector space $\Omega_c^{n-p,n-q}(M)$. 

In particular, if $Z\subseteq M$ is a closed complex analytic subvariety of complex codimension $k$, we can similarly define a closed integration current $T_Z\in D^{k,k}(M)$ by integrating over (the smooth locus of) $Z$.

\subsection{Homogeneous spaces, orientations and volume forms}\label{notations}
In this section, we introduce the notations that we will be using throughout the paper and recall some facts on volume forms on Lie groups.

Let $G$ denote a real algebraic semi-simple Lie group. Suppose we are given the following subgroups of $G$:
\begin{itemize}
    \item a lattice $\Gamma$,
    \item a semi-simple Lie subgroup $H$ without compact factor,
    \item a compact subgroup $K$.
\end{itemize}

Let $L$ be the intersection of $K$ and $H$. This is a compact subgroup of $H$. The group $H$ acts on the right on the quotient $\Gamma \backslash G$, and we will be interested in the next section in equidistribution properties of orbits of this action and their projection to $\Gamma \backslash G/K$.

Let $\fg$, $\fh$, $\fk$ and $\fl$ denote respectively the Lie algebras of $G$, $H$, $K$, and $L$. Up to taking subgroups of index $2$, we can assume that the adjoint actions of $G$, $H$, $K$, and $L$ have determinant $1$. We then fix once and for all some orientation of $\fg$, $\fh$, $\fk$ and $\fl$ and orient accordingly the quotient spaces $\fg/\fk$, $\fg/\fh$, $\fg/\fl$, $\fk/\fl$ and $\fh/\fl$. Those orientations induce orientations on $G/K$, $G/H$, $G/L$, $K/L$ and $H/L$ respectively. 

Recall that the Lie algebra $\fg$ carries a natural symmetric bilinear form called the \emph{Killing form}, which is invariant under the adjoint action and non-degenerate since $G$ is semi-simple.  Its restriction to $\fh$, $\fk$ or $\fl$ is still non-degenerate. We denote by 
$\omega_G$, $\omega_H$, $\omega_K$ and $\omega_L$ the volume forms associated to the restricted Killing metric and the prescribed orientations on $\fg$, $\fh$, $\fk$ and $\fl$ respectively, as well as the induced bi-invariant volume forms on $G$, $H$, $K$, and $L$. Finally, we denote by $\omega_{G/K}$, $\omega_{G/H}$, $\omega_{G/L}$, $\omega_{K/L}$ and $\omega_{H/L}$ the invariant volume forms on the corresponding homogeneous spaces induced by $\omega_G$, $\omega_H$, $\omega_K$ and $\omega_L$.

The volume forms $\omega_G$ and $\omega_{G/K}$ respectively factor to volume forms on $\Gamma \backslash G$ and $\Gamma \backslash G/K$, and we define respectively
\[\Vol(\Gamma \backslash G) = \int_{\Gamma \backslash G}\omega_G~,\quad\Vol(K) = \int_{K}\omega_K~,\quad
\Vol(\Gamma \backslash G/K) = \int_{\Gamma \backslash G/K}\omega_{G/K}~.\]
The compatibility of the volume forms gives the following identity:
\[\Vol(\Gamma \backslash G) = \Vol(\Gamma\backslash G/K)\cdot \Vol(K)~.\]

Similarly, if $xH$ is a closed $H$-orbit in $\Gamma\backslash G$, then $\omega_H$, $\omega_L$ and $\omega_{H/L}$ induce volume forms on $xH$, $L$ and $xH/L$. We denote respectively by $\Vol(xH)$, $\Vol(L)$ and $\Vol(xH/L)$ their total mass, and we have the following identity:
\[\Vol(xH) = \Vol(xH/L)\cdot \Vol(L)~.\]

\subsection{Equidistribution in terms of currents} \label{ss:EquidistributionCurrents}

In this section, we reformulate equidistribution results of sequences of $H$-orbits in terms of convergence of currents.
As before, for all $n\in \N$, let $\mathcal O_n$ be a finite union of closed $H$-orbits in $\Gamma \backslash G$. The starting point of our work is the remark that the equidistribution of a sequence $\{\mathcal O_n\subset \Gamma\backslash G,n\in\N\}$ can be reformulated as an equidistribution of the currents of integration over $\mathcal O_n$.

To be more precise, let $p$ denote the projection from $G/L$ to $G/H$ and $\pi$ the projection from $G/L$ to $G/K$. 

\begin{lemma} \label{l:EquidistributionCurrentsG/L}
Assume that the sequence $(\mathcal O_n)_{n\in\N}$ is equidistributed in $\Gamma \backslash G$ and let $\hat T_{\mathcal O_n/L}$ denote the integration current on $\mathcal O_n/L \subset \Gamma \backslash G/L$. Then
\[\frac{1}{\Vol(\mathcal O_n/L)}\hat T_{\mathcal O_n/L} \rightharpoonup \frac{1}{\Vol(K/L)\Vol(\Gamma \backslash G/K)}  p^*\omega_{G/H}~.\]
\end{lemma}

The form $p^*\omega_{G/H}$ can be seen as a transverse volume form to the foliation of $G/L$ by translates of $H/L$. Note that, by applying the lemma to $L= \{\mathbf{1}_G\}$, we get the following corollary:
\begin{corollaire} \label{t:EquidistributionCurrentsG}
Assume that the sequence $(\mathcal O_n)_{n\in\N}$ is equidistributed in $\Gamma\backslash G$ and let $T_{\mathcal{O}_n}$ denote the integration current on $\mathcal O_n$. Then
\[\frac{1}{\Vol(\mathcal{O}_n)}T_{\mathcal O_n} \rightharpoonup \frac{1}{\Vol(\Gamma \backslash G)}\hat p^*\omega_{G/H}~,\]
where $\hat p$ is the projection from $G$ to $\Gamma \backslash G$.
\end{corollaire}

Finally, one can push this equidistribution forward by the fibration map from $G/L$ to $G/K$. Recall that the map $\pi$ induces a push-forward map
\[\pi_*: \Omega^\bullet(G/L) \to \Omega^{\bullet-\dim(K/L)} (G/K)\]
which is by definition Poincaré dual to the pull-back map $\pi^*$, i.e.
\[\int_{G/K}(\pi_*\alpha) \wedge \beta = \int_{G/L} \alpha \wedge (\pi^*\beta)\]
for all $\beta$ with  compact support (for more details, see \Cref{pushforward} or \cite[p.37]{bottTu}). We still denote by $\pi$ and $\pi_*$ the factorization of those maps by the left action of $\Gamma$.

\begin{theorem} \label{t:EquidistributionCurrentsG/K}
Assume the sequence $(\mathcal O_n)_{n\in\N}$ equidistributes in $\Gamma\backslash G$ and let $T_{\mathcal O_n/L}$ denote the integration current on $\mathcal O_n/L \subset \Gamma \backslash G/K$. Then
\[\frac{1}{\Vol(\mathcal O_n/L)}T_{\mathcal O_n/L} \rightharpoonup \frac{1}{\Vol(K/L)\Vol(\Gamma \backslash G/K)}\pi_* p^*\omega_{G/H}~.\]

\end{theorem}

\begin{proof}
Let $\mathcal H$ denote the foliation of $\Gamma \backslash G/L$ by the left translates of $H/L$ and $T\mathcal H \subset T(\Gamma \backslash G/L)$ the tangent distribution to this foliation. The volume form $\omega_{H/L}$ on $H/L$ defines a smooth section $\omega_{\mathcal H}$ of $\Lambda^{\max}T^*\mathcal H$.

Let $\alpha$ be a form of degree $\dim(G)-\dim(H)$ with compact support on $\Gamma \backslash G/L$, and let $f$ be the smooth function with compact support such that $\alpha_{\vert T\mathcal H} = f \omega_{\mathcal H}$. The compatibility between the various volume forms gives
\[\alpha \wedge p^*\omega_{G/H} = f \omega_{G/L}~.\]

Let $p_0$ denote the projection from $\Gamma \backslash G$ to $\Gamma \backslash G/L$. We have
\begin{eqnarray*}
\frac{1}{\Vol(\mathcal{O}_n/L)}\int_{\mathcal{O}_n/L} \alpha & = & \frac{1}{\Vol(\mathcal{O}_n/L)}\int_{\mathcal{O}_n/L} f~\omega_{\mathcal H} \\
& =& \frac{1}{\Vol(\mathcal{O}_n)}\int_{\mathcal{O}_n} f\circ p_0~\omega_H \\
&\underset{n\to +\infty}{\longrightarrow}& \frac{1}{\Vol(\Gamma \backslash G)}\int_{\Gamma \backslash G} f\circ p_0~\omega_G
\end{eqnarray*}
and 
\begin{eqnarray*}
 \frac{1}{\Vol(\Gamma \backslash G)}\int_{\Gamma \backslash G} f\circ p_0~\omega_G &=& \frac{1}{\Vol(\Gamma \backslash G/L)}\int_{\Gamma \backslash G/L} f~\omega_{G/L} \\
&=&\frac{1}{\Vol(\Gamma \backslash G/L)}\int_{\Gamma \backslash G} \alpha \wedge p^*\omega_{G/H}~. \\
\end{eqnarray*}

This shows that $\frac{1}{\Vol(\mathcal{O}_n/L)}\hat T_{\mathcal{O}_n/L}$ converges weakly to $\frac{1}{\Vol(\Gamma \backslash G/L)} \omega_{G/H}$. 
\end{proof}

\begin{proof}[Proof of Theorem \ref{t:EquidistributionCurrentsG/K}]
We push forward the equidistribution of Lemma \ref{l:EquidistributionCurrentsG/L} to $\Gamma \backslash G/K$. Note that we have $T_{\mathcal{O}_n/L} = \pi_* \hat T_{\mathcal{O}_n/L}$.

Let $\alpha$ be a form of degree $\dim(H/L)$ with compact support on $\Gamma \backslash G/K$. We then have
\begin{eqnarray*}
\frac{1}{\Vol(\mathcal{O}_n/L)}\langle T_{\mathcal{O}_n/L}, \alpha\rangle & = & \frac{1}{\Vol(\mathcal{O}_n/L)} \langle \hat T_{\mathcal{O}_n/L}, \pi^*\alpha \rangle \\
&\to & \frac{1}{\Vol(\Gamma \backslash G/L)}\int_{\Gamma \backslash G/L} \pi^* \alpha \wedge p^* \omega_{G/H}\\
&=& \frac{1}{\Vol(\Gamma \backslash G/L)} \int_{\Gamma \backslash G/K} \alpha \wedge \pi_* p^* \omega_{G/H} \textrm{ (by definition of $\pi_*$).}
\end{eqnarray*}
\end{proof}

For $n\in \N$, the integration current $T_{\mathcal{O}_n/L}$ is closed since $\mathcal O_n/L$ has empty boundary. It thus has a well defined cohomology class $[\mathcal{O}_n/L]\in H^{d_H}(\Gamma\backslash G/K,\R)$ where $d_H$ is the codimension of $H/L$ in $G/K$.

\begin{corollaire}\label{cohomologyequidistribution}
Assume that the sequence $(\mathcal O_n)$ is equidistributed in $\Gamma\backslash G$ and let $[\mathcal{O}_n/L]$ denote the corresponding cohomology class in $H^{d_H}(\Gamma\backslash G/K,\R)$. Then
\[\frac{1}{\Vol(\mathcal O_n/L)}[\mathcal O_n/L] \underset{n\rightarrow \infty}{\rightarrow} \frac{1}{\Vol(K/L)\Vol(\Gamma \backslash G/K)}[\pi_* p^*\omega_{G/H}]~.\]
\end{corollaire}

\subsection{Ratner theory and its consequences}\label{ratnersetup}
In this final section, we recall some equidistribution results à la Ratner of homogeneous orbits in locally homogeneous spaces. These results will be used when studying the refined Noether--Lefschetz locus in \Cref{refinedloci}.

We place ourselves in an arithmetic setting. Though this is not a requirement of Ratner theory, it is the source of its most remarkable consequences and will be sufficient for our applications.

We thus assume that we are given an inclusion $\mathbf H\subset \mathbf G$ of semi-simple algebraic groups over $\Q$ such that
\begin{itemize}
    \item $G$ is a subgroup of $\mathbf G_\R$ containing $\mathbf G_\R^0$,
    \item $H= G\cap \mathbf H_\R$,
    \item $\Gamma$ is commensurable to $\rho^{-1}(\GL(n,\Z))$, for some faithful representation $\phi: \mathbf G \to \GL(n)$ over $\Q$.\footnote{The commensurability class of $\rho^{-1}(\GL(n,\Z))$ is independent of the representation $\rho$.}
    \end{itemize}

By a lemma of Chevalley (see \cite[Proposition 4.6]{BenoistReseaux}), we can find a free $\Z$-module $\mathbf V_\Z$ and an embedding $\mathbf G\hookrightarrow \GL(\mathbf V_\Q)$ such that $\mathbf H$ is the stabilizer in $\mathbf G$ of a vector $v_0\in \mathbf V_\Z$. We can moreover assume that the $G$-orbit of $v_0$ is Zariski closed and that $\Gamma$ preserves $\mathbf V_\Z$. For every $\lambda \in \R_{>0}$, we can now identify the homogeneous space $G/H$ with the $G$-orbit of $\lambda v_0$ in $\mathbf V_\R$.

The following classical result of homogeneous dynamics establishes the equivalence between closed $H$-orbits in $\Gamma \backslash G$ and discrete $\Gamma$-orbits in $G/H$.

\begin{lemme}
Let $g$ be an element in $G$, let $x$ be its projection to $\Gamma \backslash G$ and $v = gv_0$ its projection to $Gv_0 = G/H$. Then the following are equivalent:
\begin{itemize}
    \item the set $\Gamma g H$ is closed in $G$,
    \item the right $H$-orbit of $x$ is closed in $\Gamma \backslash G$,
    \item the left $\Gamma$-orbit of $v$ is discrete in $G/H$,
    \item the group $\Gamma \cap gHg^{-1}$ is a lattice in $gHg^{-1}$,
    \item there exists $\lambda \in \R_{>0}$ such that $\lambda v \in V_\Z$.
\end{itemize}
\end{lemme}

Now, let now $V_n$ be a finite union of discrete $\Gamma$-orbits in $G/H$ and let $\mathcal O_n$ be the corresponding finite union of closed $H$-orbits in $\Gamma \backslash G$. 
The volume form $\omega_H$ induces an $H$-invariant measure $\nu_n$ on $\Gamma \backslash G$, supported by $\mathcal O_n$, whose total mass $\Vol(\mathcal O_n)$ is finite.

\begin{definition}\label{normalization-convention}
We say that the sequence $(\mathcal O_n)_{n\in\N}$ is \emph{equidistributed} in $\Gamma \backslash G$ if the sequence of probability measures
\[\frac{1}{\Vol(\mathcal O_n)} \nu_n\] converges weakly to
\[\frac{1}{\Vol(\Gamma \backslash G)}\omega_G~.\]

We say that the sequence $(V_n)_{n\in\N}$ is \emph{equidistributed in $G/H$} if the discrete measure
\[\frac{1}{\Vol(\mathcal O_n)} \sum_{x\in V_n} \delta_x\]
converges weakly to $\omega_{G/H}$.
\end{definition}

Recall the following classical lemma from \cite[Proposition 2.2]{eskinoh}
\begin{lemme}\label{Sarnaklemma}
The sequence $(V_n)_{n\in\N}$ is equidistributed in $G/H$ if and only if the sequence $(\mathcal O_n)_{n\in\N}$ is equidistributed in $\Gamma \backslash G$.
\end{lemme}

In \cite{eskinoh}, Eskin and Oh give an equidistribution criterion for finite unions of closed orbits which relies on Ratner's groundbreaking work on unipotent dynamics as well as further developments by Mozes--Shah \cite{mozesshah} and Dani-Margulis \cite{danimargulis1,danimargulis1}.

Let $\mathcal O_n = \bigcup_{i=1}^{l_n} x_{i,n} H$ be a finite union of closed $H$-orbits  in $\Gamma \backslash G$ and let $V_n = \bigcup_{i=1}^{l_n} \Gamma v_{i,n}$ be the corresponding union of discrete $\Gamma$-orbits in $G/H$.

\begin{definition} \label{d:MassLoss}
We say that the sequence $(\mathcal O_n)_{n\in N}$ has no loss of mass if for all compact subsets $C$ of $G/H$,
\[\frac{1}{\Vol(\mathcal O_n)}\left(\sum_{\underset{x_{i,n}H\cap C = \emptyset}{i}} \Vol\left(x_{i,n}H\right)\right)\underset{n\to +\infty}{\longrightarrow} 0~.\]
\end{definition}

\begin{definition} \label{d:NonFocused}
The sequence $(\mathcal O_n)$ is called \emph{focused} if there exists $g\in G$ and a subgroup $H'$ of $G$ containing $gHg^{-1}$ and defined over $\Q$ such that 
\[\limsup_{n\to +\infty} \frac{1}{\Vol(\mathcal O_n)}\left(\sum_{\underset{\Gamma v_{i,n}\subset \Gamma H' g Z(H) v_0}{i}} \Vol(x_{i,n}H)\right) > 0~.\]
It is called \emph{non-focused} otherwise.
\end{definition}

\begin{remarque}
Eskin--Oh's original definition of being non-focused  combines both definitions \ref{d:MassLoss} and \ref{d:NonFocused}. It is more convenient to us to separate them, since we shall verify both conditions independently.
\end{remarque}

\begin{theorem}[Theorem 1.13 from \cite{eskinoh} ]\label{eskin-oh}
Assume that $H$ is a semi-simple subgroup of $G$ without compact factors. Then the sequence $(V_n)_{n\in\N}$ is equidistributed in $G/H$ if and only if it is non-focused and has no loss of mass.
\end{theorem}

Note that a sequence of closed $H$-orbits of $\Gamma \backslash G$ leaving every compact subset can only exist if $H$ is contained in a proper parabolic subgroup of $G$. We thus have the following proposition which results from Proposition 3.2 and 3.4 in \cite{eskinoh}.

\begin{proposition}
If $H$ is not contained in a proper parabolic subgroup of $G$, then any sequence of finite unions of closed $H$-orbits of $\Gamma \backslash G$ has no loss of mass.
\end{proposition}

\section{Equidistribution of intersection points}\label{equiintersectionpoints}

We consider as before a sequence $(\mathcal O_n)_{n\in\N}$ of finite unions of closed $H$-orbits of $\Gamma \backslash G$ which is assumed to be equidistributed. In this section, we want to pass from the equidistribution in terms of currents to an equidistribution of the intersection points of $\mathcal O_n$ with a subvariety of $\Gamma \backslash G/K$ of dimension $d_H$.
Though this kind of result can be expected to follow from Theorem \ref{t:EquidistributionCurrentsG/K}, some work is needed to deal with the locus where this intersection is not transverse. This will require in particular a finiteness result for maps defined in an o-minimal structure.

\subsection{Moderate geometry of locally symmetric spaces}
We recall in this section notions from o-minimal geometry in the context of locally symmetric spaces following \cite{bkt} and the structure of definable maps. For a 
general introduction to o-minimal structures, we refer to \cite{vandendries}.

A structure $\cS$ on $\R$ expanding the real field $\R$ is by definition a collection $(\cS_n)_{n\in\N^\times}$
where each $\cS_n$ is a set of subsets of $\R^n$, called the definable sets, which is a Boolean sub-algebra of the subsets of $\R^n$ containing all the algebraic subsets and which satisfy the following properties:
\begin{enumerate}
    \item If $A\in\cS_n$, $B\in\cS_m$, then $A\times B\in \cS_{n+m}$;
    \item if $p:\R^{n+1}\rightarrow \R^n$ is the projection on the first $n$-coordinates, $A\in \cS_{n+1}$, then $p(A)\in \cS_n$.
\end{enumerate}
The structure is called \emph{o-minimal}\footnote{Order-minimal.} if any element of $\cS_1$ is a finite union of points and intervals. 
Given an o-minimal structure, one can define the following notions:
\begin{enumerate}
    \item a map $f:A\rightarrow B$ between two definable sets is definable if its graph $\Gamma_f\subset A\times B$ is definable ;
    \item a $\cS$-definable manifold is a manifold having a \emph{finite} atlas of charts $(\phi_i:U_i\rightarrow \R^n)_{i\in I}$ such that the intersections $\phi_i(U_i\cap U_j)\subset \R^n$ are definable and the change of coordinates maps $\phi_i\circ\phi_j^{-1}:\phi_j(U_i\cap U_j)\rightarrow \phi_i(U_i\cap U_j)$ are $\cS$-definable maps. 
\end{enumerate}

Intuitively, definable manifolds in an o-minimal structure have reasonable geometry locally and at infinity, complex algebraic varieties being an example of definable manifolds.
The first example of an o-minimal structure is the one given by semi-algebraic subsets denoted by $\R_{alg}$. More examples of o-minimal structures have been studied during recent years and the ones relevant to Hodge theory are:
\begin{enumerate}
    \item $\R_{an}$:  the smallest o-minimal structure expanding $\R_{alg}$ and for which restricted analytic functions are definable, see \cite{gabrielov}.
    \item $\R_{exp}$ : the smallest o-minimal structure expanding $\R_{alg}$ and for which the real exponential map is definable, \cite{wilkie}.
    \item $\R_{an,exp}$: the smallest o-minimal structure expanding the two previous structures, \cite{vandendriesmiller,vandendriesmiller1}.
\end{enumerate}

One of the main theorems of \cite[Theorem 1.1]{bkt} asserts that locally symmetric spaces can be endowed with a semi-algebraic structure which is compatible in $\R_{an}$ with the analytic structure on the Borel-Serre compactification, see {\it loc. cit.} for more details. More precisely, let $G$ be a real connected semi-simple group which has a $\Q$-structure, $\Gamma$ a torsion-free arithmetic subgroup of $G$ and $K\subset G$ a compact sub-group.
\begin{theorem}\cite[Theorem 1.1]{bkt}\label{definable}
The quotients $G/K$, $\Gamma\backslash G/K$ have $\R_{alg}$ structures such that
\begin{enumerate}
    \item $G\rightarrow G/K$ is definable in $\R_{alg}$;
    \item there exists a definable fundamental domain $\mathcal{F}\subset G/K$ for the action of $\Gamma$ such that $\mathcal{F}\rightarrow \Gamma\backslash G/K$ is definable in $\R_{alg}$.
\end{enumerate}
Moreover, this structure is functorial in the triple $(G,K,\Gamma)$.
\end{theorem}
One last ingredient we need from o-minimal geometry is the structure of definable maps: a definable map $f:M\rightarrow N$ between definable sets has a \emph{definable trivialization} if there exists a pair $(F,\lambda)$ where $F$ is a definable set and $\lambda : M\rightarrow F$ is a definable map which induces a definable homeomorphism $M\rightarrow F\times N$ compatible with maps to $N$. The following theorem is from \cite[Chapter 9,Th 1.2]{vandendries}.
\begin{theorem}\label{definablemap}
Let $f:M\rightarrow N$ be a continuous definable map between definable sets in an o-minimal structure $\cS$. Then there exists a finite partition $(N_i)_i$ by definable subsets of $N$ such that $f:f^{-1}(N_i)\rightarrow N_i$ is definably trivial.
\end{theorem}
An easy corollary is that the number of connected components of the fibers of a definable map is finite and uniformly bounded.

\subsection{Counting intersection points}\label{countingpoints}
In this section, we explain our conventions for counting points of intersection of varieties mapping to the quotient $\Gamma\backslash G/K$ with closed $H$-orbits. We adopt the same notations as in \Cref{notations} to which we refer. From now on, $\cS$ will be a fixed o-minimal structure which extends $\R_{alg}$.

Let $S\subset \Gamma\backslash G/K$ be a definable real analytic subvariety of dimension $d_H = \dim(G/K)-\dim(H/L)$, and assume that the smooth locus of $S$ is oriented. Let $\mathcal O$ be a finite union of closed $H$-orbits in $\Gamma \backslash G$. Recall that we denote by $\mathcal O/L$ its projection to $\Gamma \backslash G/L$ and by $\pi(\mathcal O/L)$ its projection to $\Gamma \backslash G/K$.

In general, $\pi_{\vert \mathcal O/L}$ is not an embedding, which is why it is more convenient to count intersection at the level of $\Gamma \backslash G/L$, where $\mathcal O/L$ is a finite union of closed leaves of the foliation $\mathcal H$ introduced in Section \ref{ss:EquidistributionCurrents}.

Let $(\Gamma\backslash G/L)_S$ be the preimage of $S$ by the fibration $\pi: \Gamma\backslash G/L\rightarrow \Gamma\backslash G/K$. By assumption on $S$, $(\Gamma \backslash G/L)_S$ and $\mathcal O/L$ have complementary dimension in $\Gamma \backslash G/L$. A point $y\in (\Gamma \backslash G/L)_S \cap \mathcal O/L$ is a \emph{transverse intersection point} if $(\Gamma \backslash G/L)_S$ is smooth at $y$ (equivalently, if $S$ is smooth at $\pi(y)$) and
\begin{equation} \label{directsum}
T_y(\Gamma \backslash G/L)_S \oplus T_y \left(\mathcal O/L\right) = T_y(\Gamma \backslash G/L)~.
\end{equation}

For $y\in (\Gamma \backslash G/L)_S \cap \mathcal O$, we set $\epsilon(y) = 0$ if $y$ is not a transverse intersection point, and $\epsilon(y) =1$ if $y$ is a transverse intersection point such that the direct sum \eqref{directsum} is compatible with orientations, and $\epsilon(y) = -1$ otherwise.


 We have the following distribution on $(\Gamma\backslash G/L)_S$ given by summing over transverse intersection points : 
\begin{align}\label{currentdefinition}
\hat T_\mathcal O^S = \sum_{x\in (\Gamma\backslash G/L)_S\cap\mathcal O/L}\epsilon(x)\delta_x~.
\end{align}

Note that, since $(\Gamma\backslash G/L)_S\cap\mathcal O/L$ is definable, it has a finite number of connected components. Moreover, components of dimension $\geq 1$ consist only of non-transverse intersections points, for which $\epsilon(y) = 0$. Hence $\hat T_\mathcal O^S$ is a finite signed measure. We can finally define:
\begin{definition}
The \emph{transverse intersection measure} between $S$ and $\mathcal O/L$ is the signed measure
\[T^S_{\mathcal O} \overset{\textrm{def}}{=} \pi_* \hat T^S_{\mathcal O}~.\]
\end{definition}


Informally, $T^S_{\mathcal O}$, counts the transverse intersection points of $S$ and $\pi(\mathcal O/L)$, with a multiplicity corresponding to the (signed) number of branches of $\pi(\mathcal O/L)$ meeting $S$ at a smooth point $x$. The asymptotic behavior on $S$ of the distributions $T_{\mathcal O_n}^S$ for an equidistributing sequence $\mathcal O_n$ will be discussed in the next section. In particular, we will prove that they are equidistributed with respect to the form $\pi_*p^*\omega_{G/H}$.

\begin{remarque}
In general, the intersection $S\cap \pi(\mathcal O/L)$ could have zero dimensional components which are not transverse because they are not reduced or are not smooth points of $S$. We do not take them into account in $T_\mathcal O^S$, but we will see in the next section that they are negligible from the equidistribution point of view. 
\end{remarque}
\begin{remarque}
Working in the setting of a moderate geometry allows us to avoid topological pathologies which do not arise in the applications we intend to give and makes it possible at the same time to make statements in a more general setting.
\end{remarque}

\subsection{Equidistribution of intersection points}\label{equiinter}

We keep the notations from the previous section, namely, $G$ is a semi-simple $\Q$-group, $\Gamma\subset G$ is a torsion free arithmetic subgroup and $H\subset G$ is a semi-simple subgroup without compact factor. Let $\mathcal D=G/K$ where $K\subset G$ is a compact subgroup and $L=K\cap H$. Let $d_H$ be as before the codimension of $H/L$ in $\D$. Let $S$ be an analytic subvariety of $\Gamma \backslash \D$ of dimension $d_H$ and consider an equidistributing sequence $\left(\mathcal O_n\right)_{n\in\N}$ of finite unions of closed $H$-orbits in $\Gamma \backslash G$. In this section, we refine \Cref{t:EquidistributionCurrentsG/K} into an equidistribution theorem for the measures $T^S_{\mathcal O_n}$ introduced in the previous paragraph, which will imply \Cref{main}. \\

\begin{theorem}\label{t:interequi}
Let $S$ be an analytic subspace of $\Gamma \backslash \mathcal D$ of dimension $d_H$ and let $(\mathcal{O}_n)_{n\in\N}$ be an equidistributing sequence of finite unions of closed $H$-orbits in $\Gamma\backslash G$. Then 
\[\frac{1}{\Vol(\mathcal{O}_n)}T_{\mathcal{O}_n}^S \underset{n\to +\infty}{\rightharpoonup} \pi_*p^*\omega_{G/H}.\]
\end{theorem}

\begin{proof}
Recall that it is enough to prove the weak convergence on every open subset of a covering of $S$ by \Cref{p:ConvergenceMeasureCovering}. We can thus restrict ourselves to an open relatively compact definable subset $\Omega \subset S$ that lifts to $\mathcal D$. We still denote this lift by $\Omega$. We want to prove the equidistribution of the transverse intersection of $\Omega$ with $\pi(p^{-1}(V_n))$, where $V_n$ is the finite union of discrete $\Gamma$-orbits in $G/H$ such that $\mathcal O_n = \Gamma \backslash \pi^{-1}(V_n)$.

Denote by $(G/L)_\Omega$ the preimage of $\Omega$ by $\pi$, and consider the following open subsets of $(G/L)_\Omega$:
\begin{itemize}
    \item The domain $U_{\textit{smooth}}$ where $(G/L)_\Omega$ is smooth,
    \item the domain $U_{\textit{sub}}$ where $(G/L)_\Omega$ is smooth and $p_{\vert \Omega}$ is a submersion (hence a local diffeomorphism by equality of the dimensions). 
\end{itemize}

By definition, the distribution $\hat T_{\mathcal O_n}^\Omega$ is supported by $U_{\textit{sub}}$. We first prove that $\frac{1}{\Vol(\mathcal{O}_n)} \hat T_{\mathcal O_n}^\Omega$ converges weakly to $p^*\omega_{G/H}$ on $U_{\textit{sub}}$, then extend the weak convergence successively to $U_{\textit{smooth}}$ and $(G/L)_\Omega$ using Proposition \ref{p:ConvergenceMeasureComplement}. After pushing forward by $\pi$, we get the desired result.
\medskip

To prove the weak convergence on $U_{\textit{sub}}$, it is enough to prove it  on every open relatively compact subset $U'$ of $U_{\textit{sub}}$ such that $p_{\vert U'}$ is a diffeomorphism onto its image (since those open sets cover $U_{\textit{sub}}$). Thus, let $f$ be a continuous function with compact support on such $U'$. Let $\epsilon(U')$ be $1$ if $p_{\vert U'}$ preserves orientation and $-1$ otherwise. Then
\begin{eqnarray*}
\frac{1}{\Vol(\mathcal O_n)}\hat T_{\mathcal O_n}^{\Omega}(f) &= & \frac{\epsilon(U')}{\Vol(\mathcal O_n)} \sum_{x\in V_n\cap p(U')} f\circ p^{-1}(x)\\
&\underset{n\to+\infty}{\longrightarrow} & \epsilon(U') \int_{p(U')} f\circ p^{-1}(x) \omega_{G/H}(x)= \int_{U'} f p^*\omega_{G/H}
\end{eqnarray*}
since $V_n$ is equidistributed in $G/H$.

This shows the weak convergence of $\frac{1}{\Vol(\mathcal O_n)}\hat T_{\mathcal O_n}^{\Omega}$ on every $U'$, hence on $U_{\textit{sub}}$.

Let us now extend the weak convergence to $(G/L)_\Omega$. Since the projection map $p_{\vert (G/L)_\Omega}: (G/L)_\Omega \to G/H$ is definable in the o-minimal structure $\R_{an,exp}$, it follows from \Cref{definablemap} that the number of connected components of its fibers is uniformly bounded by some number $N$.

Let first $C$ be a compact subset of $U_{\textit{smooth}}\setminus U_{\textit{sub}}$. Then, by Sard's lemma, $p(C)$ has measure $0$. Hence, for every $\epsilon >0$, there exists an open neighborhood $U'$ of $C$ in $U_{\textit{smooth}}$ such that
\[\left \vert\int_{p(U')} \omega_{G/H}\right \vert \leq \epsilon~.\]
Since for all $x\in p(U')\cap V_n$ the set $p^{-1}(x)\cap U'$ has at most $N$ isolated points, we get that
\[\frac{1}{\Vol(\mathcal O_n)} \vert \hat T_{\mathcal O_n}^\Omega\vert(U') \leq \frac{N}{\Vol(\mathcal O_n)} \vert p(U')\cap V_n\vert ~,\]
hence
\[\frac{1}{\Vol(\mathcal O_n)}\vert\hat  T_{\mathcal O_n}^\Omega\vert(U') \leq \vert\int_{p(U')} \omega_{G/H}\vert \leq N\epsilon~.\]

We conclude that $\frac{1}{\Vol(\mathcal O_n)}\hat T_{\mathcal O_n}^{\Omega}$ converges weakly on $U_{\textit{smooth}}$ by Proposition \ref{p:ConvergenceMeasureComplement}. 

Finally, the complement of $U_{\textit{smooth}}$ in $(G/L)_\Omega$ is the singular locus, which has dimension $< d_H$. Since $p$ is definable, its image has measure $0$, and one can reproduce the previous argument to show that $\frac{1}{\Vol(\mathcal O_n)}\hat T_{\mathcal O_n}^{\Omega}$ converges weakly to $p^*\omega_{G/H}$ on $(G/L)_\Omega$.
\end{proof}

In the previous theorem, we chose to restrict to $\dim(S)= d_H$ for simplicity. We indicate briefly without proof how these statements  should be adapted when $\dim(S)> d_H$: we define the \emph{transverse intersection current} $T_S^{\mathcal O_n}$ as the integration current over the transverse intersection locus of $S\cap \mathcal O_n$ with the sign normalizations as in \Cref{countingpoints}. Intersecting further with submanifolds of dimension $d_H$ and applying our theorem, one would obtain the convergence of these intersection currents (after normalization) to the current ${\pi_*p^*\omega_{G/H}}_{\vert S}$. We thus get the following theorem. 
\begin{theorem}\label{currentsremark}
Let $S$ be an analytic subspace of dimension $d\geq d_H$ and let $(\mathcal{O}_n)_{n\in\N}$ be an equidistributing sequence of finite unions of closed $H$-orbits in $\Gamma\backslash G$. Then for every $\beta\in \Omega_{c}^{d-d_H}(S)$, we have  
\[\frac{1}{\Vol(\mathcal{O}_n)}T_{\mathcal{O}_n}^S(\beta) \underset{n\to +\infty}{\rightarrow} \pi_*p^*\omega_{G/H}\wedge\beta.\]
\end{theorem}

\section{The pull-push form}\label{s:pullpush}

In the equidistribution theorems \ref{t:EquidistributionCurrentsG/K} and \ref{cohomologyequidistribution}, we see appearing the \emph{pull-push form} $\pi_*p^*\omega_{G/H}$. This form was already studied by the second author with entirely different motivations \cite{tholozan2016volume}, while the first author proved, in the particular case of $G= \SO(2,q)$, $H= \SO(2,q-1)$ and $K= \mathrm S(O(2)\times O(q))$, that this form is the $G$-invariant K\"ahler form on $G/K$.

Here we introduce some tools to better characterize this form, most of which were already presented in \cite{tholozan2016volume}.

\subsection{Push-forward of forms}\label{pushforward}

Let us  first recall without any proof the construction of the push-forward of a form $\alpha$ under a smooth fibration $\pi: M\to B$ with compact oriented fibers. This construction can be found for instance in \cite[p.37]{bottTu}.

Let $r$ be the dimension of the fibers of $\pi$. Let $x$ be a point in $B$ and denote by $F$ the fiber of $\pi$ over $x$. Choose $\omega$ a volume form on $F$ compatible with the orientation of the fiber, and let
\[X\in \Lambda^r(TF)\subset \Lambda^r TM\]
be the multivector such that $\omega(X) = 1$.

Let now $\alpha$ be a smooth $p$-form on $M$. The contraction $\iota_X\alpha$ is a $(p-r)$ form with kernel $TF$, which can thus be seen as a section of $\Lambda^{p-r}(NF^\vee)$, where $NF = TM/TF$ is the normal bundle to $F$ and $NF^\vee = \{\varphi \in T^*M \mid \varphi_{\vert TF} = 0\}$ is its dual.

Finally, the differential of $\pi$ along the fiber $F$ induces an isomorphism $NF \simeq F\times T_xB$ and therefore $\Lambda^{p-r}(NF^\vee) = \Lambda^{p-r}(T_xB^\vee)$. With these identifications, we can now define:
\begin{definition} \label{def: Push-forward}
The push-forward of the form $\alpha$ is the $(p-r)$-form on $B$ given at $x$ by
\[(\pi_*\alpha)_x = \int_{y\in F} \iota_X(\alpha)_y  \omega~.\]
\end{definition}

Let $p+q$ be the dimension of $M$. Then we have
\begin{proposition} \label{p:PropertyPushForward}
The form $\pi_*\alpha$ is the unique $(p-r)$ form on $B$ such that for any $q$-form $\beta$ on $B$ with compact support, 
\[\int_B \pi_*\alpha \wedge \beta = \int_M \alpha \wedge \pi^*\beta~.\]

Moreover, the push-forward operation commutes with the exterior derivative. In particular, the push-forward of a closed form is closed.
\end{proposition}

\subsection{A formula for the pull-push form} \label{pullpushsection}

Let us now apply the previous general considerations in order to give a formula for the pull-push form $\pi_*p^*\omega_{G/H}$ at the Lie algebra level.

Let $G/H$ be a $G$-homogeneous space and $o$ denote the base point of $G/H$ (that is, the projection of the unit element of $G$). Then the tangent space $T_o G/H$ with the induced linear action of $H$ identifies canonically with $\fg/\fh$ endowed with the adjoint action of $H$. This identification induces an isomorphism between the differential algebra $\Omega^\bullet(G/H, \C)^G$ of smooth $G$-invariant forms on $G/H$ with complex coefficients and the differential algebra $\Lambda^\bullet((\fg/\fh)_\C^\vee)^H$ of $H$-invariant exterior forms on $\fg/\fh$, with derivative given by
\[\mathrm{d} \alpha(x_1,\ldots, x_{k+1}) = \sum_{i<j} (-1)^{i+j}\alpha([x_i,x_j],x_1,\ldots, \hat x_i, \ldots, \hat x_j, \ldots, x_{k+1})~,\]
for $\alpha\in \Lambda^k((\fg/\fh)_\C^\vee)^H$ and $x_1,\ldots,x_{k+1}\in \fg/\fh$.

Let us now come back to the setting of the previous section, where $H$ is a semi-simple subgroup of $G$, $K$ a compact subgroup of $G$ and $L=G/H$. At the Lie algebra level, the pull-back homomorphism $p^*$ (resp. $\pi^*$) identifies exterior forms on $\fg/\fh$ (resp. $\fg/\fk$) to those exterior forms on $\fg/\fl$ having $\fh/\fl$ (resp. $\fk/\fl$) in their kernel. 

Reinterpreting Definition \ref{def: Push-forward} for $G$-invariant forms on $G/L$, we get:
\begin{proposition} \label{p:PushForwardInvariantForm}
Let $\alpha$ be a $G$-invariant form on $G/L$. Then the form $\pi_*\alpha$ on $G/K$ is $G$-invariant and corresponds on $\fg/\fk$ to the $\Ad_K$-invariant form
\[\int_{k\in K/L} {\Ad_k}^* (\iota_u\alpha)  \omega_{K/L}~,\]
where $u \in \Lambda^{\max} (\fk/\fl)$ is such that $\omega_{K/L}(u) = 1$.
\end{proposition}

\begin{remarque}
Since $\fk/\fl$ and $\fh/\fl$ are both $\Ad_L$-invariant, the form $\iota_u \alpha$ is $L$-invariant and its kernel contains $\fk/\fl$. Therefore, its pull-back by $\Ad_k$ only depends on the class of $k$ in $K/L$ so that the integral makes sense. Moreover, the resulting form is obviously $K$-invariant and has $\fk/\fl$ in its kernel, so that is does identify with a $K$-invariant form on $\fg/\fk$.
\end{remarque}

\begin{corollaire} \label{c:IntegralFormulaPushPull}
Let $(\alpha_1\ldots, \alpha_{p_r})$ be an oriented orthonormal basis of $\{\phi \in (\fg/\fk)^\vee \mid \phi_{\vert\fh/\fl} = 0\}$. Then there exists a positive constant $\lambda$ such that 
\[\pi_*p^*\omega_{G/H} = \lambda \int_{k\in K/L} {\Ad_k^*}(\alpha_1\wedge \ldots \wedge \alpha_{p-r})  \omega_{K/L}(k)~.\]
If moreover $\fk/\fl$ is orthogonal to $\fh/\fl$ in $\fg/\fl$, then $\lambda = 1$.
\end{corollaire}

\begin{proof}[Proof of Proposition \ref{p:PushForwardInvariantForm}]

The $G$-invariance of $\pi_*\alpha$ easily follows from Proposition \ref{p:PropertyPushForward}, so we only need to compute $\pi_*\alpha$ at the base point of $G/K$.

Let $o$ denote the basepoint of $G/L$ and $v_1,\ldots, v_{p-r}$ be $p-r$ vectors in $T_o G/L = \fg/\fl$. Let $k$ be an element of $K$. At the point $k\cdot o$, we have
\[\iota_u\alpha(k\cdot v_1,\ldots, k\cdot v_{p-r}) = \iota_u \alpha(v_1,\ldots, v_{p-r})\] by left invariance of $\iota_u \alpha$.

Now, the differential of $\pi$ maps a vector $k\cdot v \in T_{k\cdot o}G/L$ to the vector $\Ad_{k}(v) \in T_{\pi(o)} G/K = \fg/\fk$. Applying \cref{def: Push-forward} to $\pi_*\alpha$ gives the required formula.
\end{proof}

\begin{proof}[Proof of Corollary \ref{c:IntegralFormulaPushPull}]
Complete $(\alpha_1,\ldots, \alpha_{p-r})$ into an orthonormal basis $(\beta_1,\ldots, \beta_r, \alpha_1, \ldots ,\alpha_{p-r})$ of $\{\phi\in (\fg/\fl)^* \mid \phi_{\vert \fh/\fl}=0\}$. We then have
\[p^*\omega_{G/H} = \bigwedge_{i=1}^r \beta_i \wedge \bigwedge_{i=1}^{p-r} \alpha_i~,\]
hence
\[\iota_up^*\omega_{G/H} = \lambda \bigwedge_{i=1}^{p-r} \alpha_i\]
where $\lambda = \bigwedge_{i=1}^r \beta_i(u)$, since all the $\alpha_i$ vanish on $\fk/\fl$.

If furthermore $\fk/\fl$ is orthogonal to $\fh/\fl$, then $(\beta_1, \ldots, \beta_r)$ can be chosen as an orthonormal basis of $(\fk/\fl)^\vee$, so that $\lambda =1$.

Proposition \ref{p:PushForwardInvariantForm} now concludes the proof.
\end{proof}

In practice, using this integral formula to compute the pull-push form quickly leads to rather involved computations. In \cite{tholozan2016volume}, the second author used this formula to give a sufficient vanishing criterion: if there exists $k\in K$ such that $\Ad_k$ preserves $\fh/\fl$ and reverses its orientation, then ${\Ad_k}_* (\iota_u \alpha) = - \iota_u \alpha$, and Proposition \ref{p:PushForwardInvariantForm} shows that $\pi_*p^*\omega_{G/H}$ vanishes.

In terms of our equidistribution result, this vanishing means that the positive and negative part of the intersection measure cancel out asymptotically, because a ``random'' translate of $H/L$ in $G/K$ can intersect a submanifold $S$ with two opposite equiprobable orientations.

In contrast, we can prove the following non-vanishing criterion:
\begin{corollaire}\label{complexnonvanishing}
If $G/K$ has a $G$-invariant complex structure such that $H/L$ is a complex submanifold, then $\pi_*p^*\omega_{G/H}$ does not vanish.
\end{corollaire}

\begin{proof}
Let $\alpha_1,\ldots, \alpha_{\frac{p-r}{2}}$ be a complex basis of $\{\phi\in (\fg/\fk)^\vee \mid \phi_{\vert \fh/\fl}=0\}$. Then 
$\iota_u p^*\omega_{G/H}$ is proportional to 
\[\bigwedge_{i=1}^{\frac{p-r}{2}} \alpha_i\wedge \bar\alpha_i~.\]
Which is non-negative on every complex subspace of dimension $\frac{p-r}{2}$. By invariance of the complex structure, the same holds for \[{\Ad_k}^*\left(\bigwedge_{i=1}^{\frac{p-r}{2}} \alpha_i\wedge \bar\alpha_i\right)~.\]
Finally, $\bigwedge_{i=1}^{\frac{p-r}{2}} \alpha_i\wedge \bar\alpha_i$ is positive on a complex complement of $\fh/\fl$, and so is
\[\int_{k\in K/L} {\Ad_k}^*\left(\bigwedge_{i=1}^{\frac{p-r}{2}} \alpha_i\wedge \bar\alpha_i\right)  \omega_{K/L}~.\]
\end{proof}

\subsection{Compact duality}\label{compactdualitysection}

In \cite{tholozan2016volume}, the second author investigated further the pull-push form in the case where $G/K$ and $H/L$ are symmetric spaces. There, he proved that $\pi_*p^*\omega_{G/H}$ is in some sense ``Poincaré dual'' to the inclusion $H/L\hookrightarrow G/K$, a statement which is made precise by passing to the compact dual of the symmetric space. Since cohomology classes of compact symmetric spaces are represented by a unique invariant form, this argument completely characterizes the pull-push forms and provides an efficient way to compute it in practice. The goal of this section is to extend these results to more general compact subgroups $K\subset G$. \\

Recall that a \emph{Cartan involution} of $G$ is an involutive automorphism whose fixed subgroup is a maximal compact subgroup. All the Cartan involutions of $G$ are conjugated, and every compact subgroup of $G$ is fixed by a Cartan involution.

In this section, we make the assumption (verified in all the examples we consider in this paper) that there exists a Cartan involution $\theta$ of $G$ fixing $K$ and preserving $H$. We denote by $G^\theta$ and $H^\theta$ the compact subgroups of $G$ and $H$ fixed by $\theta$. The Lie algebras $\fg$ and $\fh$ decompose respectively as
\[\fg = \fg^\theta \oplus {\fg^\theta}^\perp\]
and 
\[\fh = \fh^\theta \oplus {\fh^\theta}^\perp~.\]
We now introduce the dual Lie algebras
\[\fg_U = \fg^\theta \oplus i {\fg^\theta}^\perp \subset \fg_\C\]
and 
\[\fh_U = \fh^\theta \oplus i {\fh^\theta}^\perp \subset \fh_\C~.\]
These are respectively the Lie algebras of compact real forms $G_U$ and $H_U$ of $G_\C$ and $H_\C$, called the \emph{compact duals} of $G$ and $H$.

From $G_U$ and $H_U$, one can define compact duals to the homogeneous spaces $G/H$, $G/L$ and $G/K$, respectively given by $G_U/H_U$, $G_U/L$ and $G_U/K$. The various inclusions between these groups give the following commutative diagram:

\[
\xymatrix{
G/L \ar[r]^p \ar@{^{(}->}[rd] \ar[d]^\pi & G/H \ar@{^{(}->}[rd] & & \\
G/K \ar@{^{(}->}[rd] & G_\C /L_\C \ar[d] \ar[r] & G_\C / H_\C  & \\
& G_\C / K_\C & G_U/L \ar@{_{(}->}[lu] \ar[r]^{p_U} \ar[d]^{\pi_U} & G_U/H_U \ar@{_{(}->}[lu] \\
& & G_U/K \ar@{_{(}->}[lu]&
}
\]

Now, the inclusion of $G/H$ into $G_\C/H_\C$ induces an isomorphism of differential algebras
\[\Omega^\bullet(G/H,\C)^G = \Lambda^\bullet((\fg/\fh)^\vee)^H\otimes_\R \C \simeq \Lambda^\bullet_\C((\fg_\C/\fh_\C)^\vee)^{H_\C} = \Omega_\C^\bullet(G_\C/H_\C)^{G_\C}~,\]
where $\Lambda^\bullet_\C$ and $\Omega^\bullet_\C$ denote the complex of $\C$-multilinear forms. The same holds for all the inclusions in the above diagram. In other words, the differential algebra of real invariant forms on a homogeneous space and its compact dual are two distinct real forms of the same complex differential algebra.

\begin{proposition}\label{p:CommutativityPushPullCompactDual}
We have the following commutative diagram of differential complexes:
\[
\xymatrix{
\Omega^\bullet(G/L,\C)^G \ar[d]^{\pi_*} \ar@{<->}[dr]^{\sim} & \ar[l]^{p^*}\Omega^\bullet(G/H,\C)^G \ar@{<->}[dr]^{\sim} & \\
\Omega^{\bullet-\dim(K/L)}(G/K,\C)^G \ar@{<->}[dr]^{\sim} &\Omega^\bullet(G_U/L,\C)^{G_U} \ar[d]^{{\pi_U}_*} & \ar[l]^{p_U^*} \Omega^\bullet(G_U/H_U,\C)^{G_U} \\
& \Omega^{\bullet-\dim(K/L)}(G_U/K,\C)^{G_U} &
}
\]
\end{proposition} 

\begin{proof} 
The only non-trivial point is that $\pi_*$ and ${\pi_U}_*$ are identified as maps from $\Omega_\C^\bullet(G_\C/L_\C)^{G_\C}$ to $\Omega_\C^\bullet(G_\C/K_\C)^{G_\C}$. But this readily follows from \Cref{p:PushForwardInvariantForm} since, at the Lie algebra level, both maps are given by the contraction with $u$ followed by averaging under the adjoint action of $K$.
\end{proof}

Now, $\omega_{G/H}$ and $\omega_{G_U/H_U}$ are both generators of $\Lambda^{\max}_\C(\fg_\C/\fh_\C)$, so they are complex multiples one of the other. More precisely, we can identify 
\[\fg/\fh = \fg^\theta /\fh^\theta \oplus {\fg^\theta}^\perp/{\fh^\theta}^\perp\]
with 
\[\fg_U/\fh_U = \fg^\theta /\fh^\theta \oplus i{\fg^\theta}^\perp/{\fh^\theta}^\perp\]
via the morphism
\[\phi:(u,v)\mapsto (u,iv)\]
and normalize $\omega_{G_U/H_U}$ so that
\[\phi^*\omega_{G_U/H_U} = \omega_{G/H}~.\]
With this normalization, we have the following equality in $\Lambda^{\max}_\C(\fg_\C/\fh_\C)$
\[\omega_{G/H} = i^{\dim({\fg^\theta}^\perp/{\fh^\theta}^\perp)} \omega_{G_U/H_U}~.\]
Finally, applying Proposition \ref{p:CommutativityPushPullCompactDual}, we conclude:

\begin{corollaire}\label{equalityofforms}
The following equality holds
in $\Lambda^\bullet_\C(G_\C/K_\C)^{K_\C}$:
\[\pi_* p^*\omega_{G/H} = i^{\dim({\fg^\theta}^\perp/{\fh^\theta}^\perp)} {\pi_U}_*p_U^*\omega_{G_U/H_U}~.\]
\end{corollaire}

What we gained by switching to the compact dual space $G_U/K$ is that we can now talk about the cohomology class of the pull-push form. The following theorem was proven in \cite{tholozan2016volume} under the assumption that $K=G^\theta$, but the proof easily adapts to our more general context:
\begin{lemme}\label{dualofahomologyclass}
The de Rham cohomology class of the pull-push form $\frac{1}{\Vol(G_U/H_U)}{\pi_U}_* p_U^*\omega_{G_U/H_U}$ is Poincaré dual to the homology class of $H_U/L \subset G_U/K$; that is, for every closed form $\beta$ on $G_U/K$ of degree $\dim(H_U/L)$, we have
\[\int_{H_U/L} \beta = \frac{1}{\Vol(G_U/H_U)} \int_{G_U/K} {\pi_U}_* p_U^*(\omega_{G_U/H_U})\wedge \beta~.\]
\end{lemme}

\begin{proof}
This is essentially formal. Denote respectively by $\iota_1$ and $\iota_2$ the inclusions of $H_U/L$ in $G_U/L$ and $G_U/K$, so that we have $\iota_2 = \pi_U\circ \iota_1$ and $\iota_1(H_U/L) = p_U^{-1}(o)$, where $o$ is the basepoint of $G_U/H_U$.

Now, the form $\frac{1}{\Vol(G_U/H_U)}\omega_{G_U/H_U}$ is Poincaré dual to $[o]$ in $\H_{dR}^\bullet(G_U/H_U,\R)$, so its pull-back under $p_U$ is Poincaré dual to $[{p}_{U}^{-1}(o)] = \iota_{1*}[H_U/L]$. Finally, let $\beta$ be a closed form of degree $\dim(H_U/L)$ on $G_U/K$. We then have

\begin{eqnarray*}
\int_{H_U/L} \iota_2^* \beta &=& \int_{\iota_1(H_U/L)} \pi_U^*\beta \\
&=& \int_{G_U/L}p_U^*\omega_{G_U/H_U}\wedge \pi_U^*\beta\\
&=& \int_{G_U/K} {\pi_U}_*p_U^* \omega_{G_U/H_U} \wedge \beta~.
\end{eqnarray*}
\end{proof}

\subsubsection{The symmetric case}\label{thesymmetriccase}

When $K= G^\theta$ is a maximal compact subgroup of $G$, a theorem of Élie Cartan \cite{cartan} states that all $G$-invariant forms on $G/K$ are closed. Hence the exterior derivative is trivial on $\Omega^\bullet(G/K)^G$ and we have isomorphisms:
\[\Omega^\bullet(G/K,\C)^G \simeq \Omega^\bullet(G_U/K, \C)^{G_U} \simeq  H_{dR}^\bullet(G_U/K,\C)~.\]
In other words, a $G$-invariant form on $G/K$ is completely characterized by the cohomology class of the corresponding form on $G_U/K$.\\

This property does not hold for more general homogeneous spaces. In the next section we will see however that it remains true on Mumford--Tate domains when restricting to a variation of Hodge structure.

\section{Invariant forms on period domains}\label{definitionvhs}

In this section, we introduce Mumford--Tate domains, Hodge loci and invariant forms on period domains. Then we relate pull-push forms to Chern classes of Hodge bundles. Finally, we compute it in various cases of interest.

\subsection{Variations of Hodge structure and their period domains}

Let us first recall some definitions of Hodge theory, merely to fix notations.

Let $V$ be a free $\Z$-module of finite rank $d\in \N$ endowed with a bilinear form $B: V\times V \to \Z$. Given a field $\mathbb K$, we write $V_\mathbb{K}= V\otimes_\Z\mathbb K$ and still denote by $B$ the natural $\mathbb K$-bilinear extension of $B$ to $V_\mathbb{K}$.

A \emph{Hodge structure} of weight $k$ on $V$ polarized by $B$ is the data of a filtration of complex vector spaces
\[0\subseteq F^{k}\subseteq\cdots\subseteq F^0= \V_\C\]
such that for all $0\leq p\leq k$,
\begin{enumerate}
    \item $V_\C=F^{p}\oplus \overline{F}^{k-p+1}$,
    \item $B(u, v) = 0$ for all $(u,v)\in F^p \times F^{k-p+1}$,
    \item $i^{p-q}B(v, \bar v)>0$ for all $v\in (F^p \cap \overline{F}^{q})\backslash\{0\}$ with $p+q=k$.
\end{enumerate}

\begin{remarque}
The existence of a Hodge structure of weight $k$ implies that $B$ is non-degenerate and antisymmetric for odd $k$ or symmetric for even $k$.
\end{remarque}

For $p+q = k$, define $V^{p,q} = F^p \cap \overline F^{q}$. Then $V^{p,q}$ is a complement of $F^{p+1}$ in $F^p$. In particular, we have a decomposition
\[V_\C = \bigoplus_{p+q=k} V^{p,q},\,\textrm{with}\quad \overline{V}^{p,q}=V^{q,p}~.\]
The \emph{Hodge numbers} of the Hodge structure are the numbers $h^{p,q} \equaldef  \dim_\C(V^{p,q})$.

Let $\mathbb{S}=\mathrm{Res}_{\C/\R}\mathbb{G}_m$ denote the \emph{Deligne torus}, i.e., the restriction of scalars of the multiplicative group $\mathbb G_m$ from $\C$ to $\R$. Then $\mathbb{S}(\R)=\C^\times$ seen as an algebraic group over $\R$. Every Hodge structure on $(V_\Z, B)$ induces a representation
$\varphi: \mathbb{S}(\R) \to \GL(V_\R)$ given by
\[z\cdot u = z^{-p} \bar z^{-q} u\]
on $u\in V^{p,q}$.

Let $k\in \Z$ and $\tuple h = (h_{p,q})_{p+q = k} \in \N^{k+1}$ be such that $h_{p,q} = h_{q,p}$ and $\sum_{p=0}^k h^{p,q} = d$. The \emph{period domain} of Hodge structure of weight $k$ and Hodge numbers $(h_{p,q})_{p+q = k}$ is the set $\D$ of all filtrations $(F^{p})_{0\leq p\leq k}$ which define a Hodge structure of weight $k$ and Hodge numbers $h^{p,q}$ on $(V_\Z,B)$. It is a complex manifold homogeneous under the action of the group $\Aut_\R(B)$, and the stabilizer of a point is a compact subgroup of $\Aut_\R(B)$.

The period domain $\D$ is an open subset of the \emph{compactified period domain}
$\widehat{\D}$ of complex flags $0\subseteq F^{k}\subseteq\cdots\subseteq F^0 = V_\C$
such that $F^{k-p} = {F^p}^\bot$ and $\dim_\C(F^p/F^{p+1}) = h^{p,k-p}$ for all $p$. The compactified period domain is a \emph{flag variety} of the group $\Aut_\C(B)$ (i.e., a quotient of $\Aut_\C(B)$ by a parabolic subgroup).

Let $\cU$ denote the trivial complex vector bundle $\D \times V_\C$ equipped with the action of $\Aut_\R(B)$ given by the tautological linear action in the fibers. By construction, this bundle admits a $\Aut_\R(B)$-invariant real structure and a complex bilinear pairing $B$ as well as a \emph{universal Hodge decomposition}, i.e., a smooth decomposition as a direct sum of $G_\R$-equivariant complex vector bundles $\mathcal U^{p,q}$ such that, at a point $x$, the induced decomposition of $\cU_x = V\otimes_\Z \C$ is the Hodge decomposition associated to $x$. 

\medskip

Let now $X$ be a complex analytic variety. A \emph{(polarized) variation of Hodge structure of weight $k$} over $X$ is the data of: 
\begin{itemize}
    \item A local system $\V_\Z$ of free $\Z$-modules of finite rank $d$ with a  flat bilinear pairing $B:\mathbb{V}_\Z\otimes \mathbb{V}_\Z\rightarrow \underline{\Z}_X$,
    \item A decreasing filtration $\cF^{\bullet}\cV$ on $\mathcal{V}=\mathbb{V}_\Z\otimes_{\underline{\Z}_X} \mathcal{O}_X$ by holomorphic sub-vector bundles $0\subseteq\cF^{k}\cV\subseteq\cdots\subseteq\cF^0\cV=\cV$,
    \end{itemize}
Which satisfy the following conditions:
\begin{itemize}
    \item {\bf Hodge property:} For every $x\in X$, the flag $0\subseteq\cF^{k}\cV_x\subseteq\cdots\subseteq\cF^0\cV_x\equaldef \cV_x$ is a Hodge structure on ${\V_{\Z,x}}$ ;
    
    \item {\bf Griffiths' transversality:} The flat connection $\nabla$ associated on $\V_\Z\otimes \mathcal O_X$ satisfies
    \[\nabla(\cF^{p}\cV)\subseteq \cF^{p-1} \cV\otimes \Omega^{1}_X \textrm{ for } 0\leq p\leq k~.\]
\end{itemize}

Let $\{\V_\Z,\cF^\bullet\cV,B\}$ be a variation of Hodge structure of weight $k$ over $X$. Its \emph{Hodge decomposition} is the ($C^{\infty}$) decomposition \[\cV = \bigoplus_{p+q= k} \cV^{p,q}~,\]
where $\cV^{p,q} = \cF^p\cV \cap \overline{\cF}^{k-p} \cV$, and its \emph{Hodge numbers} are the
\[h^{p,q} = \dim_\C(\mathcal V^{p,q})~, p+q = k~.\]

Let now $\pi:\tilde X\to X$ be the universal cover of $X$ and $x$ and arbitrary point in $\tilde X$. The local system $\pi^*\V_\Z$ is trivial, and one obtains a map

\[\tilde f: \tilde X \to \D\]
such that $\tilde f(y)$ is the Hodge structure $\cF^\bullet \cV_y$ on $(\pi^*{\V_{\Z}})_y = {\V_{\Z,x}}$. This map is equivariant with respect to the monodromy $\rho:\pi_1(X) \to G_\Z = \Aut({\V_{\Z,x}})$ of the local system and thus factors to a map
\[f:X \to G_\Z \backslash \D\]
called the \emph{period map} of the variation of Hodge structure. There are canonical isomorphisms
\[\mathcal V^{p,q} \simeq f^*\mathcal U^{p,q}~.\]

In terms of the period map, Griffiths' transversality condition admits the following interpretation.
Let $x$ be a point in $\D$ and let $\varphi: \mathbb{S} \to G_\R$ be the associated representation of the Deligne torus. Then the Lie algebra $\fg_\C$ decomposes under the adjoint action as

\[\fg_\C=\bigoplus_p\fg^{p,-p},\]
where \[\fg^{p,-p}=\{\xi\in\fg,\xi\cdot V^{r,s}\subset V^{r+p,s-p}\}.\]

The subalgebra $\fg^{0,0}$ is the Lie algebra of the stabilizer of $x$, and its complement identifies with the complexified tangent space to $\D$ at $x$. The eigenspace of the complex structure on $T_x \D$ for $i$ is the subspace $\bigoplus_{p<0} \fg^{p,-p}$.

The subspace $\fg^{1,-1}\oplus \fg^{-1,1}$ is the complexification of a well defined subspace $W_x \subset T_x \D$. This defines a holomorphic $G_\R$-invariant distribution of $T_x \D$ called the \emph{Griffiths' distribution}. Now, Griffiths' transversality condition states precisely that the period map is tangent to the Griffiths' distribution.

\subsection{Hodge loci and transversality}

Let $(V, B)$ be a lattice with an integral bilinear pairing. A Hodge structure on $V$ induces a Hodge structure on $T^{k,l}V \equaldef  V^{\otimes k} \otimes {V^\vee}^{\otimes l}$ for all $k,l$, whose Hodge decomposition is given by the diagonalisation of the induced representation $\varphi: \mathbb{S}(\R) \to \End(T^{k,l}V\otimes_\Z \C)$ of the Deligne torus. Let $\mathbb{U}^1\subset \mathbb{S}(\R)$ denote the unit circle.

\begin{definition}\label{mumfordtate}
The \emph{Mumford--Tate group} $MT_{\varphi}$ of $(V,B)$ is the smallest $\Q$-algebraic subgroup of $\mathrm{GL}(V_\R)$ which contains $\varphi(\C^\times)$. The \emph{special Mumford--Tate group} is the smallest $\Q$-algebraic subgroup $sMT_{\varphi}$ which contains $\varphi(\mathbb{U}^1)$.
\end{definition}

The algebra of \emph{Hodge classes} is the bigraded $\Z$-subalgebra $\Hodge^{\bullet, \bullet}(\varphi) \subset T^{\bullet,\bullet}V$ fixed by $\varphi(\mathbb{U}^1)$.

Let now $v$ be a vector in $T^{\bullet, \bullet}V$. The \emph{Hodge domain} of $v$ is the set of variations of Hodge structure $\varphi$ on $V$ such that $\Hodge^{\bullet, \bullet}(\varphi)$ contains $v$. The connected components of the Hodge domain of $v$ are homogeneous under the stabilizer of $v$ in $G_\R$. They are called \emph{Mumford--Tate domains}, and the stabilizer of such components are Mumford--Tate groups.

\begin{remarque}
If $\Hodge^{\bullet, \bullet}(\varphi)$ contains a set $A$, then it contains the subalgebra spanned by $A$. Conversely, for every bigraded subalgebra $H^{\bullet,\bullet}$ of $T^{\bullet,\bullet}V$, there exists $v\in H^{\bullet,\bullet}$ such that 
\[v\in \Hodge^{\bullet, \bullet}(\varphi) \Longleftrightarrow H^{\bullet,\bullet} \subseteq \Hodge^{\bullet, \bullet}(\varphi)~.\]
In particular, intersections of Hodge or Mumford--Tate domains are again Hodge and Mumford--Tate domains.
\end{remarque}

Since Mumford--Tate groups are defined over $\Q$, the projection of a Mumford--Tate domain $H_\R/L_\R$ to $G_\Z \backslash G_\R/K_\R$ factors to a proper immersion of $H_\Z \backslash H_\R/L_\R$.

Let now $X$ be a connected analytic variety equipped with a variation of Hodge structure $(\V_\Z,B,\cF^\bullet \cV)$ of weight $k$ and Hodge numbers $(h^{p,q})_{p+q=k}$. We assume that the period map of $X$ is  generically immersive.

Let $G_\R/K_\R$ be a Mumford--Tate domain containing $\tilde X$. Then the monodromy representation takes values in $G_\Z$ and at every point $y\in \tilde X$ the algebra $\Hodge^{\bullet, \bullet}(\varphi)$ of Hodge classes at $y$ contains the subalgebra $H^{\bullet, \bullet}$ fixed by $G_\R$.

The variation of Hodge structure $X$ is called \emph{Hodge generic} in $G_\R/K_\R$ if there is no proper Hodge subdomain of $G_\R/K_\R$ containing $\tilde X$. In that case, at a generic point of $X$, the algebra of Hodge classes is exactly $(T^{\bullet,\bullet} \V)^{G_\R}$ and the Mumford--Tate group is a rational form of $G$. We define:

\begin{definition}
Let $G_\R/K_\R$ be the smallest Mumford--Tate domain containing $X$. The \emph{Hodge locus} of $\tilde X$ is the set of points at which the algebra of Hodge classes contains strictly $(T^{\bullet,\bullet} \V)^{G_\R}$. The Hodge locus of $X$ is its projection under the covering map.
\end{definition}

The Hodge locus of $X$ is the intersection of $X$ with the countable union of all the projections modulo $G_\Z$ of the Mumford--Tate subdomains of $G_\R/K_\R$. To be more precise, let $G_\R/K_\R$ be any Mumford--Tate domain containing $\tilde{X}$ and let $H$ be an algebraic subgroup of $G$ defined over $\Q$. We define

\begin{definition}
The \emph{Hodge locus of type $H$} is the set of points in $\tilde{X}$ whose Mumford--Tate group is conjugated over $\R$ to a subgroup of $H$. The Hodge locus of type $H$ in $X$ is its projection by the covering map.
\end{definition}

The Hodge locus of type $H$ is the intersection of $\tilde X$ with the union of Mumford--Tate domains $\bigcup_{g\in G_\R} g H_\R/L_\R$, for all $g\in G_\R$ such that $g H_\R g^{-1}$ is $\Q$-subgroup. This leads to the following definition:

\begin{definition}
The \emph{transverse Hodge locus} of type $H$ is the set of smooth points of $X$ for which there exists $g\in G_\R$ such that $gHg^{-1}$ is a $\Q$-group and $\tilde{X}$ and $g H_\R/L_\R$ intersect transversally at $x$.

If $X$ is Hodge generic in $\mathcal D$, the transverse Hodge locus (of type $H$) is called the \emph{typical Hodge locus} (of type $H$).
\end{definition}

Since Hodge loci are intersections of $X$ with locally homogeneous sub-spaces of $G_\Z \backslash \mathcal D$, we can hope to apply our equidistribution result in this setting. However, in order for it to be effective, one needs a \emph{generic transversality} property between $X$ and $H/L$:

\begin{definition}
We say that $X$ is \emph{generically transverse} to $H$-orbits at a smooth point $x$ if there exists $g\in G_\R$ such that $gH_\R/L_\R$ and $\tilde X$ intersect transversally at (some lift of) $x$.\\

We say that $X$ is \emph{generically transverse} to $H$-orbits if there exists a smooth point at which it is generically transverse.
\end{definition}

\begin{remarque}\label{r:density}
If $X$ is generically transverse to $H$-orbits, then the set of points $x$ at which it is generically transverse is an open and dense analytic subset of $X$.
\end{remarque}

\begin{proposition}\label{generictransversality}
Let $x$ be a point in $X$. Then $X$ is generically transverse to $H$-orbits at $x$ if and only if the pull-push form $\pi_*p^*\omega_{G/H}$ is non-zero at $x$.
\end{proposition}

As the consequence, we get the following density criterion for the transverse Hodge locus of type $H$.

\begin{theorem} \label{t:densitycriterion}
The following propositions are equivalent:
\begin{itemize}
    \item[$(i)$] the transverse Hodge locus of type $H$ is non-empty,
    \item[$(ii)$] the transverse Hodge locus of type $H$ is analytically dense in $X$,
    \item[$(iii)$] $X$ is generically transverse to $H$-orbits,
    \item[$(iv)$] the pull-push form $\pi_*p^*\omega_{G/H}$ is not identically $0$ on $X$.
\end{itemize}
\end{theorem}

\begin{proof}[Proof of Proposition \ref{generictransversality}]
Let $d$ be half of the degree of the form $\pi_*p^*\omega_{G/H}$.

Assume that $X$ is generically transverse to $H$-orbits at $x$, and let $g\in G$ be such that $x\in gH/L$ and \[T_xX+T_x(gH/L)=T_xG/K.\] 
Let $u$ be a multivector as in \Cref{p:PushForwardInvariantForm}. Since $T_x(gH/L)_\C$ is in the kernel of $\iota_u p^*\omega_{G/H}$, there exists holomorphic vector fields $X_1,\ldots,X_d$ on $X$ defined on a neighborhood of $x$ such that \[\iota_u\p^*\omega_{G/H}(X_1,\overline{X}_1\ldots, X_d,\overline{X}_d)>0.\]
For every $k\in K$, we have by \Cref{complexnonvanishing}, \[Ad_k^*(\iota_u\p^*\omega_{G/H})(X_1\wedge\overline{X}_1\ldots\wedge X_d\wedge \overline{X}_d)\geq 0,\]
and this inequality is strict in an open neighborhood of the base point of $K/L$. Hence, by integrating over $k$ and using \Cref{def: Push-forward}, we get $\pi_*p^*\omega_{G/H}\neq 0$ at $x$.

Conversely, assume that $X$ is not generically transverse at $x$. Then for every $g\in G$, we have $T_xX+T_x(gH/L)\subsetneq T_xG/K$. Hence for every $d$-uple of $\C$-linearly independent vectors $X_1,\ldots,X_{d}$ in $T_xX$, the intersection of the subspaces $\mathrm{span}_\R(X_1,\overline{X}_1,\ldots,X_{d},\overline{X}_{d})$ and $T_x(gH/L)_\C$ is non-empty. Hence the form $\iota_{u}\omega_{G/H}$ vanishes on the multi-vector $X_1\wedge\overline{X}_1\ldots X_{d}\wedge\overline{X}_d$. The same is true for $Ad(k)^*(\iota_u\p^*\omega_{G/H})$ for all $k\in K$. By integrating, we get that $\pi_*p^*\omega_{G/H}$
vanishes at $x$. Hence the result.
\end{proof}

\begin{proof}[Proof of Theorem \ref{t:densitycriterion}]

The implication $(ii)\Rightarrow (i)$ and $(i)\Rightarrow (iii)$ are obvious, and the equivalence $(iii)\Rightarrow (iv)$ readily follows from Proposition \ref{generictransversality}. We only have to prove $(iii)\Rightarrow (i)$.

By \Cref{r:density}, the set of points where $X$ is transverse to $H$-orbits is open and dense. Let $x$ be such a point and $g\in G_\R$ such that $\tilde X$ and $gH_\R/L_\R$ intersect transversally at $x$. By Weak Approximation, $G(\Q)$ is analytically dense in $G(\R)$. Thus there exists a sequence $g_n \in G_\Q$ converging to $g$. For $n$ large enough, by stability of transversality, $g_n H_\R/L_\H$ intersects $\tilde X$ transversally at a point $x_n$ such that $x_n \underset{n\to +\infty}{\longrightarrow} x$. Since $g_n \in G_\Q$, $g_n Hg_n^{-1}$ is a $\Q$-subgroup of $G$, hence $x_n$ belongs to the transverse Hodge locus of type $H$.
\end{proof}

Unfortunately, in many situations, variations of Hodge structure are never generically transverse to $H/L$. Indeed, Griffiths' transversality constrains their tangent space to be contained in the Griffith distribution, so that it cannot supplement $T_x(H/L)$ in other directions.

To be more precise, let $H/L \subset G/K \subset \D$ be Mumford--Tate domains. Let $\varphi:\mathbb{U}^1 \to H$ be the restriction of the representation at a point $x\in H/L$ of the Deligne torus to the unit circle. Then both $\fg$ and $\fh$ are invariant under the adjoint action of $\varphi$. We thus have decompositions 
\[\fg_\C = \bigoplus_{p=-k}^k \fg^{p,-p}~,\quad \fh_\C = \bigoplus_{p=-k}^k \fh^{p,-p}\]
with $\fh^{p,-p}\subseteq \fg^{p,-p}$. Note that $\fg^{0,0} = \fk$ and $\fh^{0,0} = \fl$.

\begin{proposition}\label{p:level2}
The following are equivalent:
\begin{itemize}
    \item There exists a smooth variation of Hodge structure in $G/K$ which is generically transverse to $H/L$.
    \item For all $|p|\geq 2$, $\fh^{p,-p} = \fg^{p,-p}$ and there exists an abelian subalgebra $\mathfrak{a} \subset \fg^{-1,1}$ such that 
    \[\mathfrak a + \fh^{-1,1} = \fg^{-1,1}~.\]
\end{itemize}
\end{proposition}

\begin{proof}
Assume that there exists a smooth variation of Hodge structure $X\subset G/K$ which is generically transverse to $H/L$. Up to left multiplication by some $g\in G$, we can assume that $x\in X$ and that 
\begin{equation} \label{eq:Transversality}
T_x^{1,0} X + T_x^{1,0} H/L = T_x^{1,0}G/K\end{equation}

Now, since $X$ is a variation of Hodge structure, $T_x^{1,0} X$ is an abelian subalgebra of $\fg^{-1,1}$, while $T_w^{1,0}H/L = \bigoplus_{p<0} \fh^{p,-p}$. The identity \eqref{eq:Transversality} thus implies that
\[T^{1,0}_x X + \fh^{-1,1} = \fg^{-1,1}\]
and
\[\fh^{-p,p} = \fg^{-p,p}\]
for all $p\geq 2$, and by Hodge symmetry, also for $p\leq -2$.\\

Conversely, assume $\fh^{-p,p} = \fg^{-p,p}$ for all $p\geq 2$ and $\mathfrak a^{-1,1} + \fh^{-1,1} = \fg^{-1,1}$ for some abelian Lie subalgebra $\mathfrak a$. Let $A$ denote the complex abelian subgroup of $G_\C$ spanned by $\mathfrak a$. Recall that $G_\C$ acts on the compactified period domain $\bar{ \mathcal D }$. If $U$ is a sufficiently small neighbourhood of the identity in $A$, then
\[X= \{ a\cdot x, a\in U\}\]
is a smooth holomorphic submanifold contained in $\mathcal D$. Since $\mathfrak a \subset \fg^{-1,1}$ and $A$ is abelian, $X$ is tangent to the Griffith distribution at every point. Hence $X$ is a smooth variation of Hodge structure transverse to $H/L$ at $x$.
\end{proof}

\begin{remarque}
Baldi--Klingler--Ullmo \cite[Prop. 6.5]{baldiklinglerullmo} proved that, when $\fg$ is simple and has components $\fg^{k,-k}$ for $k\geq 3$, then a Lie subalgebra $\fh$ can never contain $\fg^{k,-k}$, $k\geq 3$. It follows that the transverse Hodge locus is always empty in that case.
\end{remarque}

\subsection{Chern classes of the Hodge bundles}

Let $\{\V_\Z,\cF^\bullet \cV,B\}$ be a variation of Hodge structure of weight $k$ over a complex analytic variety $X$. Let $\sigma$ be the antilinear automorphism of $\mathcal V$ given by 
\[\sigma_{\vert \mathcal V^{p,q}}: v\mapsto i^{p-q} \bar v~.\]
and let $h$ be the Hermitian form 
\[h(v,v) = B(v,\sigma v)~.\]
By definition of Hodge structures, $h$ is positive definite and the decomposition 
\[\cV = \bigoplus_{p+q=k} \cV^{p,q}\]
is orthogonal for $h$.

We have $\sigma^2= (-1)^k \mathrm{Id}_{\cV}$. Define now a new linear connection $\nabla_h$ on $\mathcal V$ by
\[\nabla_h = \nabla + \frac{(-1)^k}{2}\sigma(\nabla \sigma)~.\]
(The connection $\nabla_h$ is the part of the connection $\nabla$ that preserves $\sigma$.)
Then $\nabla_h$ preserves the metric $h$ and the orthogonal decomposition 
\[\cV = \bigoplus_{p+q=k} \cV^{p,q}~.\]
Let ${\nabla_h}^{p,q}$ denote the induced Hermitian connection on $\cV^{p,q}$. and $\Theta_h^{p,q}$ denotes its curvature. One can show that $\Theta_h^{p,q}$ is of type $(1,1)$.

\begin{definition}
The Chern forms of $\cV^{p,q}$ are the $(\ell,\ell)$ forms $c_\ell(\mathcal V^{p,q})$, $1\leq l \leq h^{p,q}$ defined by
\[\det\left(\mathrm{I}_{\mathcal V^{p,q}}+\frac{i}{2\pi} \Theta_h^{p,q}\right)=1+ \sum_{\ell=1}^{h^{p,q}} c_\ell(\mathcal{V}^{p,q}).\]
\end{definition}

It is well-known that the form $c_\ell(\mathcal V^{p,q})$ represents the $\ell$\textsuperscript{th} Chern class of $\mathcal V^{p,q}$ in de Rham cohomology of $X$.

These Chern forms turn out to be pull-backs of invariant forms under the period map. Indeed, $\sigma$, $h$, and $\nabla_h$ can be defined on the universal Hodge decomposition
\[\mathcal U = \bigoplus_{p=0}^k \mathcal U^{p,q}\]
over $\D$. There, these objects are $G$-equivariant and induce $G$-invariant Chern forms 
$c_\ell(\mathcal U^{p,q}).$
These factor to the quotient $G_\Z\backslash \mathcal{D}$ and, if $f: X \to G_\Z\backslash \mathcal{D}$ denotes the period map, we have
\[c_\ell(\mathcal V^{p,q}) = f^* c_\ell (\mathcal U^{p,q})~.\]

\subsubsection{Expression at the Lie algebra level}\label{expressionliealgebra}

Let us now express the Chern forms $c_\ell(\mathcal U^{p,q})$ at the Lie algebra level. 

Let us fix a base-point $o$ in $\mathcal{D}$ with stabilizer $K$. The group $K$ decomposes as
\[K = \prod_{p=1}^{\lfloor\frac{k}{2}\rfloor}K^{p,q~,}\]
where $K^{p,q} \simeq \mathrm U(h^{p,q})$ for $p>q$ and $K^{k',k'} \simeq \mathrm{O}(h^{k',k'})$ when $k= 2k'$ is even.

When $p\geq q$ (resp. $p\leq q$), the bundle $\mathcal U^{p,q}$ is the bundle associated to the linear representation of $K$ that factors through the standard representation of $K^{p,q}$ (resp. the dual representation). Let $\fk^{p,q}$ denote the Lie algebra of $K^{p,q}$. Then, for $p\geq q$, the curvature of $\mathcal U^{p,q}$ at $x$ is the $2$-form on $\fg/\fk$ with values in $\fk^{p,q}\subset \End(\mathcal U^{p,q})$ given by
\[F_h^{p,q}(u,v) = \pi_{p,q}([u,v])- [\pi_{p,q}(u), \pi_{p,q}(v)]~,\]
where $\pi_{p,q}: \fg \to \fk^{p,q}$ denotes the orthogonal projection for the Killing metric.

\subsubsection{Chern forms on the compact dual}

Recall from \Cref{compactdualitysection} that we have an isomorphism of differential algebras
\[\psi: \Omega^\bullet(G/K,\C)^G \overset{\sim}{\to} 
\Omega^\bullet(G_U/K,\C)^{G_U}\]
which consists in identifying both spaces with $\Lambda_\C^\bullet(\fg_\C/\fk_\C)^{K_\C}$.

We now wish to identify the invariant forms on $G_U/K$ corresponding to the Chern forms on $G/K$.

Recall that $\hat{\mathcal{D}}= G_\C/P$ is the space of complex flags
\[0\subseteq F^k \subseteq \ldots \subseteq F^0 = V_\C\]
such that 
\[F^{k-p} = {F^p}^\bot\]
(where the orthogonal is intended with respect to the bilinear form $B$) and
\[\dim(F^p/F^{p+1}) = h^{p,q}~.\]
Let $\hat{\mathcal U}$ denote the trivial bundle over $\hat{\D}$ equipped with the action of $G_\C$ given by the standard linear action in the fibers. The bundle $\hat{\mathcal{U}}$ admits a tautological filtration $\cF^\bullet \hat{\mathcal U}$ by $G_\C$-equivariant vector bundles which is given at a point $x$ by the flag defining $x$.

By construction, the restriction of $\cF^\bullet \hat \cU$ to the open domain $\D$ is the filtration $\cF^\bullet \cU$ of $\cU$ given (in $C^\infty$) by
\[\cF^p \cU = \bigoplus_{p'\geq p} \cU^{p', k-p'}~.\]
\[\cF^p \cU = \bigoplus_{p'\geq p} \cU^{p', k-p'}~.\]

Let us now prove that the dual space $G_U/K$ identifies with $\hat {\D}$.

\begin{proposition}
The group $G_U$ is (conjugated to) the subgroup of $G_\C = \Aut_\C(B)$ commuting with the antilinear automorphism $\sigma$.
\end{proposition}

\begin{proof}
Let $\tau: \mathcal V_o \to \mathcal V_o$ be the complex conjugation and $\theta= \sigma \tau$. Conjugation by $\tau$ is the anti-holomorphic involution of $G_\C$ fixing $G$ and one verifies that the conjugation with $\sigma \tau$ is a Cartan involution of $G$ fixing $K$. With respect to this choice of Cartan involution, the group $G_U$ is then the fixed point set of conjugation by $\sigma$.
\end{proof}

\begin{corollaire}
The group $G_U$ acts transitively on $\hat{\D}$, and the stabilizer of $o$ in $G_U$ is $K$.
\end{corollaire}

\begin{proof}
As a maximal compact subgroup of $G_\C$, the group $G_U$ acts transitively on the flag variety $\hat{\D}$, and the stabilizer $K'$ of $o$ preserves the flag $\cF^\bullet \cU_o$. Now, since $G_U$ commutes with $\sigma$, it preserves the Hermitian form $B(\cdot, \sigma\cdot)$. Therefore, $K'$ preserves the orthogonal of $\cF^{p+1} \cU_o$ in $\cF^p \cU_o$ for $B(\cdot, \sigma\cdot)$, which is precisely $\mathcal U^{p,q}$. We conclude that $K= K'$.
\end{proof}

The $G_U$-invariant form $B(\cdot, \sigma \cdot)$ induces a flat $G_U$-invariant Hermitian metric $\hat h$ on $\hat \cU$. Let $\hat \cU^{p,q}$ denote the $\hat h$-orthogonal of $\cF^{p+1}\hat \cU$ in $\cF$. Then the the bundle $\hat \cU^{p,q}$ is $G_U$-invariant and carries a $G_U$-invariant Hermitian connection $\nabla_{\hat h}^{p,q}$ with curvature form $\Theta_{\hat h}^{p,q}$. The Chern forms of this connection define $G_U$-invariant forms 
\[c_\ell(\hat \cU^{p,q})~\]
which represent the Chern classes of $\hat \cU^{p,q} \simeq \cF^p\hat \cU/\cF^{p+1} \hat U$ on $\widehat{\D}$.

\begin{proposition}\label{samechernclass}
The isomorphism $\phi:\Omega^{2l}(G/K,\C)^G\to \Omega^{2l}(G_U/K,\C)^G$ maps $c_\ell(\cU^{p,q})$ to $c_\ell(\hat \cU^{p,q}) \in \Omega^{2l}(G_U/K,\C)^G$.
\end{proposition}

\begin{remarque}
The isomorphism $\phi$ is \emph{not} induced by the identification
\[T_o \mathcal D = T_o \widehat{\mathcal{D}}\]
coming from the inclusion $\mathcal D\subset \widehat{ \mathcal D}$ but rather from the diagram in \Cref{p:CommutativityPushPullCompactDual}.
\end{remarque}

\begin{proof}
For $p\geq q$, the bundle $\hat \cU^{p,q}$ is the vector bundle on $G_U/K$ associated to the linear representation of $K$ factoring through the standard representation of $K^{p,q}$. Hence its curvature form at $o$ is given by a formula similar to \Cref{expressionliealgebra}. 

Let now $\pi^{p,q}_\C$ denote the orthogonal projection of $\fg_\C$ to $\fk_\C^{p,q}$ for the complex Killing form. Then $\pi^{p,q}_\C$ restricts to the orthogonal projection to $\fk$ on both $\fg$ and $\fg_U$.

Therefore, both the curvature forms of $\cU^{p,q}$ and $\hat \cU^{p,q}$ at $o$ are given by
\[(u,v) \mapsto \pi_\C^{p,q}([u,v])- [\pi_\C^{p,q}(u), \pi_\C^{p,q}(v)]~,\]
hence all the symmetric polynomials in those curvature forms are identified by $\phi$.
\end{proof}

\subsubsection{Characteristic cohomology}

As mentioned in \ref{thesymmetriccase}, there might be $G$-invariant forms on $\D$ which are not closed, in which case $G$-invariant closed forms are not characterized by the corresponding cohomology class in $\H^\bullet(\hat{\D})$.

In the context of variations of Hodge structure, however, we are ultimately interested in the restriction of $G$-invariant forms to submanifolds that are tangent to the Griffiths' distribution. This motivates the introduction of the \emph{characteristic cohomology} of a period domain, which, roughly speaking, restricts the differential algebra of invariant forms to the Griffiths' distribution (see \cite[III.A]{griffithsgreenkerrarticle}).

We do not define this notion here and only mention the analogous of É. Cartan's theorem, which comforts the idea that the geometry of period domains is similar to that of symmetric spaces ``in restriction to the Griffiths' distribution''.
\begin{proposition}\label{vanishingonvhs}
Let $X$ be a complex manifold and $f:\tilde X \to G/K$ the period map of a variation of Hodge structure. Then, for every $\alpha \in \Omega^\bullet(G/K,\C)^G$, the pull-back form $f^*\alpha$ is closed of bidegree $(p,p)$ for some $p$.
\end{proposition}

\begin{proof}
Let $\alpha$ be a $G$-invariant form on $G/K$. Let $x$ be a point in $X$ and $\varphi: \C^\times \to G_\R$ the representation of the Deligne torus defining the Hodge structure $f(x)$. Then $\varphi(\mathbb{U}^1)$ is a subgroup of $\mathrm{Stab}(x)\subset G$ and acts on $W^{-1,0}$ by complex multiplication, where $W$ is the tangent space at $f(x)$. Since ${\alpha_{f(x)}}$ is $\varphi(\mathbb{U}^1)$-invariant, it must belong to $\Lambda^{p,p}(W_x^*)$ for some $p$, and we conclude that $f^*\alpha$ has bidegree $(p,p)$ since $f$ is holomorphic. In particular $f^*\alpha$ has even degree. 

Now, $d\alpha$ is also a $G$-invariant form and $f^*(d \alpha) = d(f^*\alpha)$ has odd degree. By the previous argument, it must vanish.
\end{proof}

\begin{corollaire}
Let $\alpha$ be a closed invariant form on $\mathcal D$. Then the pullback of $\alpha$ by any variation of Hodge structure is completely determined by the cohomology class $[\phi(\alpha)]\in \H^\bullet(\widehat{\D})$.
\end{corollaire}

\begin{proof}
Let $\alpha'$ be another closed $G$-invariant form on $G/K$ such that $\phi(\alpha- \alpha')$ is exact on $G_U/K$. We can write $\alpha-\alpha'= d \beta$, where $\beta$ is $G$-invariant. Let now $f:X\to G/K$ be a variation of Hodge structure. Then $f^*\beta$ is closed by the previous proposition, hence 
\[f^*\alpha - f^*\alpha' = d(f^*\beta) = 0~.\]
\end{proof}

\begin{remarque}
We only mentioned these results for period domains, but one can prove that they remain true on every Mumford--Tate domain.
\end{remarque}


\subsection{Examples}

We now apply the previous considerations to compute the pull-push form in various examples.

\subsubsection{Noether--Lefschetz loci in weight $2$}

Assume in this section that $\mathcal{D}$ is the period domain for a polarized variation of Hodge structure of weight $2$ on a quadratic lattice $(V,B)$ which is assumed to be of signature $(p,2q)$. Let $R$ be a rational subspace of $V \otimes_\Z \Q$ of rank $r\leq h^{1,1}$ such that $B$ is positive definite in restriction to $R\otimes_\Q \R$, and let $\mathcal D_R \subset \mathcal{D}$ be the set of Hodge structures $x\in\mathcal{D}$ such that $R\subset \mathcal V^{1,1}_x$.

Choose a basepoint $o$ in $\mathcal D_R$. Let $K$ be the stabilizer of $o$ in $G$, $H$ be the subgroup of $G$ fixing $R$ and $L=K\cap H$. Then $H$ is a Mumford--Tate group and
$\mathcal D_R \subset \mathcal{D}$ is the Mumford--Tate domain $H/L\subset G/K$.

Denoting as before by $p$ and $\pi$ the respective projections from $G/L$ to $G/H$ and $G/K$, we can now prove the following:

\begin{theorem}\label{pullpushweight2}
Let $X$ be an smooth complex analytic manifold, let $\cV = \cV^{2,0} \oplus \cV^{1,1} \oplus \cV^{0,2}$ be the $C^\infty$ Hodge decomposition of a variation of Hodge structure of weight $2$ and Hodge numbers $(q,p,q)$ on $X$ and let $\tilde f:\tilde X \to \mathcal{D}$ be the corresponding period map. Then
\[ f^*\left(\pi_*p^* \omega_{G/H}\right) = \Vol(G_U/H_U)\cdot c_q(\cV^{2,0})^r~.\]
\end{theorem}

\begin{proof}

Let $\sigma$ be the antilinear automorphism defined in the previous section. Since we are in even weight, $\sigma$ is an involution which fixes a real form $\cU^\sigma$ of $\cU_o$ on which the symmetric form $B$ is real and positive definite. Since $\sigma$ coincides with the standard complex conjugation on $\mathcal U^{1,1}_o$, the subspace $R$ is contained in $\cU^\sigma$. 

Now, $G_U$ is the subgroup of $G_\C = \Aut_\C(B)$ preserving $\cU^\sigma$ and $H_U= G_U\cap H_\C$ is the subgroup of $G_U$ fixing $R$. Therefore, $H_U/L$ is the domain $\hat {\mathcal D}_R \subset \hat {\mathcal{D}}$ where $\hat \cU^{1,1}$ contains $R$. Since $R$ is $\sigma$-invariant and 
\[\cF^1\hat \cU \cap \sigma(\cF^1 \hat \cU) = \hat \cU^{1,1}~,\]
we also have that
\[\hat{\mathcal {D}}_R = \{x\in \hat{\mathcal{D}} \mid \cF^1 \hat \cU_x \supseteq R\}~.\]

Let $(u_1,\ldots, u_r)$ be a basis of $R$. The projection of $u_\ell$ into $\cF^0\hat \cU/\cF^1 \hat \cU$ defines a holomorphic section $s_\ell$ of $\cF^0\hat \cU/\cF^1 \hat \cU$, and $\widehat{\mathcal {D}}_R$ is the transverse intersection of the vanishing loci of all the $s_\ell$. We conclude that $\hat{\mathcal{D}}_R$ is Poincaré dual to the $r$\textsuperscript{th} power of the Euler class of $\cF^0\hat \cU/\cF^1\hat \cU$ i.e.,
\[c_q(\hat \cU^{0,2})^r~.\]

By \Cref{dualofahomologyclass}, we have
\[\pi_*p^*\omega_{G_U/H_U} = \Vol(G_U/H_U)\cdot c_q(\hat \cU^{0,2})^r + \d \alpha\]
for some invariant form $\alpha$. 

By \Cref{equalityofforms} and \Cref{samechernclass} we have
\begin{eqnarray*}
\pi_*p^*\omega_{G/H} &=& i^{2qr} \pi_*p^*\omega_{G_U/H_U}\\
&=& (-1)^{qr}c_q(\hat \cU^{0,2})^r + (-1)^{qr}\d \alpha\\
&=&c_q(\hat \cU^{2,0})^r + (-1)^{qr}\d \alpha\\
&=& c_q(\cU^{2,0})^r + (-1)^{qr} \d \alpha~.
\end{eqnarray*}

Finally, by \Cref{vanishingonvhs}, the pull-back of $\d \alpha$ by the period map of a variation of Hodge structure vanishes, and the conclusion follows.
\end{proof}

\subsubsection{Diagonal embedding of Shimura varieties}

Let $G_1$ be a semi-simple Lie group of Hermitian type and $K_1$ a maximal compact subgroup, so that $\mathcal D \equaldef G_1/K_1$ is a Hermitian symmetric space of non-compact type. We apply the results of previous sections to $G=G_1\times G_1$ and $H=\Delta(G_1)$, the diagonal embedding of $G_1$. Let $\Delta: \mathcal D\hookrightarrow \mathcal D\times \mathcal D$ be the corresponding diagonal embedding of symmetric spaces.

First, remark that since $\D$ and $\D\times \D$ are Hermitian symmetric, their tangent space is equal to the Griffiths' distribution. Hence by \Cref{vanishingonvhs}, the complex $\Omega^\bullet(\D\times \D,\C)^G$ is supported in even degrees and is isomorphic to the complex
$H^{\bullet}_{dR}(\hat \D\times \hat \D,\C)$.

Recall the following classical result. Let $X$ be a closed orientable smooth manifold of dimension $n$.
For all $0\leq k \leq n$, let us fix a basis $\left([\alpha^{k,i}]\right)_{i\in J_k}$ of $H^k(X,\C)$ and denote by $\left(\alpha_{k,i}^\vee\right)_{i\in J_k}$ the dual basis of $H^{n-k}(X,\C)$ with respect to Poincar\'e pairing. Let $\pi_1,\pi_2:X\times X\rightarrow X$ denote respectively the projections onto the first and the second factor.
Then by \cite[Lemma 1.22]{bottTu}, the cycle class of the diagonal $\Delta(X)\hookrightarrow X\times X$ is Poincaré-dual to the de Rham cohomology class  \begin{align}\label{dualofdiagonal}
\gamma_X=\sum_{k=0}^{n} (-1)^{n(n-k)} \sum_{i\in J_k}\pi_1^*[\alpha_{k,i}^\vee]\wedge\pi_2^*[\alpha_{k,i}] \in H_{dR}^n(X\times X,\C)~.\end{align}

We can now state the main theorem of this section. Let $(\gamma_{k,i})_{i\in J_k}$ be a basis of $\Omega^{2k}(\mathcal D, \C)^{G_1}$ and let $(\gamma_{k,i}^\vee)_{i\in J_k}$ denote the dual basis of $\Omega^{2d-2k}(\mathcal D,\C)^{G_1}$, i.e. such that
\[\gamma_{k,i}\wedge \gamma_{k,j}^\vee =  \frac{\delta_{i,j}}{\Vol(\hat{\mathcal D})} \omega_{\mathcal D}~,\]
where $\omega_{\mathcal D}$ denotes the invariant volume form of $\mathcal D\simeq G/K$ as in \Cref{notations}.

\begin{theorem}\label{diagonalrestriction}
Set $G= G_1\times G_1$, $K = K_1\times K_1$ where $K_1$ is a maximal compact subgroup of $G$ and take $H$ the diagonal embedding of $G_1$ into $G$. Let $\pi_1$ and $\pi_2$ denote the projections of $G/K= \mathcal D \times \mathcal D$ on the first and second factor.
Then  \[\frac{1}{\Vol(G_U/H_U)}\pi_*p^*\omega_{G/H}=\sum_{\underset{i\in J_k}{0\leq k\leq d}}\pi_1^*(\gamma_{k,i}^\vee)\wedge \pi_2^*(\gamma_{k,i})~. \]
In particular, its pull-back by the diagonal embedding $\Delta$ of $\mathcal D$ is given by
\[\frac{1}{\Vol(G_U/H_U)}\Delta^{*}\pi_*p^*\omega_{G/H}=\frac{\chi(\hat{ \mathcal D})}{\Vol(\hat{ \mathcal{D}})}\,\omega_{\mathcal D}~.\]
where $\chi(\hat{\mathcal D})>0$ is the Euler characteristic of $\widehat{\mathcal D}$.
\end{theorem}

\begin{proof}
 By \Cref{equalityofforms} and \Cref{dualofahomologyclass}, it is enough to determine the cohomology class of the corresponding pull-push form on the compact dual $\hat{ \mathcal D} \times \hat {\mathcal D}$. By \Cref{dualofahomologyclass}, this cohomology class is Poincaré dual to the diagonal embedding of $\hat{\mathcal D}$. The conclusion now follows from \Cref{dualofdiagonal}.
\end{proof}

\subsubsection{Hodge locus in Shimura varieties}
In this section, we prove  \Cref{t:shimura_simple} and \Cref{C:unitary-shimura}.

Let $G$ be a semi-simple Lie group of Hermitian type, $\D$ the associated Hermitian symmetric space, $\Gamma\subset G$ an arithmetic subgroup and $S=\Gamma\backslash D$. Let $(H,\D_H)$ be a Shimura subdatum such that $\pi_*p^*\omega_{G/H}$ is a positive form of type $(k,k)$. In particular, its restriction to any subvariety of $S$  of dimension at least $k$ is non-zero. Hence by the equivalences from \Cref{generictransversality}, the Hodge locus in $X$ is analytically dense and equidistributed with respect to $\pi_*p^*\omega_{G/H}$.

For the second part of the theorem, the pull-push form associated to $G/H$ is a $(1,1)$-form and since $G$ is absolutely irreducible, there is, up to a scalar, a unique $(1,1)$-from on $\D$ which is given as the Chern form of the canonical bundle on $\D$. The latter is known to be K\"ahler. Hence $\pi_*p^*\omega_{G/H}$ is K\"ahler and we conclude as before. 

We now  prove \Cref{C:unitary-shimura}. Let $n\geq 1$ and let $G=\mathrm{SU}(n,1)$. For $1\leq r\leq n$, let $H=\mathrm{SU}(n-r,1)$. Let $K=S(U(1)\times U(n))$ be the maximal compact subgroup of $G$ and let $\D=\mathbb{B}^{n}\simeq G/K$  be the unit ball which is isomorphic to the symmetric space of $G$. The natural $1$-dimensional representation of $U(1)$ on the determinant of the cotangent bundle of $\D$ determines a Hermitian line bundle on $\mathbb{B}^{n}$ with first Chern form $\omega$. Using a similar method as in \Cref{pullpushweight2}, one can prove easily that: 
\begin{proposition}
We have $\pi_*p^*\omega_{G/H}=\Vol(G_U/H_U)\omega^r$. 
\end{proposition}
Let $\Gamma$ be a neat arithmetic subgroup of $G$.  The quotient $S=\Gamma\backslash \mathbb{B}^n$ is a unitary Shimura variety and the form $\omega$ is K\"ahler on $S$. If $r=1$, then we are in the situation of \Cref{t:shimura_simple}. Hence  the Hodge locus is dense and equidistributed in any subvariety of $S$ of positive dimension.

\subsubsection{Hodge locus in $\mathcal{A}_g$}\label{s:hodge-locus-ag}
In this section, we prove \Cref{T:densityabelianvarieties}.

Let $(V,\Psi)$ be a rational vector space of dimension $2g$ endowed with a non-degenerate symplectic form $\Psi$. Let $G=\mathrm{Sp}_{2g}(\R)$ and let $\mathbb H_g$ be the Siegel upper-half space which is the Hermitian space associated to $G$. 

For $1\leq k\leq \frac{g}{2}$, let $V_k\subseteq V$ be a non-degenerate rational subspace of rank $2k$ and let $H\subseteq G$  be the stabilizer of $V_k$ in $G$. Then $H\simeq \mathrm{Sp}_{2k}(\R)\times \mathrm{Sp}_{2g-2k}(\R)$ and its symmetric space is equal to  $\mathbb{H}_k\times\mathbb{H}_{g-k}$. 
\medskip

To compute the pull-push form $\pi_*p^*\omega_{G/H}$, we follow the general method explained in \Cref{s:pullpush} and \Cref{definitionvhs}. The compact dual $Y_g$ of $\mathbb{H}_g$ is equal to the space of Lagrangian subspaces of $V_\C$ and the compact dual of $\mathbb{H}_k\times\mathbb{H}_{g-k}$ is the subspace $Y_k\times Y_{g-k}$ where $Y_k$ and $Y_{g-k}$ are the space of Lagrangian subspaces of $V_{k,\C}$ and $V_{k,\C}^{\perp}$ respectively. The inclusion $Y_k\times Y_{g-k}\hookrightarrow Y_g$ is then given by taking direct sums of Lagrangians in a compatible way with the decomposition $V_\C=V_{k,\C}\oplus V_{k,\C}^\bot$. 

Let $\mathcal{V}\rightarrow Y_g$ be the trivial vector bundle of rank $2g$ determined by $V$ and let $\widehat{\cF}^{1}\rightarrow Y_g$ be the Hodge vector bundle whose restriction to $\mathbb{H}_g$ will be denoted simply $\cF^{1}$. 
Let $\mathcal{V}_k$ be the trivial vector bundle determined by $V_k$. Then we have a natural map of vector bundles: 
\[f:\mathcal{V}_k\rightarrow \mathcal{V}/\widehat{\cF}^{1}.\]
The locus where this map has rank at most $k$ corresponds to the locus where the rank of the kernel is at least $k$. Since the kernel is Lagrangian in $V_k$, it is also the locus where the rank is exactly $k$ and hence it is equal to $Y_k\times Y_{g-k}$. 
By the Giambelli--Porteous--Thom formula \cite{kempf-laksov}, the locus $Y_k\times Y_{g-k}$ is Poincaré  dual to the class \[\det\left((c_{g-k+i-j}(\mathcal{V}/\widehat{\cF}^{1}))_{1\leq i,j\leq k}\right)\]

By \Cref{equalityofforms} and \Cref{samechernclass}, we have: 
\begin{align*}
    \pi_*p^*\omega_{G/H}&=i^{2k(g-k)}\pi_*p^*\omega_{G_U/H_U}\\
    &=(-1)^{k(g-k)}\Vol(G_U/H_U)\det\left(((-1)^{g-k+i-j}c_{g-k+i-j}({\cF}^{1}))_{1\leq i,j\leq k}\right) \\
    &=\Vol(G_U/H_U)\det\left(((-1)^{i-j}c_{g-k+i-j}(\cF^{1}))_{1\leq i,j\leq k}\right)\\
    &=\Vol(G_U/H_U)s_k.
\end{align*}

Now combining \Cref{generictransversality}, and \Cref{main}, the first part of \Cref{T:densityabelianvarieties} easily follows. For the second part, we use the Main Theorem of \cite{keel-sadun} which stipulates that $s_1=c_{g-1}$ is non-zero restricted to any compact subvariety of dimension $>\frac{(g-1)(g-2)}{2}$. 
\medskip

\section{Applications}\label{vhs}

We discuss in this section various applications of \Cref{main}. They concern mainly equidistribution of Hodge loci in variations of Hodge structures, in particular in the context of weight $2$ Hodge structures and Hodge structures parametrized by Shimura varieties. For an introduction to these topics, we refer to \cite[III,VI]{voisin} and \cite{griffiths}.

\subsection{Refined Noether--Lefschetz loci in $\Z$-PVHS of weight $2$.}\label{refinedloci}
Let $\{\mathbb{V}_\Z,\mathcal{F}^\bullet\mathcal{V},B\}$ be a polarized variation of Hodge structure of weight $2$ over a smooth complex quasi-projective variety $X$. 
Assume as before that the local system $(\V_\Z,B)$ has fibers  isomorphic to a quadratic lattice $(V_\Z,B)$ equipped with a bilinear form \[B:V_\Z\times V_\Z\rightarrow \Z\]
with associated quadratic form $\frac{B(y,y)}{2}\in\Z$ for $y\in V_\Z$. Let $(q,p,q)$ be the Hodge numbers and $d=p+2q$ the rank of $\V_\Z$.

For $x\in X$, let $\rho(x)$ be the rank of the \emph{Picard group} $\mathrm{Pic}(x)=\cF^1\mathcal{V}_x\cap \V_{\Z,x}$. Assume that the variation is simple, i.e., $\rho(x)=0$ at a very general point. For $r\geq 1$, we introduce the \emph{refined Noether--Lefschetz locus}\footnote{Historically, Hodge loci are referred to as Noether--Lefschetz loci in weight $2$.}
\[NL^{\geq r}=\{x\in X, \rho(x)\geq r\}.\]
It corresponds to a sub-Hodge locus for $\mathbb{V}_\Z^{r}$. It can be written as the union over algebraic subvarieties in the two following ways: 
\begin{enumerate}
    \item It is the union, over all integers $N\geq 1$, of the sets
    \[\{x\in X\mid \exists P\subseteq \mathrm{Pic}(x) \textrm{ of rank $r$},  \disc(P)\leq N \}.\]
    \item It is the union, over all positive definite symmetric matrices $M\in M_r(\Z)$, of the sets \[\{x\in X\mid  \exists (\lambda_1,\cdots,\lambda_r)\in \mathrm{Pic}(x), \left(B(\lambda_i.\lambda_j)\right)_{1\leq i,j\leq r}=M \}.\]
    
\end{enumerate}
We will prove in the next two subsections that both formulations give equidistribution results and hence we prove \Cref{maink3-alt} and  \Cref{maink3} by using different techniques in each case.

\subsubsection{Equidistribution on average}\label{S:averageequidistribution}

Let $x$ be a point in $X$ and denote by $(V_\Z,B,F^\bullet)$ the fiber at $x$ of $\V$. As in the previous section, we set $G= \mathrm O(V_\R, B)$, $\Gamma = \mathrm O(V_\Z,B)$ and let $K\subset G$ be the stabilizer of $F^\bullet$. Finally, let $H$ denote the stabilizer of a positive definite subspace of $V_\R$ of dimension $r$, so that $G/H$ identifies with the space $\mathcal P$ of positive definite real subspaces of dimension $r$.

Define the \emph{discriminant} of a rational subspace $W\in \mathcal P$ as the determinant of the intersection matrix
\[I(W) = \left(B(v_i,v_j)\right)_{1\leq i,j\leq r}\]
where $(v_i)$ is a basis of $W\cap V_{\Z}$. We denote by $\mathcal P_n$ the discrete subset of $\mathcal P$ consisting of rational subspaces of discriminant $\leq n$. The set $\mathcal P_n$ is a finite union of $\Gamma$-orbits of $G/H$, corresponding to a finite union $\mathcal O_n$ of closed $H$-orbits in $\Gamma \backslash G$. We prove here that $\mathcal P_n$ is equidistributed in $G/H$.

\begin{theorem}\label{T:averageequi}
The sequence $(\mathcal P_n)_{n\in \N}$ is equidistributed in $G/H$.
\end{theorem}

The proof we give here is a refinement of the fact that integral vectors of length at most $n$ equidistribute. Some trick is needed in order to get rid of multiplicities, but the proof is ``elementary'' in the sense that it does not rely on any involved argument such as the circle method, automorphic functions or Ratner theory. Of course, by counting all rational subspaces of length \emph{less than $n$}, we avoid all the difficult arithmetic questions that arise when looking at rational subspaces of a fixed discriminant.

Let $\Omega$ denote the open cone of $V_\R^r$ consisting of tuples of vectors spanning a positive definite subspace of dimension $r$. The function
\begin{align*}
    h:\Omega &\to \R_{>0} \\
     \underline v=(v_1,\ldots, v_r)&\mapsto  \det\left(B(v_i,v_j)\right)_{1\leq i,j\leq r}
\end{align*}
is homogeneous of degree $2r$. We denote by $\Omega_1$ the hypersurface $\{\underline v\in \Omega \mid h(\underline v) = 1\}$ and by $\Omega_{\leq 1}$ the subset
$\{\underline v \mid h(\underline v) \leq 1\}$. We denote by 
\[pr:\Omega \to \Omega_1\]
the projection map
\[\underline v \mapsto h(\underline v)^{-\frac{1}{2r}} \underline v~.\]

We endow $V_\R^r$ with the Lebesgue measure for which $V_\Z^r$ has covolume $1$ and denote by $\omega$ the push-forward by $pr$ of the Lebesgue measure restricted to $\Omega_{\leq 1}$ (i.e. the volume form such that
\[\int_U \omega = \mathrm{Leb}\left(\bigcup_{0<t \leq 1} tU\right)~.\]

Define
\[\mathcal Q_n =\{\underline v\in \Omega\cap V_\Z^r\mid h(\underline v) \leq n\}\]
and let $\mu_n$ be the counting measure of $pr(\mathcal Q_n)$, i.e.,
\[\mu_n = \sum_{\underline v\in \mathcal Q_n} \delta_{pr(v)}~.\]
We first prove the following elementary counting result:

\begin{proposition} \label{p:Equidistribution sublattices}
The sequence of measures $n^{-\frac{d}{2}}\mu_n$ converges weakly to the smooth measure $\omega$.
\end{proposition}

\begin{proof}
Let $f:\Omega_1 \to \R$ be a continuous function with compact support which we extend on $\Omega$ by setting $f(t\underline v)=f(\underline v)$ for all $t\in \R_{>0}$ and all $\underline v\in \Omega_1$. We have
\begin{align*}
n^{-\frac{d}{2}}\mu_n(f)& =n^{-\frac{d}{2}} \sum_{\underline v\in V_\Z^r \cap \Omega\mid h(\underline v)\leq n} f(\underline v)\\
&=n^{-\frac{d}{2}} \sum_{\underline v\in n^{-\frac{1}{2r}}V_\Z^r \cap \Omega\mid h(\underline v)\leq 1} f(\underline v) \quad \textrm{by homogeneity of $f$}\\
& \underset{n\to +\infty}{\longrightarrow} \int_{\Omega_{\leq 1}} f = \omega(f)\quad \textrm{by Riemann summation.}
\end{align*}
\end{proof}

The group $G \times \SL(r,\R)$ acts transitively on $\Omega_1$ by
\[(g_1,g_2)\cdot \underline v = g_1 \cdot \underline v \cdot g_2^{-1}\]
and preserves the measure $\omega$. The restriction of this action to $\SL(r,\R)$ is proper and the quotient $\SL(r,\R) \backslash \Omega_1$ is the space $\mathcal P$ or positive definite $r$-subspaces of $V_\R$.

Now, the subgroup $\SL(r,\Z)$ preserves the set $\mathcal Q_n$ and acts properly discontinuously on $\Omega_1$ so that the quotient set $\bar {\mathcal Q}_n = \SL(r,\Z) \backslash \mathcal Q_n$ still equidistributes in $\SL(r,\Z)\backslash \Omega_1$. Let us consider the projection
\[\pi: \SL(r,\Z)\backslash \Omega_1 \to \SL(r,\R) \backslash \Omega_1 = \mathcal P~.\]
We still denote $\omega$ the volume form induced on $\SL(r,\Z)\backslash \Omega_1$.
The push-forward measure $\pi_*\omega$ is $G$-invariant (since $\omega$ is $G$-invariant and $\pi$ is $G$-equivariant), non-zero, and locally finite since $\SL(n,\Z)\backslash \SL(n,\R)$ has finite volume. We hence deduce from Proposition \ref{p:Equidistribution sublattices}:

\begin{corollaire}
Define the measure
\[\nu_n = \sum_{\underline v\in \bar Q_n} \delta_{\pi(\underline v)}~.\]
Then
\[n^{-\frac{d}{2}} \nu_n \rightharpoonup \lambda \omega_{G/H}\]
for some $\lambda\neq 0$.
\end{corollaire}

\begin{remarque}
The multiplicative constant $\lambda$ could be computed in terms of the volume of $\SL(n,\R)/\SL(n,\Z)$.
\end{remarque}

The measure $\nu_n$, however, is not the counting measure of $\mathcal P_n$. To be more precise, note that $\bar {\mathcal Q}_n$ is the set positive definite sublattices of $V_\Z$ of discriminant $\leq n$ and $\pi$ maps $\Lambda\in \bar {\mathcal Q}_n$ to $\Lambda\otimes \R$. Therefore, we have
\[\nu_n= \sum_{W\in \mathcal P_n} m_n(W) \delta_W\]
where
\[m_n(W) = |\{\Lambda\subset W\cap V_\Z\mid h(W)  [W\cap V_Z:\Lambda] \leq n\}|~.\]
In other words, $\nu_n$ counts a rational subspace $W$ with a weight equal to the number of sublattices of $W\cap V_\Z$ with discriminant $\leq n$. 

Let $\nu_n'$ be the counting measure of $\mathcal P_n$. To relate $\nu_n$ and $\nu_n'$, let us introduce
\[b_k= |\{\Lambda \subset \Z^r\mid [\Z^r: \Lambda]= k\}|~.\]
We have the following estimate on $b_k$:
\begin{proposition}
\[b_k \ll k^r~.\]
\end{proposition}
\begin{proof}
We will prove a sharper estimate. Consider the zeta function which converges for large $s$: 
\[\zeta_{\Z^r}(s)=\sum_{k\geq 1}\frac{b_k}{k^s}=\sum_{g}|\det(g)|^{-s}\]
where $g$ runs through $\mathrm{GL}_{r}(\Z)\backslash\left(\mathrm{M}_{r}(\Z)\cap \mathrm{GL}_r(\Q)\right)$. By \cite[Equation (15.10)]{lubotzky}, we have the equality 
\[\zeta_{\Z^r}(s)=\prod_{i=0}^{r-1}\zeta(s-i).\]
Hence by identifying the coefficients, we get
\begin{align*}
    b_k&=\sum_{\underset{k_0\cdots k_{r-1}=k}{(k_0,\cdots,k_{r-1})} }k_0\cdots k_{r-1}^{r-1}\\
    &\leq k^{r-1}|\{(k_0,\cdots,k_{r-1}), k_0\cdots k_{r-1}=k\}|\\
    &\ll_\epsilon k^{r-1+\epsilon}
\end{align*}
for every $\epsilon >0$. Hence the result.
\end{proof}

We now have
\begin{align*}
    \nu_n &= \sum_{P\in \mathcal P_n} \left(\sum_{k\leq \sqrt \frac{n}{h(P)}} b_k\right) \delta_P\\
    &= \sum_{k\leq \sqrt{n}} \sum_{h(P)\leq \frac{n}{k^2}} b_k \delta_P
\end{align*}
and we conclude that
\begin{equation} \label{eq:nu_n nu'_n}
\nu_n = \sum_{k\leq \sqrt{n}} b_k \nu'_{\lfloor\frac{n}{k^2}\rfloor}~.
\end{equation}

Set 
\[\alpha = \sum_{k} \frac{b_k}{k^{d}}~.\]

\begin{proposition}
The measure $\nu'_n$ converges weakly to
\[\frac{\lambda}{\alpha}\omega_{G/H}~.\]
\end{proposition}

\begin{proof}
Remark first that, under the hypothesis $1\leq r \leq p = d-2q \leq d-2$, we have 
\[\frac{b_k}{k^{d}} \ll \frac{1}{k^{2}}~,\]
hence
\[\alpha \leq \zeta(2) < 2~.\]

Let $f$ be a continuous function with compact support on $G/H$. Set
\[s_n = n^{-\frac{d}{2}}\frac{\nu_n(f)}{\lambda\int_{G/H} f \omega_{G/H}}\] and 
\[s'_n = n^{-\frac{d}{2}}\frac{\nu'_n(f)}{\lambda\int_{G/H} f \omega_{G/H}}~.\]
By \eqref{eq:nu_n nu'_n}, we have
\begin{equation} \label{eq:s_n s'_n}s_n = \sum_{k\geq 1} n^{-\frac{d}{2}}\left \lfloor\frac{n}{k^2}\right \rfloor^\frac{d}{2} b_k s'_{\left \lfloor\frac{n}{k^2}\right \rfloor}~.\end{equation}
By Proposition \ref{p:Equidistribution sublattices}, we have 
\[s_n\underset{n\to +\infty}{\longrightarrow} 1~.\]
In particular, $(s_n)$ is bounded by a constant $c$. Since $s'_n \leq s_n$ the sequence $s'_n$ is also bounded and, by convergence of the series
\[\sum_{k} \frac{b_k}{k^{d}}~.\]
we can find for all $\epsilon >0$ some $k_0$ such that 
\[\sum_{k\geq k_0} n^{-\frac{d}{2}}\left \lfloor\frac{n}{k^2}\right \rfloor^\frac{d}{2} b_k s'_{\left \lfloor\frac{n}{k^2}\right \rfloor}\leq \epsilon\]
for all $n$.

Set $m= \liminf s'_n$ and $M=\limsup s'_n$, and let $n_i$ be a subsequence such that $s'_{n_i} \underset{n\to +\infty}{\longrightarrow} m$. For all $2\leq k \leq k_0$ we have
\[\limsup_{i\to +\infty} n_i^{-\frac{d}{2}}\left \lfloor\frac{n_i}{k^2}\right \rfloor^\frac{d}{2} s'_{\left \lfloor\frac{n_i}{k^2}\right \rfloor} =
\leq \frac{b_k}{k^{d}} M~.\]
Hence, taking the limsup of \ref{eq:s_n s'_n} along $n_i$, we get
\[1 \leq m +\sum_{k=2}^{k_0} \frac{b_k}{k^{rd}}M + \epsilon = m + (\alpha-1)M + \epsilon~.\]
Similarly taking the liminf along a subsequence $n_i$ such that $s'_{n_i} \underset{n\to +\infty}{\longrightarrow} \to M$, we obtain
\[1\geq M + (\alpha-1)m~.\]
Combining the two, we get
\[ M + (\alpha - 1)m \leq m + (\alpha-1)M + \epsilon~,\]
hence
\[M-m \leq \frac{\epsilon}{2-\alpha}\]
since $\alpha< 2$.

Since this is true for all $\epsilon >0$, we conclude that $M=m$. Hence $s'_n$ converges to $m=M$. Taking the limit in \eqref{eq:s_n s'_n} gives $1 = \alpha m$, and we conclude that \[s'_{n} \underset{n\to +\infty}{\longrightarrow} \frac{1}{\alpha}~.\]

Going back to the definition of $s'_n$ we have proved that
\[n^{-\frac{d}{2}}\nu'_n(f) \underset{n\to +\infty}{\longrightarrow} \frac{\lambda}{\alpha} \int_{G/H} f\omega_{G/H}~.\]
\end{proof}

\subsubsection{Equidistribution along intersection matrices}\label{S:levelequidistribution}
To prove the second version of the equidistribution theorem which yields \Cref{maink3}, we can restrict to matrices $M$ which are primitively represented by $(V_\Z,B)$, i.e., for which there exists $(\lambda_1,\cdots,\lambda_r)\in V_\Z^r$ generating a primitive sublattice of $\V_\Z$ and with intersection matrix $M$. For simplicity, if $\underline{\lambda}\in V_\R^r$, we denote by $I(\underline{\lambda})$ the intersection matrix $\left(B(\lambda_i.\lambda_j)\right)_{1\leq i,j\leq r}$. Let \[V_{\R,\mathrm{I}_r}^r=\{\underline{\lambda}=(\lambda_i)_{1\leq i\leq r},\, I(\underline{\lambda})=\mathrm{I}_r\}.\]

Then $V_{\R,\mathrm{I}_r}^r$ is an affine homogeneous variety under the action of the group $G=\mathrm{O}(V_\R,B)\simeq \mathrm{O}(p,2q)$ and letting $H\simeq \mathrm{O}(p-r,2q)$ be the stabilizer of a point $\underline{\lambda}_0\in V_{\R,\mathrm{I}_r}^r$, then $V_{\R,\mathrm{I}_r}^r\simeq G/H$. When $r<p$ and $q\geq 1$ the group $H$ is simple without compact factors, so that Ratner theory can be applied as explained in \ref{eskin-oh}. Finally, $H$ is not contained in any proper parabolic subgroup of $G$, so that  sequences of closed $H$-orbits of $\Gamma \backslash G$ do not have loss of mass.

There is a right action of a $r\times r$-matrix $A=(a_{i,j})$ on a vector $\underline{u}=\left(u_1,\cdots,u_r\right)\in V_\R^r$ given by matrix product:
\[\underline{u}\cdot A=\left(\sum_{j=1}^ra_{1,j}u_1,\ldots ,\sum_{j=1}^ra_{r,j}u_j\right).\]
Notice that this action commutes with the diagonal action of $\mathrm{GL}(V_\R)$ and that the components of $\underline{u}\cdot A$ span a subspace of the vector space spanned by components of $\underline{u}$ in $V_\R$, and they are equal if $A$ is invertible. Their intersection matrix are related by 
\[I(\underline{u})= {^tA}I(\underline{u})A~.\]
Let $M$ be a positive definite integral matrix of size $r$ and let \[W_M=\{\underline{\lambda}=(\lambda_i)_{1\leq i\leq r}\in V_\Z^r,\, I(\underline \lambda)=M\}.\] 

In order to study equidistribution of $W_M$, it is natural to first project it to $V_{\R,\mathrm{I}_r}^r$. We have thus a $G$-equivariant projection map 
\begin{align*}
    pr:W_M&\rightarrow V_{\R,\mathrm{I}_r}^r\\
    \underline{\lambda}&\mapsto \underline{\lambda}\cdot \sqrt{M}^{-1},
\end{align*}
where $\sqrt{M}$ is the unique positive definite matrix such that $\sqrt{M}^2=M$. 

By a theorem of Borel and Harish-Chandra \cite[Theorem 6.9]{borelharish}, the projection $pr(W_M)$ is a finite union of discrete $\Gamma$-orbits of $G/H$, which thus corresponds to a finite union of closed $H$-orbits $\mathcal O_M \subset \Gamma \backslash G$. The volume of $\mathcal O_M$ is finite by Borel--Harish-Chandra's theorem \cite[theorem 9.4]{borelharish} since $H$ is semi-simple, and the following lemma gives an estimate for its volume:

\begin{lemma}{\ }\label{totalvolume}
\begin{enumerate}
    \item Let $M$ be a positive definite matrix of rank $r\leq p$ represented by the lattice $(V_\Z,B)$. Then there exists $c>0$ depending only on $(V_\Z,B)$ and $r$ such that:
\[a(M)\equaldef\Vol(\mathcal{O}_M)=c\det(M)^{\frac{p+2q-r-1}{2}}\prod_{a\,\textrm{prime}}\beta_a(M),\]
where for a prime number $a$, the local density at $a$ is expressed as \[\beta_{a}(M)\equaldef\lim_{s\rightarrow \infty}a^{-s\left(r(p+2q-\frac{r+1}{2})\right)}|\{\underline{\lambda}\in V_\Z^r/a^sV_\Z^r, I(\underline{\lambda})=M\}|~.\]

\item If $(M_n)_{n\in \N}$ is a sequence of positive definite matrices primitively represented by $(V_\Z,B)$, then \[a(M_n) \underset{n\to \infty}{\geq}\det(M_n)^{\frac{p+2q-r-1}{2}-\epsilon}\]
for any $\epsilon>0$. In particular, $a(M_n)$ goes to $+\infty$, as $\det(M_n)$ goes to $+\infty$.

\item If moreover $r\leq \frac{p+2q-3}{2}$, then 
\[a(M_n) \underset{n\to \infty}{\asymp}\det(M_n)^{\frac{p+2q-r-1}{2}}~.\]
\end{enumerate}
\end{lemma}

\begin{proof}
The first assertion is simply the Siegel-Weil formula, which is valid in this setting by \cite{weilsiegel}. To prove the second statement, we need to find a lower bound on the growth of the product of the local densities $\beta_a(M_n)$ assuming that $M_n$ is primitively represented by $(V_\Z^r,B)$. Let $n\in \N$ and let $P(n)$ be the set of odd primes $a$ which are coprime to $\det(M_n).\det(V_\Z)$. By \cite[Proposition 5.6.2(ii)]{kitaoka}, there exists two positive numbers $c_1,c_2$ depending only on $V_\Z$ such that 
\[c_1<\prod_{a\in P(n)}\beta_a(M_n)<c_2.\]

If $r\leq \frac{p+2q-3}{2}$, then by Corollary 5.6.2 {\it loc.cit.}, the above estimate on the product is true for $a$ ranging over all primes, proving the third statement.

Otherwise, since we assumed that $M_n$ is represented by a sublattice of $V_\Z$ which is primitive \footnote{Even weaker assumption such as locally bounded imprimitivity is enough, see \cite[Theorem 5.6.5]{kitaoka}.}, then \cite[Theorem 5.6.5,(a)]{kitaoka} yields that there exists a constant $c_3<1$ such that for any prime $a$ dividing $\det(V_\Z)\cdot\det(M_n)$ we have
\[\beta_a(M_n)\geq c_3 \]
Since the number of prime divisors of $\det(M_n)$ is $O(\frac{\log(\det(M_n))}{\log\log(\det(M_n))})$, we obtain that
\[ \Vol(\mathcal O_{M_n}) \geq  c c_1 \det(M_n)^{\frac{p+2q-r-1}{2}} c_3^\frac{c_4\log(\det(M_n))}{\log\log(\det(M_n))} = \det(M_n)^{\frac{p+2q-r-1}{2} + o(1)}~.\]
\end{proof}

For a positive definite integral matrix $M$, let $\mu_1(M)$ be the square root of the smallest non-zero integer represented by $M$.
\begin{theorem}\label{T:levelequi}
Let $(M_n)_{n\in N}$ be a sequence of positive definite matrices primitively represented by $(V_\Z,B)$ and such that $\mu_1(M_n)\rightarrow\infty$ as $n\rightarrow \infty$. Then the sequence of subsets $\{pr(W_{M_n}),n\in\N\}$ is equidistributed in $V_{\R,\mathrm{I}_r}^r$ in the sense of \Cref{eskin-oh}.

\end{theorem}

\begin{proof}
Note first that, since $M$ is positive definite, we have $\det(M)\geq c \mu_1(M)^2$ where $c$ depends only on the rank of $M$, see \cite[Equation (5)]{EK}. Hence $\det(M)$ goes to $\infty$ and so does $\Vol(\mathcal O_{M_n})$ by Lemma \ref{totalvolume}.

To prove the equidistribution, we apply \Cref{eskin-oh}. Since $H$ is not contained in a proper parabolic subgroup, the sequence $\mathcal O_{M_n}$ has no loss of mass, and we need to prove that it is non-focused see (\Cref{d:NonFocused}). We are in the easy situation where any sequence of $\Gamma$-orbits $\Gamma \cdot \underline \lambda_n \subset  pr(W_{{M_n}})$ is non-focused.

To prove this, assume by contradiction that, up to taking a subsequence, there exists a proper subgroup $H'$ of $G$ defined over $\Q$, an element $g\in G$ such that $gH^0g^{-1}\subset H'$ and a sequence $\underline \lambda_{n} \in \mathcal E_{M_{n}}$
such that $pr(\underline \lambda_n)\in H' g Z(H^0)\cdot \underline pr(\lambda_0)\lambda_0$. 

Set $V_n = \mathrm{Span}_\R(\underline \lambda_n)$ and let $H_n \subset G$ be the subgroup fixing $V_n$. Then $H_n$ is conjugate to $H$ and contained in $H'$ for all $n$ by assumption on $\underline \lambda_n$. In particular, by \Cref{maxsubgroup} below, $H'$ preserves a rational subspace $W$ contained in $V_0$. Hence every $H_n$ preserves $~W$. Since the action of $H_N$ on $V_n^\perp$ is irreducible, we deduce that $W\subset V_n$ for all $n$.

Since $W$ is rational, it intersects $V_\Z$ in a lattice which is contained in $\mathrm{Span}_\Z(\underline \lambda_n)$ for all $n$ since $\mathrm{Span}_\Z(\underline \lambda_n)$ is primitive. This contradicts the assumption that $\mu_1(M_n)\to +\infty$.

\end{proof}

\begin{lemma}\label{maxsubgroup}
Let $V_0$ be a positive definite rational subspace of $V_\Q$, let $H_0$ be the subgroup of $G$ fixing $V_0$ and let $H$ be a proper connected subgroup of $G$ defined over $\Q$ and containing $H_0$. Then $H$ leaves invariant a rational subspace of $V_0$.
\end{lemma}

\begin{proof}
As a representation of $H_0$, the Lie algebra $\fg$ decomposes orthogonally with respect to the Killing form as
\[\fg=\fh_0 \oplus \so(V_0) \oplus \fp\]
where $\fp = \{u\in \fg \mid u(V_0)\subset V_0^\perp\}$ is isomorphic to $\Hom(V_0,V_0^\perp)$. Note that the representation of $H_0$ on $V_0^\perp$ is irreducible and $H_0$ acts trivially on $V_0$, hence also on $\so(V_0)$.

Since $H$ contains $H_0$, its Lie algebra $\fh$ is a subrepresentation of $\fg$ and thus decomposes as
\[\fh_0 \oplus \fk \oplus \fp'~,\]
where $\fk$ is a Lie subalgebra of $\so(V_0)$ and $\fp'$ is a subrepresentation of $\fp \simeq \Hom(V_0,V_0^\perp)$. By elementary representation theory, there exists a subspace $W\subset V_0$ such that 
\[\fp' = \{u\in \fp \mid u_{\vert W} = 0\}~.\] 
We have
\[W= \{x\in V_0 \mid u(x) \in V_0 \textrm{ for all } u\in \fh\}~,\]
in particular, $W$ is a rational subspace since $\fh$ and $V_0$ are defined over $\Q$.

We claim that $\fk$ preserves $W$. Indeed, assume by contradiction that there exists $u\in \fk$ and $x\in W$ such that $u(x)\in V_0\setminus W$. Then there is $v\in \fp'$ such that $vu(x)\notin V_0$. Since $v(x)=0$, we obtain that
\[[u,v](x) = uv(x)\notin V_0~,\]
contradicting $x\in W$.

In conclusion, the Lie algebra $\fh$ preserves $W\subset V_0$. If $W$ were trivial, then we would have $\fh \supset \fp\oplus \fh_0$, hence $\fh = \fg$ since $[\fp,\fp]\supset \so(V_0)$. Since $H$ is a proper subgroup, $W$ is non-trivial.
\end{proof}

\subsubsection{Proof of Theorems \ref{maink3-alt} and \ref{maink3}} \label{sss:ProofEquidistributionLN}

Gathering together the results of the previous sections, we can finally prove our equidistribution theorems for refined Noether--Lefschetz loci in weight $2$. Let us first state them more precisely.

Let $\{\mathbb{V}_\Z,\mathcal{F}^\bullet\mathcal{V},B\}$ be a $\Z$-PVHS of weight $2$ over a complex manifold $S$ of dimension $rq$ as in \Cref{maink3-alt}. Let $s\in S$ and let $P\subseteq \Hodge(s)$ be a subspace of rank $r$. Then the pair $(s,P)$ is a transverse intersection point of $S$ with a $H$-orbit, where $H$ is the stabilizer of a positive definite subspace of $V_\R$ as in \Cref{S:averageequidistribution}, if $s$ does not admit first order deformations such that $P$ still embeds in the group of Hodge classes. Similarly, if $(\lambda_1,\cdots,\lambda_r)\in \Hodge(s)^r$ have intersection matrix $M$, then the tuple $(s,\lambda_1,\cdots,\lambda_r)$ is a transverse intersection point with a $H$-orbit, $H$ being now the stabilizer of an orthonormal $r$-tuple as in \cref{S:levelequidistribution}, if $s$ does not admit first order deformations such that $\lambda_1,\cdots,\lambda_r$ all remain Hodge classes.
\medskip

We can now prove the main theorems in \Cref{section:hodgeloci}. Notations and hypothesis are as in \Cref{maink3-alt}.
\begin{theorem}\label{correct:maink3-alt}
There exists a constant $\lambda>0$ such that, for every relatively compact open subset $\Omega\subset S$ with boundary of measure $0$, we have
\begin{align*}
    n^{-\frac{p+2q}{2}}\left \vert \{(s,P)|\, s\in \Omega, P\subseteq \Hodge(s), \rank(P)=r,  (s,P)\textrm{ regular}, \mathrm{disc}(P)\leq n \} \right \vert\\ \underset{n\to+\infty}{\longrightarrow} \lambda \int_\Omega c_q(\mathcal F^2\cV)^r~,
\end{align*}
where $c_q$ denotes the $q$\textsuperscript{th} \emph{Chern form} of the bundle $\mathcal F^2\cV$ endowed with the Hodge metric.
\end{theorem}
\begin{proof}
We use the notations from \Cref{S:averageequidistribution}. By \Cref{T:averageequi}, the sequence $(\mathcal{P}_n)_{n\in\N}$ is equidistributed in $G/H$. We are now in the setting of \Cref{t:interequi}. By \Cref{pullpushweight2}, the pull-push form $\pi_*p^* \omega_{G/H}$ is equal to $\Vol(G_U/H_U).c_q(\cF^2\cV)^r$, where $c_q(\cF^2\cV)$ is the $q^{th}$-Chern form of $\cF^2\cV$. We can hence apply \Cref{t:interequi} to deduce \Cref{maink3}.
\end{proof}
Notations and hypothesis  are now as in \Cref{maink3}.
\begin{theorem}
For every relatively compact open subset $\Omega \subset S$ with boundary of measure $0$, we have:
    \begin{align*}
        \frac{1}{a(M_n)}|\{(s,\lambda_1,\cdots,\lambda_r))\in \Omega\times \V^r_{\Z,s}\, \textrm{regular tuple },(B(\lambda_i.\lambda_j))_{1\leq i,j\leq r}=M_n, \lambda_i \in \mathrm{Hdg}(s)\}|\\
        \underset{n\rightarrow \infty}{\longrightarrow} \Vol(G_U/H_U) \int_{\Omega}c_q(\mathcal{F}^2\cV)^r.
    \end{align*}
\end{theorem}

\begin{proof}
We use the notations from \Cref{S:levelequidistribution}. By \Cref{T:levelequi}, the sequence $(pr(W_{M_n}))_{n\in\N}$ is equidistributed in $G/H$. We are now again in the setting of \Cref{t:interequi}. Here $H$ is the stabilizer of an orthonormal $r$-tuple of vectors, and if we denote by $H'$ the stabilizer of the rank $r$ subspace they generate in $\V_\R$, then $H'/H$ is compact and one easily checks then that $\pi_*p^* \omega_{G/H}=\Vol(H'/H)\pi_*p^* \omega_{G/H'}$.

By \Cref{pullpushweight2}, the pull-push form $\pi_*p^*=\omega_{G/H'}$ is equal to $\Vol(G_U/H'_U)\cdot c_q(\cF^2\cV)^r$, where $c_q(\cF^2\cV)$ is the $q^{th}$-Chern form of $\cF^2\cV$. Moreover, 
\[\Vol(G_U/H'_U)=\frac{\Vol(G_U/H_U)}{\Vol(H'_U/H_U)}=\frac{\Vol(G_U/H_U)}{\Vol(H'/H)}.\]
Hence $\pi_*p^* \omega_{G/H}=\Vol(G_U/H_U)c_q(\cF^2\cV)^r$. We can hence apply \Cref{t:interequi} to deduce \Cref{maink3}. Indeed, one can easily see again that regular points in our definition above correspond to transverse intersection points defined there.
\end{proof}

Finally we mention briefly how to prove \Cref{kudlamillson}.
\begin{proof}[Proof of \Cref{kudlamillson}]
Combining \Cref{cohomologyequidistribution}, \Cref{Sarnaklemma} and Proposition 6.8, we get \Cref{kudlamillson}. 
\end{proof}

\subsection{Equidistribution of families of CM points in Shimura varieties}\label{shimura}

In this section, we use \Cref{main} to study the  equidistribution of transverse intersection loci of Hecke correspondences on Shimura varieties and deduce the equidistribution of some families of CM points in average. We recall first the definition of a CM Hodge structure, see \cite[V]{griffiths}.

A CM field is a totally imaginary number field which is a quadratic extension of a totally real number field. A CM algebra is a finite product of CM number fields.
\begin{definition}\label{CM}
Let $(V,B,F^\bullet)$ be a pure polarized integral Hodge structure and let $d=\rank_\Z V$. We say that $(V,B,F^\bullet)$ has \emph{Complex Multiplication} (``CM'' for short) if one of the following equivalent conditions hold:
\begin{enumerate}
    \item Its algebraic Mumford--Tate group $\mathrm{MT}_{\varphi}$ is a torus ;
    \item  The ring $\mathrm{End}_\Q(V,F^\bullet)$ contains an étale CM sub-algebra of dimension $2d$.
\end{enumerate} 
\end{definition}

We refer to \cite[(IV.B.1)]{griffithsgreenkerrarticle} for the equivalence in the definition above.
\begin{example}{\ }

\begin{enumerate}
    \item Let $(A,\lambda)$ be a polarized complex abelian variety and let  $\mathrm{End}(A)_\Q=\mathrm{End}(A)\otimes \Q$. Recall that $A$ has complex multiplication if $\mathrm{End}_\Q(A)$ contains an étale sub-algebra of degree $2\dim(A)$ over $\Q$. This is equivalent to the polarized Hodge structure $(H^1(X,\Z),\psi)$ being CM in the sense of \Cref{CM}. 
    \item Let $(X,\ell)$ be a complex polarized K3 surface and let $T(X)$ be the transcendental lattice  of $X$, i.e., the orthogonal complement of $\mathrm{Pic}(X)$ inside $H^2(X,\Z)$ with respect to the Poincaré form. Then $X$ has CM if $E\equaldef\mathrm{End}(T(X))_\Q$ is a CM field and $T(X)_\Q$ is of dimension $1$ over $E$. If $PH^2(X,\Z)$ is the primitive cohomology of $X$, then this is equivalent to $PH^2(X,\Z)$ being CM in the sense of \Cref{CM}.
\end{enumerate}
\end{example}


Let $(\widetilde{G},\mathcal{D})$ be a Shimura datum, see \cite{deligneshimura,deligneshimura2},  and let $\D^+$ be a connected component of $\D$ ; $\D⁺$ is a $\widetilde{G}^{ad}(\R)^+$-conjugacy class of a morphisms $h^{ad}:\mathbb{S}\rightarrow \widetilde{G}^{ad}_\R$ and it is a Hermitian symmetric domain. Let $K_\infty$ be the stabilizer of $h^{ad}(i)$ in $\widetilde{G}^{ad}(\R)^+$. If $K_{\infty,+}$ is its preimage by the adjoint map, then we have an isomorphism $\D^+=\widetilde{G}(\R)_+/K_{\infty,+}\simeq \widetilde{G}^{ad}(\R)^+/K_\infty$. Let $\Gamma \subset \widetilde{G}(\Q)$ be an arithmetic subgroup and let $X=\Gamma\backslash \D^+$. Then $X$ is a \emph{(connected) Shimura variety}.

\begin{definition}
Let $g\in \widetilde{G}^{ad}(\Q)$ and $\Gamma_g=g^{-1}\Gamma g\cap \Gamma$. The Hecke correspondence $\mathcal{C}_g\hookrightarrow X\times X$ is the image of $\Gamma_g\backslash \D$ by the embedding 
\begin{align*}
    \Gamma_g\backslash \D&\hookrightarrow X\times X\\
    [x]&\mapsto ([x],[gx]).
\end{align*}
\end{definition}

If $g=1$ is the identity of $\widetilde{G}(\mathbb{Q})$, then $\mathcal{C}_1$ is simply the diagonal embedding of $X$ in $X\times X$. 

\begin{proposition}\label{p:regularity}
For $f\in \tilde G(\R)$, the following properties are equivalent:
\begin{enumerate}
    \item[(i)] $f$ has a unique fixed point in $\mathcal D$,
    \item[(ii)] The centralizer of $f$ is compact in $\tilde G^{ad}(\R)$,
    \item[(iii)] The intersection of the graph of $f$ and the diagonal in $\D\times \D$ is transversal and non-empty.
\end{enumerate}
If $f$ satisfies those properties, we say that $f$ is \emph{regular}.
\end{proposition}
\begin{proof}
$(1)\Rightarrow(2)$: let $x$ be the unique fixed point of $f$, then $f$ is contained in the stabilizer of $x$ in $\widetilde{G}^{ad}(\R)$ which is a compact subgroup. Moreover, for any $g\in Z(f)$, $g\cdot x$ is also a fixed point for $f$, hence equal to $x$ and thus $Z(f)\subseteq K$.

$(2)\Rightarrow(3)$: Since $Z(f)$ is compact, it is contained in a maximal compact subgroup $K$ of $\widetilde{G}^{ad}(\R)$. Hence $f$ fixes a point $x$ and the differential $df_x$ on $T_x\mathcal{D}$
 identifies to the action of $Ad(f)$ on $\mathfrak{p}$, the orthogonal complement of $\mathfrak{p}$ in $\widetilde{\mathfrak{g}}^{ad}$ with respect to the Killing form. Then $Ad(f)$ does not have $1$ as eigenvalue in $\mathfrak{p}$, as $Z(f)\subseteq K$. This will hold true at any fixed point $f$ in $\mathcal{D}$. Let $(x,x)$ be an intersection point of the graph $\mathcal{C}_f$ of $f$ and the diagonal $\Delta$ in $\mathcal{D}\times\mathcal{D}$, then 
 the tangent spaces of $\mathcal{C}_f$ and $\Delta$ at $(x,x)=(x,f\cdot x)$ inside $T_{(x,x)}(\mathcal{D}\times \mathcal{D})$ are given respectively by $\{(X,X)|\,X\in \mathfrak{p}\}$ and $\{(X,Ad(f)\cdot X)|\, X\in \mathfrak{p}\}$. Their sum is equal to $\mathfrak{p}\oplus{\mathfrak{p}}$ if and only if their intersection is zero, which is true as $1$ is not an eigenvalue of $Ad(f)$ in $\mathfrak{p}$.

 $(3)\Rightarrow(1)$: if the intersection is transverse, then by the previous computation, for any fixed $x$ point of $f$, the eigenvalues of $df_x$ in $\mathfrak{p}$ are different from $1$. If $f$ fixes another point $y\in\mathcal{D}$, then it fixes the geodesic line $\gamma:\R\rightarrow \D$ linking $\gamma(0)=x$ to $y$ and hence acts trivially on this line. Hence $df_x(\dot{\gamma}(0))=\dot{\gamma}(0)$ which is a contradiction. Thus, $f$ has a unique fixed point.
\end{proof}

\begin{lemma}\label{transverselocusCM}
A point $[x] \in X$ is CM if 
and only if there exists $g\in \tilde G(\Q)$ such that $\mathcal C_1$ and $\mathcal C_g$ intersect transversally at $x$.
\end{lemma}

\begin{proof}
The transverse intersection locus of $\mathcal{C}_g$ and $\mathcal{C}_1$ inside $X\times X$ is necessarily of dimension $0$ by dimension count. Let $[x]$ be a point in this intersection and $h_x:\mathbb{S}\rightarrow \widetilde{G}(\R)$ a lift to $\D$. Then there exists $y \in \D$, $\gamma_1,\gamma_2\in \Gamma$ such that $x=\gamma_1\cdot y$ and $x=\gamma_2g\cdot y$. Hence $x=\gamma_2g\gamma_{1}^{-1}\cdot x$ which implies that $MT(x)\subseteq Z(\gamma_2g\gamma_{1}^{-1})$. Since the intersection is transverse at $x$, then by the Lemma above, $\gamma_2g\gamma_{1}^{-1}$ is regular and contained in $K$ by regularity. Hence by \cite[IV.B.1]{griffithsgreenkerrarticle}, $x$ is a CM point. 

Conversely, let $x\in\mathcal{D}$ be a CM point. Then $L\equaldef Z(MT(x))^{\circ}$, the connected component of the Mumford Tate group in its centralizer, is defined over $\Q$. Then $ L(\R)\subset K$ and the function $u:f\mapsto \det(Ad(f)_{\vert \mathfrak p}-\mathrm{Id}_{\mathfrak{p}})$ is well defined and does not vanish as $u(h_x(i))\neq 0$. By \cite[\S 18, Corollary 18.3]{borelbook}, $L(\Q)$ is Zariski dense in $L(\R)$, hence there exists an element $f\in L(\Q)$ which is regular and $MT(x)\subseteq Z(f)$. Hence $x$ is a transverse intersection point of $\mathcal{C}_f$ and $\mathcal{C}_1$. 

\end{proof}

For $g\in \widetilde{G}^{ad}(\mathbb{Q})$, let $\deg(g)=[\Gamma:\Gamma_g]$. More generally, if $V\subset\widetilde{G}^{ad}(\mathbb{Q})$ is a $\Gamma$-double class with finitely many left $\Gamma$-orbits, we let $\deg(V)$ be the number of distinct left $\Gamma$-orbits, in particular $\deg(g)=\deg(\Gamma g\Gamma)$. If we set $G\equaldef \widetilde{G}^{ad}(\R)^+\times\widetilde{G}^{ad}(\R)^+$ and $H=\widetilde{G}^{ad}(\R)^+$ embedded diagonally in $G$, then $G/H\simeq \widetilde{G}^{ad}(\R)$ and we are in the situation of \Cref{ratnersetup}. Then we denote by $\mathcal{O}\subset (\Gamma\times\Gamma)\backslash G$ the corresponding finite union of closed  $H$-orbits and by $\mathcal{C}_{\mathcal{O}}\hookrightarrow X\times X$ the associated Hecke correspondence. 

\begin{theorem}\label{equidistributionCM}
Let $X$ be a Shimura variety associated to  a Shimura datum $(\widetilde{G},\mathcal{D})$ such that $\widetilde{G}^{ad}$ is connected and $\Q$-simple. Let $\left(V_n\right)_{n\in\N}$ be a sequence of $\Gamma$-double classes in $\widetilde{G}(\Q)$ such that $\deg(V_n)\rightarrow \infty$. Then for every $\Omega\subset X$  open relatively compact subset with zero measure boundary 
\[|\{(x,f)|\, x\in \Omega, f\in V_n, MT(x) \subset Z(f), f \textrm{ regular}\}|\underset{n\rightarrow \infty}{\sim}\frac{\deg(V_n)\cdot \chi(\widehat{\D})}{\Vol(\widehat{\D})}\int_{\Omega}\omega_\D~.\]
\end{theorem}
\begin{proof}
Let $H=\widetilde{G}^{ad}\hookrightarrow G\equaldef \widetilde{G}^{ad}\times \widetilde{G}^{ad}$ and by assumption $\widetilde{G}^{ad}$ is simple. Then the quotient $G/H$ is isomorphic to $H$ via the map $p:(x,y)\mapsto yx^{-1}$. The preimage by $p$ of an element $a\in G$ is equal to $(1,a).H\hookrightarrow G$.

Then the sequence of $\Gamma$-double class $(V_n)_{n\N}$ are equidistributed in $G/H$. This result has been proved by \cite{clozelohullmo} in the following cases: $\widetilde{G}$ is connected, almost simple simply connected and $\textrm{rank}_\Q(\widetilde{G})\geq 1$ (see Theorem 1.6 {\it loc. cit.}) and for $G=\mathrm{GSp}_{2g}$ (Remark (3) page 332 {\it loc. cit.}). More generally, Eskin and Oh \cite{eskinoh1} proved this result for any $\widetilde{G}$ connected and simple over $\Q$.\footnote{The simplification comes at a cost of not having an error term.} Hence we are in the setting of \Cref{main}.

By \Cref{diagonalrestriction}, the restriction of the form $\pi_*p^*\omega_{G/H}$ to $\D$ is equal to $\frac{\chi(\widehat{\D})}{\Vol(\widehat{\D})}\cdot \omega_{\D}$. Hence by \Cref{generictransversality}, $\Gamma\backslash \D$ is generically transverse to $H$-orbits and by \Cref{main}, the transverse intersection locus of $\mathcal{C}_{V_n}$ with $X=\mathcal{C}_1$ is equidistributed in $X$ with respect to $\frac{\chi(\widehat{\D})}{\Vol(\widehat{\D})}\cdot \omega_{\D}$ as $n\rightarrow \infty$. By \Cref{transverselocusCM}, this transverse locus is formed by elements $x$ where $x$ is a CM point with Mumford--Tate group $MT(x)\subset Z(f)$ where $f$ is regular and $f\in V_n$. Hence the result. 
\end{proof}
In this next section, we will give examples in situations where the Shimura variety $X$ receives an immersive dense map from a moduli space of algebraic varieties, namely principally polarized abelian varieties and polarized K3 surfaces.

\subsubsection{Equidistribution in average of CM abelian varieties}

We now apply \Cref{equidistributionCM} to study equidistribution of CM principally polarized abelian varieties. Let $g\geq 1$, $\widetilde{G}=GSp_{2g}$ the standard symplectic group over $\Q$ and $\mathbb{H}_g$ the Siegel upper half-space. Then $(\widetilde{G},\mathbb{H}_g)$ is a Shimura datum and for $\Gamma=\mathrm{Sp}(2g,\Z)$, the quotient $\mathcal{A}_g\equaldef \Gamma\backslash \mathbb{H}_g$ is in bijection with the set of isomorphism classes of principally polarized abelian varieties over $\C$.
\medskip

For every $N\geq 1$, we have a double class $V_N=\{f\in \mathrm{GL}_{2g}(\Z)\cap \widetilde{G}(\Q), f^\dagger\circ f=N\cdot \mathrm{Id}\}$ where $f^\dagger$ is the adjoint with respect to the symplectic form. The Hecke correspondence $\mathcal{C}_N$ given by this double class has the following modular interpretation : $\mathcal{C}_N\hookrightarrow \mathcal{A}_g\times \mathcal{A}_g$ is the moduli of pairs $(A_1,A_2,f)$ where $(A_1,A_2)$ are principally polarized abelian varieties of dimension $g\geq 1$ and $f:A_1\rightarrow A_2$ is an isogeny satisfying $f^\dagger\circ f =N\cdot\mathrm{Id}_{A_1}$ where $f^\dagger:B\rightarrow A$ is the dual isogeny. Note that $\mathcal{C}_1$ is the diagonal embedding of $\cA_g$ in $\cA_g\times \cA_g$.
In this context, the transverse intersection locus of $\mathcal{C}_N$ with $\mathcal{C}_1$ corresponds to principally polarized CM abelian varieties $A$ endowed with an isogeny $f:A\rightarrow A$ whose homology class is a regular element of $\widetilde{G}(\R)$ and lies in $V_N $. Then the Mumford--Tate group of $A$ is a subgroup of $Z(f)$. 
\begin{lemma}
The transverse intersection loci of $\mathcal{C}_N$ and $\mathcal{C}_1$ is set theoretically formed by triples $(A,\lambda,f)$ where $(A,\lambda)$ is a principally polarized abelian variety, $f:A\rightarrow A$ is an isogeny whose homological realization is regular and lies in $V_N$ and $MT(A)\subseteq Z(f)$. In particular, $A$ is a CM abelian variety.
\end{lemma}
By applying \Cref{equidistributionCM}, we get \Cref{CMabelianvarieties}  in the introduction.

\subsubsection{Equidistribution in average of CM K3 surfaces}

We now discuss the second example which is the equidistribution of CM points in the moduli space of polarized K3 surfaces. Let $d\geq 1$ and let $\mathcal{F}_{2d}$ be the moduli space of complex polarized K3 surfaces of degree $2d$. Then $\mathcal{F}_{2d}$ can be embedded into an orthogonal Shimura variety which is given as follows. Let $V_{K3}$ be the K3 lattice, $V_{K3}=U^{3}\oplus E_{8}(-1)^2$, $\ell_{2d}\in V_{K3}$ a primitive class of self-intersection $2d$ (it is unique up to the action of $\mathrm{O}(V_{K3})$) and let $V_{2d}=\ell_{2d}^{\bot}$. Let $\widetilde{G}=\mathrm{GO}(V_{2d})$ and $\D=\{x\in \mathbb{P}(V_{2d,\C}, (x,x)=0, (x,\overline{x})=0)\}$. Then $(G,\D)$ is a Shimura datum and for $\Gamma=\mathrm{Ker}(\mathrm{O}(V_{2d})\rightarrow \mathrm{O}(V^{\vee}_{2d}/V_{2d}))$, we have a period map $\mathcal{F}_{2d}\rightarrow \Gamma\backslash \D$ which is a local embedding by Torelli theorem \cite{huybrechts} and the complement of the image is a finite union of Cartier divisors. Under this map, K3 surfaces with CM, in the sense of \cite[Remark 3.10]{huybrechts} correspond to CM points of the orthogonal Shimura variety $\Gamma\backslash \D$. Let $\omega_{\D}$ be the volume form on $\D$ as in \Cref{notations}, and for $N\geq 1$, let $V_N$ be the double class of integral elements $f\in \widetilde{G}(\Q) $ which scale the bilinear form by $N$. 
\begin{theorem}
Let $N\geq 1$ and let $CM(N)$ be the set of pairs $(X,\ell_{2d},f)$ where $(X,\ell_{2d})$ is a CM polarized K3 surface of degree $2d$, $f\in End(PH^2(X,\Z))$ an isogeny with $f^\dagger\circ f=N$ and $f\in\mathrm{GO}(\Q)$ is regular. Then for every open relatively compact subset with zero measure boundary $\Omega \subset \mathcal{F}_{2d}$, we have 
\[|\{(X,\ell_{2d},f)\in CM(N), (X,\ell_{2d})\in \Omega\}|\underset{N\rightarrow \infty}{\sim}\frac{\chi(\widehat{\D}).\deg(V_N)}{\Vol(\widehat{\D})}\int_{\Omega}\omega_{\D}.\]
\end{theorem}
\subsection{Equidistribution of Hecke translates of the Torelli locus}\label{hecketranslatestorelli}
We prove in this section \Cref{Torelliequidistribution} and  \Cref{Corollary:Torelli}. As the reader will notice, this is a statement about the dynamics of Hecke operators rather than the varieties themselves. 
\medskip 

Let $S$ and $D$ be two subvarieties of $\mathcal{A}_g$ of complimentary dimensions and let $d$ be the dimension of $S$. Let $\omega_{G/H}$ be the pull-push form on $\mathcal{A}_g\times \mathcal{A}_g$ as constructed in \Cref{diagonalrestriction} with respect to the groups $G=\mathrm{PGSp}_{2g}\times \mathrm{PGSp}_{2g}$ and $H=\mathrm{PGSp}_{2g}$ embedded diagonally. We have an inclusion $\iota:S\times D\hookrightarrow \mathcal{A}_{g}\times\mathcal{A}_{g}$.

Let as before $\cF^{1}\rightarrow \mathcal{A}_g$ be the Hodge bundle of the universal abelian scheme $\mathcal{A}_g$ and let $\omega$ be its first Chern form with respect to the Hodge metric. Finally, let $\omega_{S}$ and $\omega_{D}$ be its restriction to $S$ and $D$ respectively.
\begin{lemma}
We have $\iota^{*}\omega_{G/H}=\frac{1}{Vol(\widehat{\mathbb{H}}_g)}\omega_{S}^d\wedge\omega_{D}^{\frac{g(g+1)}{2}-d}$.
\end{lemma}
\begin{proof}
This is a consequence of \Cref{diagonalrestriction} as the only non-vanishing differential forms are the product of a form of degree $2d$ and a form of degree $g(g+1)-2d$, as the others vanish on $S\times D$, combined with the fact that $H^{2d}(\mathcal{A}_g,\R)$ is generated by $\omega^d$ for $d\leq 2$, see \cite[Prop. 2.2]{vandergeer}, and whose dual form (in the sense preceding \Cref{diagonalrestriction}) is simply $\frac{1}{Vol(\widehat{\mathbb{H}}_g)}\omega^{g(g+1)/2-d}$ since the volume form on $\mathcal{A}_g$ is $\omega^{g(g+1)/2}$, hence the result.
\end{proof}
It is well-know that $\omega$ is K\"ahler form on  $\mathcal{A}_g$ and hence the integration of $\omega_{S}^d$ and $\omega_{D}^{\frac{g(g+1)}{2}-1}$ define Lebesgue measures on $S$ and $D$ respectively, which are in fact finite by \cite[Main Theorem 3.1]{mumford}. Finally, for $(s,d)\in S\times D$, an isogeny $f:\mathcal{A}_{S,s}\rightarrow \mathcal{A}_{D,d}$ is said to be \emph{regular} if it does not admit first order deformation, or, equivalently, $S\times D$ intersects $\mathcal{C}_f$ transversely at $(s,d)$.
\begin{theorem}\label{actualtorelli}
For every open relatively compact subsets with zero measure boundary $\Omega \subset S$ and $\Omega'\subset D$, we have
\begin{align*}
    \left|\{(s,d,f)|\, (s,d)\in \Omega\times \Omega', f\in \mathrm{Isog}^N(\mathcal{A}_{S,s}, \mathcal{A}_{D,d}),f\,\textrm{is regular}\}\right|\\
    \underset{N\rightarrow \infty}{\sim}\deg(V_N)\cdot \frac{\chi(\widehat{\mathbb{H}}_g)}{\Vol(\widehat{\mathbb{H}}_g)}\int_\Omega\omega_{S}^d\int_{\Omega'}\omega_D^{\frac{g(g+1)}{2}-d}.
\end{align*}
In particular, the locus of points in $S$ isogenous to  a point in $D$ is analytically dense in $S$.
\end{theorem}

\begin{proof}[Proof of \Cref{Torelliequidistribution}]
For every $N\geq 1$, we have defined the Hecke correspondance $\mathcal{C}_N\hookrightarrow \mathcal{A}_g\times \mathcal{A}_g$ which parameterizes pairs of principally polarized abelian varieties together with a polarized isogeny with similitude factor equal to $N$. By the Lemma above, the restriction of the pull-push form $\omega_{G/H}$ is K\"ahler, hence the generic transversality assumption in \Cref{generictransversality} is satisfied and we are in the setting of \Cref{main}: the transverse intersection locus of $S\times D$ and $\mathcal{C}_N$ is equidistributed with respect to $\frac{\chi(\widehat{\mathbb H}_g)}{\Vol(\widehat{\mathbb H}_g)}\omega_{S}^d\wedge\omega^{g(g+1)/2-d}_{D}$ as $N\rightarrow \infty$. By the discussion preceding the theorem, the transverse locus is exactly given by regular isogenies. Hence by choosing subsets of the form $\Omega\times \Omega'$, we get the desired equidistribution result. One has also similar equidistribution results on $D$.
\end{proof}

Let $g\geq 2$ and let $\mathscr{M}_g$ be the coarse moduli space parameterizing smooth projective genus $g$ curves over $\C$. Recall that for any such curve $C$, one can associate its Jacobian $J(C)$, which is a principally polarized abelian variety of dimension $g$ over $\C$. This construction can be done in families so that we get a map, the \emph{Torelli map}, $\iota_g:\mathscr{M}_g\hookrightarrow\mathcal{A}_g$ between coarse moduli spaces. This map is injective by \cite{oortsteenbrink} and its image  is called the \emph{Torelli locus}. For $g=4$, the Torelli locus is a divisor in $\mathcal{A}_4$. Hence \Cref{Corollary:Torelli} follows by applying the previous theorem to $D=\mathscr{M}_4$.

\bibliographystyle{alpha}
\bibliography{bibliographie}
\end{document}